\documentclass[11pt,reqno]{amsart}

\usepackage{geometry}
\usepackage[utf8]{inputenc}
\usepackage{amstext,latexsym,amsbsy,amsmath,amssymb,amsthm,mathtools,relsize,geometry,enumerate}
\usepackage{hyperref}
\hypersetup{
	colorlinks   = true, 
	urlcolor     = blue, 
	linkcolor    = red, 
	citecolor   = red 
}
\usepackage{mathtools, enumerate,enumitem}

\usepackage{mathrsfs}
\usepackage{booktabs}
\usepackage{caption}
\usepackage{float}
\usepackage{placeins}

\usepackage{comment}
\usepackage{graphicx}
\usepackage{epstopdf}
\setkeys{Gin}{width=\linewidth,totalheight=\textheight,keepaspectratio}
\graphicspath{{./images/}}

\usepackage[normalem]{ulem}

\usepackage{multicol}
\usepackage{booktabs}
\usepackage{mathrsfs}
\usepackage{color}

\numberwithin{equation}{section}

\newtheorem{remark}{Remark}[section]
\newtheorem{theorem}{Theorem}[section]

\newtheorem{lemma}{{Lemma}}[section]

\allowdisplaybreaks[4]

\newcommand{\mE}{\mathbb{E}}

\begin{document}
	
	\title{Strong error estimates for finite element approximations of the stochastic Burgers--Huxley equation with multiplicative noise}

	\author{
		Tin Tran$^1$ and 
		Liet Vo$^2$ 
	}

	\thanks{$^1$ School of Mathematical and Statistical Sciences, The University of Texas Rio Grande Valley, Edinburg, TX 78539, U.S.A.  ({\tt tin.tran01@utrgv.edu}).}
	
	\thanks{$^2$ School of Mathematical and Statistical Sciences, The University of Texas Rio Grande Valley, Edinburg, TX 78539, U.S.A.  ({\tt liet.vo@utrgv.edu}). This author was partially supported by the NSF grant DMS-2530211. Corresponding authors: Tin Tran and Liet Vo}

	\begin{abstract}
		This paper is concerned with strong error estimates for fully discrete finite element approximations of the stochastic Burgers--Huxley equation with multiplicative noise. We use the implicit Euler--Maruyama scheme for time discretization and continuous piecewise linear finite elements for spatial discretization. Moment bounds and exponential integrability estimates are established and used with the discrete stochastic Gronwall inequality to derive error estimates in full expectation. Sub-second-moment error estimates are first obtained in the discrete maximum-in-time $L^2$-norm and the discrete $L^2$-in-time $H^1$-norm. Second-moment error estimates in the discrete maximum-in-time $L^2$-norm are then established using a higher-moment argument. The convergence rates are of order $k^{1/2}$ in time and $h$ in space, with the spatial rate being optimal in the $H^1$-norm. Finally, numerical experiments are provided to examine the convergence rates and the influence of multiplicative noise on the solution.
	\end{abstract}
	\maketitle

	
	\medskip
	

	\section{Introduction}\label{sec1}
	We consider the stochastic Burgers--Huxley equation with multiplicative It\^o noise:
	\begin{subequations}\label{eq:spde-system}
		\begin{align}
			du &= \big(\nu\,\partial_x^2u-\alpha\,u\,\partial_xu+\beta\,u(1-u)(u-\gamma)\big)\,dt+\sigma(u)\,dW(t),
			&& \text{in }(0,T)\times D \quad \text{a.s.}, \label{eq:spde-system-a}\\
			u(0) &= u_0,
			&& \text{in }D \quad \text{a.s.}, \label{eq:spde-system-b}
		\end{align}
	\end{subequations}
	where $D=(0,L)$ with $L>0$ and $T>0$ is a fixed final time. Here, $u=u(t,x,\omega)$ is the unknown real-valued random field, and $u_0$ is the initial condition. We impose periodic boundary conditions with period $L$, that is, $u(t,x+L)=u(t,x)$ almost surely for almost every $(t,x)\in(0,T)\times\mathbb R$. The constants $\nu>0$, $\alpha>0$, and $\beta>0$ denote the viscosity, nonlinear advection, and reaction coefficients, respectively, while $\gamma\in(0,1)$ is the unstable threshold between the two stable equilibria of the deterministic reaction dynamics. In \eqref{eq:spde-system-a}, $\{W(t)\}_{t\geq0}$ is a real-valued Wiener process, and $\sigma$ is the noise coefficient whose definition and assumptions are given in Section~\ref{subsec2.5}.

	This paper is devoted to the strong error analysis of finite element approximations for the stochastic Burgers--Huxley equation with multiplicative noise. Although several results on well-posedness and qualitative properties are available, the numerical analysis remains limited, particularly in the multiplicative-noise case. Our primary goal is to establish polynomial-order error estimates in full expectation for a fully discrete finite element method.
	
	The stochastic Burgers--Huxley equation can be viewed as a random force perturbation of the deterministic Burgers--Huxley equation, which combines nonlinear convection, diffusion, and cubic reaction and admits traveling-wave solutions \cite{wang1990solitary,mohan2021generalized}. The interpretation of the solution $u$ depends on the application. When $\beta=0$, the equation reduces to the Burgers equation, where $u$ represents a scalar velocity field in a simplified model of wave steepening, viscous dissipation, and turbulence \cite{da1994stochastic,bec2007burgers,chekhlov1995kolmogorov}. When $\alpha=0$, it reduces to the Huxley reaction--diffusion equation, which is associated with simplified models of nerve excitation \cite{mohan2026stochastic,stannat2013stability,geldhauser2024travelling}. The full equation therefore provides a prototype for studying the interaction between nonlinear transport and bistable reaction. In \eqref{eq:spde-system}, the random forcing is represented by $\sigma(u)\,dW(t)$, where $W$ is a real-valued Wiener process. When $\sigma$ is independent of $u$, the noise is called additive; otherwise, it is called multiplicative. The latter case is the focus of this paper.
	
	Several works have studied the mathematical properties of this equation. Mohan \cite{mohan2026stochastic} established global strong solvability under multiplicative Gaussian noise using local monotonicity and a localized Minty--Browder argument, and investigated large deviations and ergodicity under additive Gaussian noise. For the generalized equation driven by space--time white noise, Mohan \cite{mohan2022mild} proved existence and uniqueness of global mild solutions using fixed-point and stopping-time arguments. Kumar and Mohan \cite{kumar2024absolute} considered multiplicative space--time white noise and established global solvability, a comparison theorem, Malliavin differentiability, and, under suitable assumptions, absolute continuity of the solution law. On the numerical side, Wang, Cao, and Cao \cite{wang2025strong} established strong convergence rates for a spectral Galerkin method combined with a tamed exponential integrator under additive space--time white noise. To the best of our knowledge, strong error estimates in full expectation for finite element approximations of \eqref{eq:spde-system} with multiplicative noise have not previously been established.
	
	Although numerical methods for the deterministic equation can be adapted to construct stochastic approximations, their error analyses cannot be extended directly. There are several additional difficulties. First, the objective is to estimate errors in stochastic quantities, such as expectations and higher moments, rather than errors for a single realization. Second, the interaction between nonlinear convection, cubic reaction, and solution-dependent noise must be controlled. Finally, the standard deterministic Gronwall argument does not directly apply to error inequalities containing martingale terms and random coefficients.
	
	To explain the main difficulty, notice that the convection term $u\partial_xu$ gives rise to products involving the solution and the approximation error. After applying Sobolev embedding and Young's inequality, these products lead to terms of the form $\|u\|_{H^1}^2\|e\|_{L^2}^2$, where $e$ denotes an error function. In the deterministic case, the solution norm can serve as a coefficient in a Gronwall argument. In the stochastic case, both $u$ and $e$ are random fields, and the expectation $\mathbb E[\|u\|_{H^1}^2\|e\|_{L^2}^2]$ cannot generally be bounded by a deterministic constant times $\mathbb E[\|e\|_{L^2}^2]$. The multiplicative noise introduces further stochastic error terms involving $\sigma(u)-\sigma(u_h)$. These difficulties require estimates that control the nonlinear terms and preserve the martingale structure of the stochastic terms.
	
	One approach is to estimate the error on sub-sample spaces where the relevant solution norms remain bounded. This localization approach was used by Carelli and Prohl \cite{carelli2012rates} and in subsequent works \cite{breit2021convergence,vo2026mixed,thieu2026splitting,li2022finite} for stochastic Navier--Stokes equations, Keller--Segel systems, and wave equations with multiplicative noise. The resulting localized estimates yield convergence in probability, which is weaker than convergence in $L^p(\Omega)$. However, they do not alone provide error bounds in full expectation, which are useful for quantifying errors in computed statistical quantities. In \cite{bessaih2019strong,bessaih2022space}, estimates on the complementary sub-sample spaces were also obtained, leading to strong convergence in full expectation. For the multiplicative-noise cases considered there, the convergence rates are logarithmic and therefore slower than polynomial rates.
	
	More recently, polynomial-order strong error estimates in full expectation have been obtained for stochastic Navier--Stokes, Kuramoto--Sivashinsky, Benjamin--Bona--Mahony, and Keller--Segel equations with multiplicative noise \cite{feng2025full,nguyen2026fully,nguyen2026finite,vo2026optimal}. These analyses use a framework based on exponential integrability estimates, the discrete stochastic Gronwall inequality \cite{kruse2018discrete}, and higher-moment arguments. The martingale terms are treated directly, while exponential integrability controls the random coefficients in the error inequalities. In this paper, we develop this approach for the stochastic Burgers--Huxley equation, where the convection and reaction nonlinearities must be estimated together with the multiplicative noise.
	
	Specifically, we consider a fully discrete method combining continuous piecewise linear finite elements in space with an implicit Euler--Maruyama scheme in time. We establish moment bounds, spatial regularity, temporal H\"older continuity, and exponential integrability estimates for the exact solution, together with moment bounds for the numerical solution. These estimates provide the analytical foundation for the error analysis. We first apply the discrete stochastic Gronwall inequality to obtain sub-second-moment error estimates in the discrete maximum-in-time $L^2(D)$-norm and the discrete $L^2$-in-time $H^1(D)$-norm. We then use a higher-moment argument to establish second-moment error estimates in the discrete maximum-in-time $L^2(D)$-norm. The convergence rates are of order $k^{1/2}$ in time and $h$ in space, with the spatial rate being optimal in the $H^1(D)$-norm. Thus, the main contribution is a finite element error analysis with polynomial rates in full expectation for the multiplicative-noise problem. Numerical experiments examine these rates and illustrate the behavior of the expected solution and individual sample solutions.
	
	The remainder of this paper is organized as follows. In Section~\ref{sec2}, we introduce the notation, assumptions, and variational formulation, and establish stability and regularity estimates for the exact solution. In Section~\ref{sec3}, we present the fully discrete method and prove its stability and high-order moment bounds. Section~\ref{sec4}, which is the main section, is devoted to the strong error analysis. We first derive sub-second-moment error estimates and then establish second-moment estimates through a higher-moment argument. Section~\ref{sec5} presents numerical experiments, and Section~\ref{sec6} concludes the paper. The appendix provides auxiliary results and proofs of selected estimates for the continuous problem.
	
	\section{Preliminaries}\label{sec2}
	\subsection{Notations}
	Throughout the paper, we use standard notation for function spaces. For $1\leq p\leq\infty$ and any integer $m\geq0$, we denote the Lebesgue and Sobolev spaces of $L$-periodic functions by $L^p(D)$ and $H^m(D)$, respectively, with their usual norms $\|\cdot\|_{L^p}$ and $\|\cdot\|_{H^m}$. The inner product in $L^2(D)$ is denoted by $(\cdot,\cdot)$, and $\mathcal L(X)$ denotes the space of bounded linear operators from a normed space $X$ into itself.
	
	For $s\geq0$, we define the periodic Bessel potential space by
	\[
	\mathcal H^s(D):=\left\{u\in L^2(D):L\sum_{k\in\mathbb Z}\left(1+\left|\frac{2\pi k}{L}\right|^2\right)^s|\hat u(k)|^2<\infty\right\},
	\]
	with inner product
	\[
	(u,v)_{\mathcal H^s}:=L\sum_{k\in\mathbb Z}\left(1+\left|\frac{2\pi k}{L}\right|^2\right)^s \hat u(k)\overline{\hat v(k)}
	\]
	and norm $\|u\|_{\mathcal H^s}:=(u,u)_{\mathcal H^s}^{1/2}$. Here,
	\[
	\hat u(k):=\frac1L\int_0^Le^{-2\pi ikx/L}u(x)\,dx
	\]
	is the $k$th Fourier coefficient of $u$. For every integer $m\geq0$, $\mathcal H^m(D)$ coincides with the periodic Sobolev space $H^m(D)$, and their norms are equivalent \cite{adams2003sobolev}. We also recall the one-dimensional Sobolev embedding
	\begin{equation}\label{ass:sobolev-embedding}
		\|u\|_{L^\infty}\leq C_e\|u\|_{H^1},
		\qquad \forall u\in H^1(D),
	\end{equation}
	where $C_e>0$ depends only on $D$. Let $\mathcal{A}$ be the linear operator on $H:=L^2(D)$ defined by
	\[
	\mathcal{A}u:=-\partial_x^2u,\qquad D(\mathcal{A})=H^2(D).
	\]
	
	The Fourier modes
	\[
	e_k(x):=\frac1{\sqrt L}e^{2\pi ikx/L},\qquad k\in\mathbb Z,
	\]
	form an orthonormal basis of the complexification of $H$ and satisfy
	\begin{equation*}
		\mathcal{A}e_k=\lambda_ke_k,\qquad \lambda_k=\left(\frac{2\pi k}{L}\right)^2.
	\end{equation*}
	
	The operator $\mathcal{A}$ is nonnegative and self-adjoint. Since constants belong to its kernel, we introduce $\Lambda:=(\mathcal{I}+\mathcal{A})^{1/2}$, where $\mathcal{I}$ denotes the identity operator. The Bessel potential norm can then be written as
	\[
	\|u\|_{\mathcal H^s}=\|\Lambda^su\|_{L^2}=\|(\mathcal{I}+\mathcal{A})^{s/2}u\|_{L^2},
	\qquad s\geq0;
	\]
	see \cite[Section~1.2.3]{adams2012function}. Moreover, $-\nu \mathcal{A}$ generates the analytic semigroup $S(t):=e^{-\nu t\mathcal{A}}$ on $H$. We use the semigroup estimate
	\begin{equation}\label{eq:semigroup-increment}
		\|\Lambda^{-\rho}(\mathcal{I}-S(t))\|_{\mathcal L(H)}
		\leq Ct^{\rho/2},
		\qquad t>0,\quad \rho\in[0,2];
	\end{equation}
	see \cite[Lemma~3.22]{lord2014introduction}.
	
	We define the Huxley reaction term by $c(u):=u(1-u)(u-\gamma)$, which satisfies the one-sided Lipschitz condition
	\begin{equation}\label{eq:one-sided-lipschitz}
		\bigl(c(a)-c(b)\bigr)(a-b)\leq\mu|a-b|^2,
		\qquad \forall a,b\in\mathbb R.
	\end{equation}
	
	In addition, 
	we also obtain
	\begin{equation}\label{eq:reaction-growth}
		\|c(u)-c(v)\|_{L^2}
		\leq C\bigl(1+\|u\|_{L^\infty}^2+\|v\|_{L^\infty}^2\bigr)\|u-v\|_{L^2},
		\qquad \forall u,v\in L^\infty(D).
	\end{equation}
	
	By the Sobolev product estimates and \eqref{ass:sobolev-embedding}, for every integer $m\geq0$,
	\begin{equation}\label{eq:cubic-estimate}
		\begin{aligned}
			\|c(u)\|_{H^m}
			&\leq C\bigl(1+\|u\|_{L^\infty}^2\bigr)\|u\|_{H^m} \leq C\bigl(1+\|u\|_{H^{m+1}}^3\bigr),
			\qquad \forall u\in H^{m+1}(D).
		\end{aligned}
	\end{equation}
	
	Similarly, the Burgers nonlinearity satisfies
	\begin{equation}\label{eq:burger-estimate}
		\|u\partial_xu\|_{H^m}\leq C\|u\|_{H^{m+1}}^2,
		\qquad \forall u\in H^{m+1}(D).
	\end{equation}
	
	Here, the constants may depend on $m$, $\gamma$, and $D$, but are independent of $u$.
	
	Throughout the paper, $C>0$ denotes a generic constant that may vary from line to line but is independent of the spatial and temporal discretization parameters.
	
	\subsection{Assumptions}\label{subsec2.5}
	For any integer \(m\geq 0\), we assume that \(\sigma:H^m(D)\to H^m(D)\) satisfies the following conditions:
	There exists a constant $L_{\sigma} > 0$ such that, for all $u,v\in H^m(D)$,
	\begin{align}\label{ass:noise-lipschitz}
		\|\partial_x^m\sigma(u)-\partial_x^m\sigma(v)\|_{L^2}^2 \le L_{\sigma}\|\partial_x^m u-\partial_x^mv\|_{L^2}^2.
	\end{align}
	
	This implies that there exists a constant $K_{\sigma} > 0$ such that
	\begin{align}\label{ass:noise-linear-growth}
		\|\partial_x^m\sigma(u)\|_{L^2}^2 \le K_{\sigma}(1 + \bigl\|\partial_x^m u\|_{L^2}^2),
	\end{align}
	where $K_{\sigma} = 2\max\{L_{\sigma},\|\partial_x^m\sigma(0)\|_{L^2}^2\}$.

	\subsection{Variational solution}
	We first introduce the variational solution of \eqref{eq:spde-system} in the following theorem. The proof follows \cite[Theorem~3.5]{mohan2026stochastic}, with minor modifications for periodic boundary conditions.
	\begin{theorem}\label{theo:unique-sol}
		Let $u_0\in L^{2p}(\Omega;L^2(D))$, $p\geq2$, be $\mathcal F_0$-measurable, and assume that $\sigma$ satisfies \eqref{ass:noise-lipschitz}. Then, there exists a unique adapted variational solution
		\[
		u\in L^{2p}\big(\Omega;L^\infty(0,T;L^2(D))\big)\cap L^2\big(\Omega;L^2(0,T;H^1(D))\big)\cap L^4\big(\Omega;L^4(0,T;L^4(D))\big)
		\]
		to \eqref{eq:spde-system}, which has almost surely continuous paths in $L^2(D)$. Moreover, for all $t\in[0,T]$ and $\phi\in H^1(D)$, $u$ satisfies, $\mathbb P$-a.s.,
		\begin{equation}\label{eq:week-solution}
			\begin{aligned}
				(u(t),\phi)&=(u_0,\phi)-\nu\int_0^t(\partial_xu(s),\partial_x\phi)\,ds-\alpha\int_0^t(u(s)\partial_xu(s),\phi)\,ds\\
				&\quad+\beta\int_0^t(u(s)(1-u(s))(u(s)-\gamma),\phi)\,ds+\left(\int_0^t\sigma(u(s))\,dW(s),\phi\right).
			\end{aligned}
		\end{equation}
	\end{theorem}
	
	\begin{remark}
		In \cite[Theorem~3.5]{mohan2026stochastic}, the problem is considered on $(0,1)$ with homogeneous Dirichlet boundary conditions. The same argument applies here by replacing the Dirichlet spaces with periodic spaces and using the full $H^1(D)$-norm. The boundary terms from integration by parts vanish by periodicity. We omit the proof to save space. We impose periodic boundary conditions to derive the higher-order regularity estimates in Lemma~\ref{lem:Hm-q-moment}. For well-posedness, homogeneous Dirichlet boundary conditions can also be used, as in the cited reference.
	\end{remark}
	
	We first establish an exponential stability estimate for $u$, which plays an important role in the error analysis.
	\begin{lemma}\label{lem:exponential-L2-estimate}
		Let $u_0\in L^2(\Omega;L^2(D))$ and assume that $\sigma$ satisfies \eqref{ass:noise-linear-growth}. For $0<\kappa\leq\frac{\beta}{8K_\sigma L}$, suppose that $\mathbb E[\exp(2\kappa\|u_0\|_{L^2}^2)]<\infty$. Then, for all $t\in[0,T]$,
		\[
		\mathbb E\left[\exp\left(\kappa\|u(t)\|_{L^2}^2+2\beta\gamma\kappa\int_0^t\|u(s)\|_{L^2}^2\,ds+2\nu\kappa\int_0^t\|\partial_xu(s)\|_{L^2}^2\,ds\right)\right]
		\leq C_{\ref{lem:exponential-L2-estimate}},
		\]
		where
		\[
		C_{\ref{lem:exponential-L2-estimate}}
		:=\exp\left(\left[\frac{(\beta(1+\gamma)^2+K_\sigma)^2L}{\beta}+K_\sigma\right]T\kappa+2\kappa^2K_\sigma T\right)
		\left(\mathbb E\left[\exp\left(2\kappa\|u_0\|_{L^2}^2\right)\right]\right)^{1/2}.
		\]
	\end{lemma}
	\begin{proof}
		See Appendix~\ref{AppendixB1}.
	\end{proof}
	
	We next establish higher-order regularity and moment estimates for $u$ under suitable assumptions on the initial data.
	\begin{lemma}\label{lem:Hm-q-moment}
		Let $m\geq0$ and $q\geq2$ be integers, and assume that $\sigma$ satisfies \eqref{ass:noise-linear-growth}. Suppose that $u_0\in\bigcap_{i=0}^{m}L^{3^{m-i}q}(\Omega;H^i(D))$. Then,
		\begin{align}\label{eq:Hm-moment-estimate}
			\sup_{s\in[0,T]}\mathbb E\big[\|u(s)\|_{H^m}^q\big]
			+\frac{\nu q}{2}\int_0^T\mathbb E\big[\|u(s)\|_{H^m}^{q-2}\|\partial_xu(s)\|_{H^m}^2\big]\,ds
			\leq C_{\ref{lem:Hm-q-moment}}^{m,q}
		\end{align}
		and
		\begin{equation}\label{eq:Hm-sup-moment-estimate}
			\mathbb E\left[\sup_{s\in[0,T]}\|u(s)\|_{H^m}^q\right]
			\leq C_{\ref{lem:Hm-q-moment}}^{m,q},
		\end{equation}
		where $C_{\ref{lem:Hm-q-moment}}^{m,q}=C\left(1+\sum_{i=0}^{m}\mathbb E[\|u_0\|_{H^i}^{3^{m-i}q}]\right)$.
	\end{lemma}
	\begin{proof}
		See Appendix~\ref{AppendixB2}.
	\end{proof}
	
	Finally, we establish temporal H\"older estimates for $u$, which will be used in the error analysis.
	\begin{lemma}\label{lem:time-regularity}
		Let $m\geq0$ and $q\geq2$ be integers. Suppose that $u_0\in\bigcap_{i=0}^{m+1}L^{6q\cdot3^{m+1-i}}(\Omega;H^i(D))$ and that $\sigma$ satisfies \eqref{ass:noise-linear-growth}. Then, for all $0\leq s\leq t\leq T$,
		\begin{equation}\label{eq:H1-time-regularity}
			\mathbb E\left[\|u(t)-u(s)\|_{H^m}^{2q}\right]
			\leq C_{\ref{lem:time-regularity}}^{m,q}(t-s)^q,
		\end{equation}
		where $C_{\ref{lem:time-regularity}}^{m,q}=C(T,m,q,u_0,K_\sigma)>0$ is independent of $s$ and $t$.
	\end{lemma}
	\begin{proof}
		See Appendix~\ref{AppendixB3}.
	\end{proof}

	\section{Fully discrete finite element method}\label{sec3}
	\subsection{Formulation and stability of the fully discrete finite element method}
	Let $N$ and $M$ be positive integers. We partition $D=(0,L)$ into intervals with nodes $0=x_0<x_1<\cdots<x_N=L$ and assume that the mesh is quasi-uniform. We set $h:=\max_{1\leq i\leq N}(x_i-x_{i-1})$, $k:=T/M$, and $t_n:=nk$ for $n=0,\dots,M$. The finite element space is defined by
	\begin{equation}
		V_h:=\Bigl\{\varphi_h\in H^1(D):\varphi_h|_{(x_{i-1},x_i)}\in \mathcal P_1,\;1\leq i\leq N\Bigr\},
	\end{equation}
	where $\mathcal P_1$ denotes the space of polynomials of degree at most one. Here, periodicity is included in the definition of $H^1(D)$.
	
	We denote by $P_h:L^2(D)\to V_h$ the $L^2$-projection defined by
	\begin{equation}\label{eq:projector}
		(u-P_hu,\phi_h)=0,\qquad \forall\phi_h\in V_h.
	\end{equation}
	
	The following approximation estimate holds \cite{brenner2008mathematical}:
	\begin{equation}\label{eq:projector-inequality}
		\|u-P_hu\|_{L^2}+h\|\partial_x(u-P_hu)\|_{L^2}
		\leq C_Ph^2\|u\|_{H^2},\qquad \forall u\in H^2(D),
	\end{equation}
	where $C_P>0$ is independent of $h$.
	
	\medskip
	\noindent\textbf{Main Algorithm.} Set $u_h^0=P_hu_0$. For $n=0,1,\dots,M-1$, find $u_h^{n+1}\in V_h$ such that, $\mathbb P$-a.s.,
	\begin{equation}\label{scheme:implicit-euler}
		\begin{aligned}
			(u_h^{n+1}-u_h^n,\phi_h)
			&+\nu k(\partial_xu_h^{n+1},\partial_x\phi_h)
			+\alpha k(u_h^{n+1}\partial_xu_h^{n+1},\phi_h)\\
			&=\beta k(u_h^{n+1}(1-u_h^{n+1})(u_h^{n+1}-\gamma),\phi_h)
			+(\sigma(u_h^n)\Delta W_n,\phi_h),
		\end{aligned}
	\end{equation}
	for all $\phi_h\in V_h$, where $\Delta W_n:=W(t_{n+1})-W(t_n) \sim \mathcal{N}(0,k)$.
	\medskip
	
	We first establish the following stability estimates for the fully discrete solution.
	\begin{lemma}\label{lem:stability}
		Suppose that $u_0\in L^2(\Omega;L^2(D))$ and that $\sigma$ satisfies \eqref{ass:noise-linear-growth}. Then, there exists $k_0>0$ such that, for $0<k\leq k_0$, the solution $\{u_h^n\}_{n=0}^{M}$ generated by the Main Algorithm satisfies
		\[
		\mathbb E\left[\max_{1\leq n\leq M}\|u_h^n\|_{L^2}^2\right]
		+\sum_{n=0}^{M-1}\mathbb E\left[\|u_h^{n+1}-u_h^n\|_{L^2}^2\right]
		+\nu k\sum_{n=1}^{M}\mathbb E\left[\|\partial_xu_h^n\|_{L^2}^2\right]
		\leq C_{\ref{lem:stability}},
		\]
		where $C_{\ref{lem:stability}}=C(u_0,K_\sigma,T,\beta,\gamma)>0$ is independent of $h$ and $k$.
	\end{lemma}
	\begin{proof}
		Taking $\phi_h=u_h^{n+1}$ in \eqref{scheme:implicit-euler} and using the identity $2(x-y,x)=\|x\|_{L^2}^2-\|y\|_{L^2}^2+\|x-y\|_{L^2}^2$, we obtain
		\begin{align}
			&\frac12\left[\|u_h^{n+1}\|_{L^2}^2-\|u_h^n\|_{L^2}^2+\|u_h^{n+1}-u_h^n\|_{L^2}^2\right]
			+\nu k\|\partial_xu_h^{n+1}\|_{L^2}^2+\beta k\|u_h^{n+1}\|_{L^4}^4+\beta\gamma k\|u_h^{n+1}\|_{L^2}^2 \nonumber\\
			&=-\alpha k\big(u_h^{n+1}\partial_xu_h^{n+1},u_h^{n+1}\big)
			+\beta k(1+\gamma)\int_D\big(u_h^{n+1}\big)^3\,dx
			+\big(\sigma(u_h^n)\Delta W_n,u_h^{n+1}\big) \nonumber\\
			&=:I_1+I_2+I_3. \label{eq:disc-energy}
		\end{align}
		
		By integration by parts and periodicity, we have $I_1=0$. Using H\"older's and Young's inequalities, we obtain
		\begin{align*}
			I_2&:=\beta k(1+\gamma)\int_D\big(u_h^{n+1}\big)^3\,dx\\
			&\leq\beta k(1+\gamma)\|u_h^{n+1}\|_{L^4}^2\|u_h^{n+1}\|_{L^2}\\
			&\leq\beta k\left(\frac14\|u_h^{n+1}\|_{L^4}^4+(1+\gamma)^2\|u_h^{n+1}\|_{L^2}^2\right).
		\end{align*}
		
		Using the Cauchy--Schwarz inequality and assumption \eqref{ass:noise-linear-growth}, we have
		\begin{align*}
			I_3 &:= \big( \sigma(u_h^n)\,\Delta W_n, u_h^{n+1}-u_h^{n} \big) + \big( \sigma(u_h^n)\,\Delta W_n, u_h^{n} \big) \\
			&\le \tfrac14 \,\|u_h^{n+1}-u_h^{n}\|_{L^2}^{2} + \| \sigma(u_h^n)\,\Delta W_n \|_{L^2}^{2} + \big( \sigma(u_h^n)\,\Delta W_n,u_h^{n} \big) \\
			&\le \tfrac14 \,\|u_h^{n+1}-u_h^{n}\|_{L^2}^{2} + \| \sigma(u_h^n) \|_{L^2}^{2}\,|\Delta W_n|^{2} + \big( \sigma(u_h^n)\,\Delta W_n,u_h^{n} \big) \\
			&\le \tfrac14 \,\|u_h^{n+1}-u_h^{n}\|_{L^2}^{2} + K_{\sigma}\bigl(1 + \| u_h^{n} \|_{L^2}^{2}\bigr)\,|\Delta W_n|^{2} + \big( \sigma(u_h^n)\,\Delta W_n, u_h^{n} \big).
		\end{align*}
		
		Substituting the estimates for \(I_1,I_2\), and \(I_3\) into \eqref{eq:disc-energy}, we arrive at
		\begin{align*}
			\frac{1}{2} &\left[\|u^{n+1}_h\|^2_{L^2} - \|u^n_h\|^2_{L^2} + \frac{1}{2}\|u_h^{n+1} - u_h^n\|^2_{L^2}\right] + \nu k\|\partial_x u_h^{n+1}\|^2_{L^2}  + \frac{3\beta k}{4} \|u^{n+1}_h\|^4_{L^4} + \beta k \gamma \|u^{n+1}_h\|^2_{L^2} \nonumber \\
			&\le \beta\,k(1+\gamma)^2\,\|u_h^{n+1}\|_{L^2}^{2} + K_{\sigma}\bigl(1 + \| u_h^{n} \|_{L^2}^{2}\bigr)\,|\Delta W_n|^{2} +\bigl(\sigma(u_h^{n})\,\Delta W_n,u_h^{n}\bigr)\,.
		\end{align*}
		
		Summing the above inequality over \(n=0,\ldots,\ell\), where \(0\le \ell\le M-1\), and defining
		\[
		\mathcal{M}_\ell := \sum_{n=0}^{\ell} \bigl(\sigma(u_h^n)\Delta W_n,u_h^n\bigr),
		\]
		we obtain
		\begin{align*}
			\frac{1}{2} & \|u^{\ell+1}_h\|^2_{L^2} + \frac{1}{4}\sum^{\ell}_{n=0} \|u_h^{n+1} - u_h^n\|^2_{L^2} + \nu k \sum^{\ell}_{n=0} \|\partial_x u_h^{n+1}\|^2_{L^2}  + \frac{3\beta k}{4} \sum^{\ell}_{n=0}  \|u^{n+1}_h\|^4_{L^4} + \beta k \gamma \sum^{\ell}_{n=0}  \|u^{n+1}_h\|^2_{L^2} \nonumber \\
			&\le \frac{1}{2} \|u^{0}_h\|^2_{L^2} +\beta\,k(1+\gamma)^2 \sum^{\ell}_{n=0} \|u_h^{n+1}\|_{L^2}^{2} + K_{\sigma} \sum^{\ell}_{n=0} \bigl(1 + \| u_h^{n} \|_{L^2}^{2}\bigr)\,|\,\Delta W_n|^{2} + \mathcal{M}_\ell.
		\end{align*}
		
		Multiplying the above inequality by $4$ and using
		\[
		\sum_{n=0}^{\ell}\|u_h^{n+1}\|_{L^2}^2 = \sum_{n=0}^{\ell}\|u_h^{n}\|_{L^2}^2 -\|u_h^{0}\|_{L^2}^2 +\|u_h^{\ell+1}\|_{L^2}^2,
		\]
		we obtain
		\[\begin{aligned}
			&\Bigl(2-4\beta k(1+\gamma)^2\Bigr)\|u_h^{\ell+1}\|_{L^2}^2 +\sum_{n=0}^{\ell}\|u_h^{n+1}-u_h^n\|_{L^2}^2 +4\nu k\sum_{n=0}^{\ell}\|\partial_xu_h^{n+1}\|_{L^2}^2 \\
			&\le 2\|u_h^0\|_{L^2}^2 + 4\beta k(1+\gamma)^2\sum_{n=0}^{\ell}\|u_h^n\|_{L^2}^2 +4K_\sigma\sum_{n=0}^{\ell}(1+\|u_h^n\|_{L^2}^2)\,|\,\Delta W_n|^2 +4\mathcal{M}_\ell.
		\end{aligned}\]
		
		Choosing $k_0>0$ such that $4\beta k(1+\gamma)^2\le 1$, $\forall k\in(0,k_0]$, yields
		\begin{equation}\label{eq:L2-energy-pre-max}
			\begin{aligned}
				&\|u_h^{\ell+1}\|_{L^2}^2 +\sum_{n=0}^{\ell}\|u_h^{n+1}-u_h^n\|_{L^2}^2+4\nu k\sum_{n=0}^{\ell}\|\partial_xu_h^{n+1}\|_{L^2}^2 \\
				\le&\,2\|u_h^0\|_{L^2}^2 + 4\beta k(1+\gamma)^2\sum_{n=0}^{\ell}\|u_h^n\|_{L^2}^2 +4K_\sigma\sum_{n=0}^{\ell}(1+\|u_h^n\|_{L^2}^2)\,|\,\Delta W_n|^2 +4\mathcal{M}_\ell.
			\end{aligned}
		\end{equation}
		
		Next, applying the expectation to \eqref{eq:L2-energy-pre-max}, and using the facts that
		\begin{align*}
			\mathbb{E}\!\left[1 + \| u_h^{n} \|_{L^2}^{2}\,|\,\Delta W_n|^2\right] &= \mathbb{E}\!\left[1 + \|u_h^{n}\|_{L^2}^2\right]\, \mathbb{E}\!\left[|\Delta W_n|^2\right] = k\bigl(1 + \mathbb{E}\!\left[\|u_h^{n}\|_{L^2}^2\right]\bigr),
		\end{align*}
		and that $\mathbb{E}\!\left[(\,\sigma(u_h^{n})\,\Delta W_n,u_h^{n}\,)\right] = 0$. We obtain
		\begin{equation}\label{eq:lemma31-before-discrete-gronwall}
			\begin{aligned}
				&\mE\big[\|u^{\ell+1}_h\|^2_{L^2}\big] +\sum_{n=0}^{\ell}\mathbb{E}\big[\|u_h^{n+1}-u_h^n\|_{L^2}^2\big]+ 4\nu k \sum^{\ell}_{n=0} \mE\big[\|\partial_x u_h^{n+1}\|^2_{L^2}\big]  \\
				&\le 2\mE\big[\|u^{0}_h\|^2_{L^2}\big] +4\left(\beta(1+\gamma)^2 + K_{\sigma}\right) k\sum^{\ell}_{n=0} \bigl(1+\mE\big[\|u_h^{n}\|_{L^2}^{2}\big]\bigr).
			\end{aligned}
		\end{equation}
		
		By the discrete Gronwall inequality, we obtain
		\begin{equation}\label{eq:L2-energy-after-gronwall}
			\begin{aligned}
				\max_{1\leq n \leq M} \mathbb{E}\left[\|u_h^{n}\|_{L^2}^{2}\right]  +\sum_{n=0}^{\ell}\mathbb{E}\big[\|u_h^{n+1}-u_h^n\|_{L^2}^2\big]+ 4\nu k\sum_{n=1}^{M} \mathbb{E}\left[\|\partial_xu_h^{n}\|_{L^2}^{2}\right]  \le C.
			\end{aligned}
		\end{equation}
		
		Next, returning to \eqref{eq:L2-energy-pre-max}, taking the maximum over $0 \leq \ell \leq M-1$ and then the expectation on both sides of \eqref{eq:L2-energy-pre-max}, we obtain
		\begin{equation}\label{eq:L2-energy-max}
			\begin{aligned}
				&\mathbb{E} \left[\max_{0\leq \ell \leq M-1} \|u_h^{\ell+1}\|_{L^2}^2\right] +\sum_{n=0}^{M-1}\mathbb{E}\bigl[\|u_h^{n+1}-u_h^n\|_{L^2}^2\bigr] +4\nu k\sum_{n=0}^{M-1}\mathbb{E}\bigl[\|\partial_xu_h^{n+1}\|_{L^2}^2\bigr] \\ 
				&\le 2\mathbb{E}\bigl[\|u_h^0\|_{L^2}^2\bigr] + 4k\bigl(K_{\sigma} + \beta (1+\gamma)^2\bigr) \sum_{n=0}^{M-1}\mathbb{E}\bigl[\|u_h^n\|_{L^2}^2\bigr] +4\mathbb{E}\left[\max_{0 \leq \ell \leq M - 1}\sum_{n=0}^{\ell}\bigl(\sigma(u_h^n)\Delta W_n,u_h^n\bigr)\right].
			\end{aligned}
		\end{equation}
		
		We apply the Burkh\"older--Davis-Gundy and Young's inequalities as follows:
			\begin{align}
				\mathbb{E}\left[\max_{0\le \ell\le M-1} \left|\sum_{n=0}^{\ell} \bigl(\sigma(u_h^n)\Delta W_n,u_h^n\bigr)\right|  \right] &\le C\mathbb{E}\left[ \left( k\sum_{n=0}^{M-1} \|\sigma(u_h^n)\|_{L^2}^{2}\|u_h^n\|_{L^2}^{2} \right)^{\frac12} \right] \notag\\
				&\le C\mathbb{E}\left[ \max_{0\le n\le M-1}\|u_h^n\|_{L^2} \left( k\sum_{n=0}^{M-1} \|\sigma(u_h^n)\|_{L^2}^{2} \right)^{\frac12} \right] \notag\\
				&\le \frac18\mathbb{E}\left[\|u_h^0\|_{L^2}^{2}\right] + \frac18\mathbb{E}\left[\max_{1\le n\le M-1}\|u_h^n\|_{L^2}^{2}\right] \label{eq:martingale-estimate}\\
				&\quad+ C\mathbb{E}\left[ k\sum_{n=0}^{M-1}\|\sigma(u_h^n)\|_{L^2}^{2}\right] \notag\\
				&\le \frac18\mathbb{E}\left[\|u_h^0\|_{L^2}^{2}\right] + \frac18\mathbb{E}\left[\max_{1\le n\le M}\|u_h^n\|_{L^2}^{2} \right] \notag\\
				&\quad + C\mathbb{E}\left[ k\sum_{n=0}^{M} \left(1+\|u_h^n\|_{L^2}^{2}\right)\right]\notag.
			\end{align}
			
			Substituting this bound into \eqref{eq:L2-energy-max}, multiplying both sides by \(2\), we obtain, after simplification,
			\begin{align*}
				&\mathbb{E} \left[\max_{0\leq \ell \leq M-1} \|u_h^{\ell+1}\|_{L^2}^2\right] +2\sum_{n=0}^{M-1}\mathbb{E}\bigl[\|u_h^{n+1}-u_h^n\|_{L^2}^2\bigr] +8\nu k\sum_{n=0}^{M-1}\mathbb{E}\bigl[\|\partial_xu_h^{n+1}\|_{L^2}^2\bigr]\\
				&\leq 5\mathbb{E}\bigl[\|u_h^0\|_{L^2}^2\bigr] + 8k\bigl(K_{\sigma} + \beta (1+\gamma)^2\bigr) \sum_{n=0}^{M}\mathbb{E}\bigl[\|u_h^n\|_{L^2}^2\bigr] + C\mathbb{E}\left[ k\sum_{n=0}^{M} \left(1+\|u_h^n\|_{L^2}^{2}\right)\right]\\
				&\leq 5\mathbb{E}\bigl[\|u_h^0\|_{L^2}^2\bigr] + 8\bigl(K_{\sigma} + \beta (1+\gamma)^2\bigr) T\mathbb{E}\bigl[\|u_h^0\|_{L^2}^2\bigr] + CT \mathbb{E}\bigl[\|u^0_h\|^2_{L^2}\bigr] +CT \\
				&\leq C_{\ref{lem:stability}},
			\end{align*}
			where the second inequality follows from \eqref{eq:L2-energy-after-gronwall}. The proof is complete.
		\end{proof}

		We next establish high-moment stability estimates for the fully discrete solution.
		\begin{lemma}\label{lem:stability-2}
			Let $q\geq1$. Suppose that $u_0\in L^q(\Omega;L^2(D))$ and that $\sigma$ satisfies \eqref{ass:noise-linear-growth}. Then, there exists $k_0>0$ such that, for $0<k\leq k_0$, the solution $\{u_h^n\}_{n=0}^{M}$ generated by the Main Algorithm satisfies
			\[
			\mathbb E\left[\max_{1\leq n\leq M}\|u_h^n\|_{L^2}^q\right]\leq C_{\ref{lem:stability-2}}^{0,q},
			\]
			where $C_{\ref{lem:stability-2}}^{0,q}=C(q,u_0,K_\sigma,T,\beta,\gamma)>0$ is independent of $h$ and $k$.
		\end{lemma}
		\begin{proof}
			It suffices to prove the estimate for \(q\ge 4\). Indeed, once the case \(q=4\) is established, the case \(1\le q<4\) follows from Jensen's inequality.
			\begin{equation*}
				\mathbb{E}\left[\max_{1\le n\le M}\|u_h^n\|^q_{L^2}\right] \le \left(\mathbb{E}\left[\max_{1\le n\le M}\|u_h^n\|_{L^2}^{4}\right]\right)^{q/4}\le C.
			\end{equation*}
			
			Starting from \eqref{eq:L2-energy-pre-max} and dropping the nonnegative terms on the left-hand side, we obtain
			\begin{equation*}
				\|u_h^{\ell+1}\|_{L^2}^2 \le C\|u_h^0\|_{L^2}^2 + Ck\sum_{n=0}^{\ell}\|u_h^n\|_{L^2}^2 + C\sum_{n=0}^{\ell} (1+\|u_h^n\|_{L^2}^2)\,|\,\Delta W_n|^2 + C\mathcal{M}_\ell .
			\end{equation*}
			
			Taking the maximum over \(0\le \ell\le m-1\), where \(1\le m\le M\), gives
			\begin{equation*}
				\max_{1\le n\le m}\|u_h^n\|_{L^2}^2 \le C\|u_h^0\|_{L^2}^2 + Ck\sum_{n=0}^{m-1}\|u_h^n\|_{L^2}^2 + C\sum_{n=0}^{m-1}(1+\|u_h^n\|_{L^2}^2)\,|\,\Delta W_n|^2 + C\max_{0\le \ell\le m-1}|\mathcal{M}_\ell|.
			\end{equation*}
			
			Set $r := q/2 \geq 2$. Raising both sides to the power \(r\), taking expectations, and using \((a_1+\cdots+a_4)^r\le 4^{r-1}(a_1^r+\cdots+a_4^r)\), we get
			\[
			\begin{aligned}
				\mathbb{E}\left[\max_{1\le n\le m}\|u_h^n\|_{L^2}^{q}\right] &\le C\mathbb{E}\left[\|u_h^0\|_{L^2}^{q}\right]  + C\mathbb{E}\left[\left(k\sum_{n=0}^{m-1}\|u_h^n\|_{L^2}^2\right)^r\right] \\
				&\quad + C\mathbb{E}\left[\left(\sum_{n=0}^{m-1}(1+\|u_h^n\|_{L^2}^2)\,|\,\Delta W_n|^2\right)^r\right] + C\mathbb{E}\left[\max_{0\le \ell\le m-1}|\mathcal{M}_\ell|^r\right].
			\end{aligned}
			\]
			
			We now estimate the three terms on the right-hand side. First, by H\"older's inequality,
			\[
			\left(k\sum_{n=0}^{m-1}\|u_h^n\|_{L^2}^2\right)^r \le T^{r-1}k\sum_{n=0}^{m-1}\|u_h^n\|_{L^2}^{2r}.
			\]
			
			Since \(2r=q\), we have
			\[
			\mathbb{E}\left[\left(k\sum_{n=0}^{m-1}\|u_h^n\|_{L^2}^2\right)^r\right]\le Ck\sum_{n=0}^{m-1} \mathbb{E}\left[\|u_h^n\|_{L^2}^{q}\right].
			\]
			
			Next, using the independence of \(u_h^n\) and \(\Delta W_n\), together with $\mathbb{E}\left[|\Delta W_n|^{2r}\right]=C_r k^r$, we obtain
			\[\begin{aligned}
				\mathbb{E}\left[\left(\sum_{n=0}^{m-1}(1+\|u_h^n\|_{L^2}^2)\,|\,\Delta W_n|^2\right)^r\right] &\le m^{r-1}\sum_{n=0}^{m-1}\mathbb{E}\left[(1+\|u_h^n\|_{L^2}^2)^r\,|\,\Delta W_n|^{2r}\right] \\
				&\quad\le C k\sum_{n=0}^{m-1}\left(1+\mathbb{E}\left[\|u_h^n\|_{L^2}^{2r}\right]\right) \\
				&\quad\le C\left(1 + k\sum_{n=0}^{m-1}\mathbb{E}\left[\|u_h^n\|_{L^2}^{q}\right]\right).
			\end{aligned}\]
			
			Finally, by the discrete Burkh\"older--Davis--Gundy inequality and assumption \eqref{ass:noise-linear-growth},
			\begin{align*}
				\mathbb{E}\left[\max_{0\le \ell\le m-1}|\mathcal{M}_\ell|^r\right] &\le C\mathbb{E}\left[\left(k\sum_{n=0}^{m-1} \|\sigma(u_h^n)\|_{L^2}^2\|u_h^n\|_{L^2}^2 \right)^{r/2} \right] \\
				&\le C\mathbb{E}\left[\left( k\sum_{n=0}^{m-1}\big(1+\|u_h^n\|_{L^2}^2\big)\|u_h^n\|_{L^2}^2\right)^{r/2} \right] \\ 
				&\le C\mathbb{E}\left[\left(k\sum_{n=0}^{m-1}\big(1+\|u_h^n\|_{L^2}^4\big)\right)^{r/2}\right].
			\end{align*}
			
			Since \(r\ge2\), H\"older's inequality yields
			\[
			\left(k\sum_{n=0}^{m-1}\big(1+\|u_h^n\|_{L^2}^4\big)\right)^{r/2} \le C\left(1+k\sum_{n=0}^{m-1}\|u_h^n\|_{L^2}^{2r}\right).
			\]
			
			Combining the above estimates, we arrive at
			\begin{equation*}
				\mathbb{E}\left[\max_{1\le n\le m}\|u_h^n\|_{L^2}^{q}\right]\le C\left(1+\mathbb{E}\left[\|u_h^0\|_{L^2}^{q}\right]+k\sum_{n=0}^{m-1}\mathbb{E}\left[\|u_h^n\|_{L^2}^{q}\right]\right).
			\end{equation*}
			
			Since $
			\mathbb{E}\left[\|u_h^n\|_{L^2}^{q}\right] \le \mathbb{E}\left[ \max_{0\le j\le n}\|u_h^j\|_{L^2}^{q}\right]$,
			the discrete Gronwall inequality gives
			\[
			\mathbb{E}\left[\max_{1\le n\le M}\|u_h^n\|_{L^2}^{q}\right]\le C\left(1+\mathbb{E}\left[\|u_h^0\|_{L^2}^{q}\right]\right).
			\] The proof is complete.
		\end{proof}
		
		\section{Full expectation error estimates}\label{sec4}
		In this section, we establish strong error estimates for the Main Algorithm with multiplicative noise. We first obtain sub-second-moment error estimates and then use a higher-moment argument to derive second-moment estimates.
		\subsection{Sub-second-moment error estimates}\label{sec4.1}
		\begin{theorem}\label{thm:sub-second-moment-error}
			Let $q\in(0,1)$, and let $u$ be the variational solution of \eqref{eq:week-solution}. Suppose that $\{u_h^n\}_{n=0}^{M}$ is generated by the Main Algorithm with $u_h^0=P_hu_0$, $u_0\in\bigcap_{i=0}^{2}L^{12\cdot3^{2-i}}(\Omega;H^i(D))$, and $\sigma$ satisfies \eqref{ass:noise-lipschitz}. Assume further that $\mathbb E[\exp(2\kappa\|u_0\|_{L^2}^2)]<\infty$ and $K_\sigma\leq\frac{\beta}{8L\kappa}$, where
			\[
			\kappa=\frac{80\alpha^2C_e^2q(1+q)}{\nu(1-q)}\max\left\{\frac{1}{\beta\gamma},\frac{1}{\nu}\right\}.
			\]
			Then
			\[
			\left(\mathbb E\left[\max_{1\leq n\leq M}\|u(t_n)-u_h^n\|_{L^2}^{2q}\right]\right)^{\frac1{2q}}
			+\left(\mathbb E\left[\left(\nu k\sum_{n=1}^{M}\|\partial_x(u(t_n)-u_h^n)\|_{L^2}^2\right)^q\right]\right)^{\frac1{2q}}
			\leq\widehat C\left(k^{1/2}+h\right),
			\]
			where $\widehat C=C(q,u_0,T,L_\sigma)>0$ is independent of $h$ and $k$.
		\end{theorem}
		\begin{proof}
			Define $e^n := u(t_n) - u^n_h = \theta^n + \varepsilon^n$, where
			\begin{align*}
				\theta^n := u(t_n) - P_hu(t_n) \quad\text{and} \quad \varepsilon^n := P_h u(t_n) - u^n_h.
			\end{align*}
			
			Subtracting the fully discrete scheme \eqref{scheme:implicit-euler} from the variational formulation \eqref{eq:week-solution}, we obtain, for all \(\phi_h\in V_h\),
			\begin{equation}\label{eq:subtract-equations}
				\begin{aligned}
					&(e^{n+1} - e^n, \phi_h) + \nu k\big(\partial_x e^{n+1}, \partial_x \phi_h\big)\\
					&= \nu\int^{t_{n+1}}_{t_n} \big(\partial_x (u(t_{n+1}) - u(s)),\partial_x \phi_h \big)\,ds + \alpha\int^{t_{n+1}}_{t_n} \big(u(t_{n+1})\partial_x u(t_{n+1})-u(s)\partial_x u(s),\phi_h \big)\,ds \\
					&\qquad - \alpha k \big(u(t_{n+1}) \partial_x u(t_{n+1}) - u^{n+1}_h \partial_x u_h^{n+1},\phi_h\big) \\
					&\qquad- \beta \int^{t_{n+1}}_{t_n} \big(u(t_{n+1})(1-u(t_{n+1}))(u(t_{n+1}) - \gamma) - u(s)(1-u(s))(u(s) - \gamma),\phi_h\big)\,ds \\
					&\qquad+ \beta k \big(u(t_{n+1})(1-u(t_{n+1}))(u(t_{n+1}) - \gamma) - u^{n+1}_h(1-u^{n+1}_h)(u^{n+1}_h - \gamma),\phi_h \big) \\
					&\qquad+ \left( \int^{t_{n+1}}_{t_n} \big(\sigma(u(s)) - \sigma(u^n_h)\big)\,dW(s),\phi_h \right).
				\end{aligned}
			\end{equation}
			
			Using the $L^2$-projection orthogonality \eqref{eq:projector}, we can recast the left-hand side of \eqref{eq:subtract-equations} as follow
			\begin{align*}
				(\varepsilon^{n+1} - \varepsilon^n, \phi_h) + \nu k\left(\partial_x \varepsilon^{n+1}, \partial_x \phi_h\right) =  - \nu k\left(\partial_x \theta^{n+1}, \partial_x \phi_h\right).
			\end{align*}
			
			Substituting this into \eqref{eq:subtract-equations}, taking $\phi_h = \varepsilon^{n+1} \in V_h$ and using the identity $2a(a-b) = a^2 - b^2 + (a-b)^2$, we obtain        
			\begin{equation}\label{eq:Z-equations}
				\begin{aligned}
					&\frac{1}{2} \left(\|\varepsilon^{n+1}\|^2_{L^2} - \|\varepsilon^{n}\|^2_{L^2} + \|\varepsilon^{n+1} - \varepsilon^n\|^2_{L^2}\right) + \nu k \|\partial_x \varepsilon^{n+1} \|^2_{L^2} \\
					&= - \nu k\left(\partial_x \theta^{n+1}, \partial_x \varepsilon^{n+1}\right) +\nu\int^{t_{n+1}}_{t_n} \left(\partial_x (u(t_{n+1}) - u(s)),\partial_x \varepsilon^{n+1} \right)ds \\
					&\qquad- \alpha\int^{t_{n+1}}_{t_n} \left(u(s)\partial_x u(s)-u(t_{n+1})\partial_x u(t_{n+1}),\varepsilon^{n+1} \right)ds \\
					&\qquad- \alpha k \left(u(t_{n+1}) \partial_x u(t_{n+1}) - u^{n+1}_h \partial_x u_h^{n+1},\varepsilon^{n+1}\right) \\
					&\qquad- \beta \int^{t_{n+1}}_{t_n} \left(u(t_{n+1})(1-u(t_{n+1}))(u(t_{n+1}) - \gamma) - u(s)(1-u(s))(u(s) - \gamma),\varepsilon^{n+1}\right)ds \\ 
					&\qquad+ \beta k \left(u(t_{n+1})(1-u(t_{n+1}))(u(t_{n+1}) - \gamma) - u^{n+1}_h(1-u^{n+1}_h)(u^{n+1}_h - \gamma),\varepsilon^{n+1} \right) \\
					&\qquad+ \left(\int^{t_{n+1}}_{t_n} \sigma(u(s)) - \sigma(u^n_h))\,dW(s),\varepsilon^{n+1}\right) \\
					&=: Z_1 + Z_2 + Z_3 + Z_4 + Z_5 + Z_6 + Z_7.
				\end{aligned}
			\end{equation}
			
			We now estimate the individual terms on the right-hand side of \eqref{eq:Z-equations}. For $Z_1$, by applying the Cauchy--Schwarz inequality, Young's inequality, and \eqref{eq:projector-inequality}, we obtain 
			\begin{align*}
				Z_1 &\leq \nu k \|\partial_x \theta^{n+1}\|_{L^2}\|\partial_x \varepsilon^{n+1}\|_{L^2} \\ &\leq 2\nu k \|\partial_x \theta^{n+1}\|^2_{L^2} + \frac{\nu k}{8} \|\partial_x \varepsilon^{n+1}\|^2_{L^2} \\
				&\leq 2\nu k C_{P}^2h^{2}\|u(t_{n+1})\|^2_{H^2}+ \frac{\nu k}{8} \|\partial_x \varepsilon^{n+1}\|^2_{L^2}.
			\end{align*}
			
			For $Z_2$, also using the Cauchy--Schwarz inequality, we obtain
			\begin{align*}
				Z_2 &= \nu \int_{t_n}^{t_{n+1}} \left( \partial_x\bigl(u(t_{n+1})-u(s)\bigr), \partial_x \varepsilon^{n+1}\right)ds \\
				&\leq 2\nu \int^{t_{n+1}}_{t_n} \|\partial_x\bigl(u(t_{n+1})-u(s)\bigr)\|_{L^2}^2\,ds + \frac{\nu k}{8} \|\partial_x\varepsilon^{n+1}\|^2_{L^2}.
			\end{align*}
			
			Next, to estimate $Z_3$, we add and subtract the term  $u(t_{n+1})\partial_x u(s)$ and then use Sobolev embedding inequality \eqref{ass:sobolev-embedding} as follows
			\begin{align*}
				Z_3 &= -\alpha\int_{t_n}^{t_{n+1}} \Bigl(\big(u(s)-u(t_{n+1})\big)\partial_x u(s) + u(t_{n+1})\partial_x\big(u(s)-u(t_{n+1})\big), \varepsilon^{n+1} \Bigr)\,ds \\ 
				&\le \frac{\alpha^2}{4} \int_{t_n}^{t_{n+1}}\big\|\big(u(t_{n+1}) - u(s) \big)\partial_x u(s) + u(t_{n+1})\partial_x\big(u(t_{n+1}) - u(s))\big)\big\|_{L^2}^2\,ds + k\|\varepsilon^{n+1}\|_{L^2}^2 \\
				&\le \frac{\alpha^2}{2} \int_{t_n}^{t_{n+1}} \Bigl(\big\|\big(u(t_{n+1})-u(s)\big)\partial_x u(s)\big\|_{L^2}^2 + \big\|u(t_{n+1})\partial_x\big(u(t_{n+1})-u(s)\big)\big\|_{L^2}^2 \Bigr)\,ds + k\|\varepsilon^{n+1}\|_{L^2}^2 \\
				&\le \frac{\alpha^2}{2} \int_{t_n}^{t_{n+1}}\Bigl(\big\|u(t_{n+1})-u(s)\big\|_{L^\infty}^2\|\partial_x u(s)\|_{L^2}^2+ \|u(t_{n+1})\|_{L^\infty}^2\big\|\partial_x\big(u(t_{n+1})-u(s)\big)\big\|_{L^2}^2\Bigr)\,ds + k\|\varepsilon^{n+1}\|_{L^2}^2 \\
				&\le \frac{\alpha^2 C_e^2}{2} \int_{t_n}^{t_{n+1}}\Bigl(\|\partial_x u(s)\|_{L^2}^2 + \| u(t_{n+1})\|_{H^1}^2\Bigr)\big\|u(t_{n+1})-u(s)\big\|_{H^1}^2\,ds + k\|\varepsilon^{n+1}\|_{L^2}^2.
			\end{align*}
			
			We decompose and estimate \(Z_4\) as follows:
			\begin{align*}
				Z_4 &=-\alpha k \left(u(t_{n+1}) \partial_x u(t_{n+1}) - u^{n+1}_h \partial_x u_h^{n+1}, \varepsilon^{n+1}\right) \\
				&= -\alpha k \left((u(t_{n+1}) - u^{n+1}_h) \partial_x u(t_{n+1}) + u^{n+1}_h \partial_x (u(t_{n+1}) - u_h^{n+1}), \varepsilon^{n+1}\right) \\
				&= -\alpha k \left(e^{n+1}\partial_x u(t_{n+1}) + u^{n+1}_h \partial_x e^{n+1}, \varepsilon^{n+1}\right) \\
				&= -\alpha k \left(e^{n+1}\partial_x u(t_{n+1}), \varepsilon^{n+1}\right)  -\alpha k \left(u^{n+1}_h \partial_x e^{n+1}, \varepsilon^{n+1}\right) \\
				&=: Z_{4,1} + Z_{4,2}.
			\end{align*}
			
			Using the Sobolev embedding \eqref{ass:sobolev-embedding}, the $L^2$-projection estimate \eqref{eq:projector-inequality}, and Young's inequality, we estimate $Z_{4,1}$
			\begin{align*}
				Z_{4,1} &= -\alpha k \left(e^{n+1}\partial_x u(t_{n+1}), \varepsilon^{n+1}\right) \\
				&= -\alpha k \big(\varepsilon^{n+1}\partial_x u(t_{n+1}),\varepsilon^{n+1}\big) - \alpha k \big(\theta^{n+1}\partial_x u(t_{n+1}),\varepsilon^{n+1}\big)\\
				&\leq \alpha k C_e \big\|\varepsilon^{n+1}\big\|_{H^1} \|\partial_xu(t_{n+1})\|_{L^2} \big\|\varepsilon^{n+1}\big\|_{L^2} + \alpha k C_e C_P h^2 \|u(t_{n+1})\|_{H^2} \|\partial_xu(t_{n+1})\|_{L^2} \big\|\varepsilon^{n+1}\big\|_{H^1} \\
				&\leq \frac{\nu k}{8}\big\|\varepsilon^{n+1}\big\|_{H^1}^2 +  Ckh^4 \|u(t_{n+1})\|_{H^2}^4 + \frac{4\alpha^2C^2_ek}{\nu}\|\partial_xu(t_{n+1})\|_{L^2}^2 \big\|\varepsilon^{n+1}\big\|_{L^2}^2 \\
				&\leq \frac{\nu k}{8}\big\|\partial_x\varepsilon^{n+1}\big\|^2_{L^2} + \frac{\nu k}{8}\big\| \varepsilon^{n+1}\big\|^2_{L^2} + Ckh^4\|u(t_{n+1})\|^4_{H^2} \\
				&\quad\, + \frac{8\alpha^2C^2_e}{\nu}\Bigg(\int^{t_{n+1}}_{t_n} \big\|\partial_x \big(u(t_{n+1}) - u(s)\big)\big\|_{L^2}^2\,ds+\int^{t_{n+1}}_{t_n} \|\partial_x u(s)\|^2_{L^2}\,ds\Bigg) \big\|\varepsilon^{n+1}\big\|^2_{L^2}.
			\end{align*}
			
			With regard to $Z_{4,2}$, noting that
			\[
			(e^{n+1}\partial_xe^{n+1},e^{n+1})=0,
			\]
			which implies that
			\[
			(u(t_{n+1})\partial_xe^{n+1},e^{n+1}) = (u^{n+1}_h\partial_xe^{n+1},e^{n+1}).
			\]
			
			Thus, $Z_{4,2}$ can be rewritten as
			\begin{align*}
				Z_{4,2} &= -\alpha k\left(u(t_{n+1})\partial_x e^{n+1},e^{n+1}\right) + \alpha k\left(u_h^{n+1}\partial_x e^{n+1},\theta^{n+1}\right) \\
				&\leq \alpha k \|u(t_{n+1})\|_{L^\infty} \|\partial_x e^{n+1}\|_{L^2}\|e^{n+1}\|_{L^2} + \alpha k \|u_h^{n+1}\|_{L^2} \|\partial_x e^{n+1}\|_{L^2} \|\theta^{n+1}\|_{L^\infty}\\
				&\leq \alpha \Bigg(\Bigg\|\int^{t_{n+1}}_{t_n} \big(u(t_{n+1}) - u(s)\big)\,ds\Bigg\|_{L^{\infty}} + \Bigg\|\int^{t_{n+1}}_{t_n}  u(s)\,ds\Bigg\|_{L^{\infty}} \Bigg) \|\partial_x e^{n+1}\|_{L^2}\|e^{n+1}\|_{L^2} \\
				&\qquad\, +\alpha k C_e C_P h \|u_h^{n+1}\|_{L^2} \|\partial_x e^{n+1}\|_{L^2} \|u(t_{n+1})\|_{H^2} \\
				&\leq \alpha C_e \Bigg(\int^{t_{n+1}}_{t_n} \big\|u(t_{n+1}) - u(s)\big\|_{H^1}\,ds + \int^{t_{n+1}}_{t_n}  \big\|u(s)\big\|_{H^1}\,ds \Bigg) \|\partial_x e^{n+1}\|_{L^2}\|e^{n+1}\|_{L^2} \\
				&\qquad\, + \frac{\nu k}{32} \|\partial_x e^{n+1}\|_{L^2}^2 + Ckh^2 \|u_h^{n+1}\|_{L^2}^2 \|u(t_{n+1})\|_{H^2}^2 \\
				&\leq \frac{\nu k}{16} \|\partial_x e^{n+1}\|^2_{L^2} + \frac{16 \alpha^2 C_e^2}{\nu} \Bigg(\int^{t_{n+1}}_{t_n} \big\|u(t_{n+1}) - u(s)\big\|_{H^1}^2\,ds+\int^{t_{n+1}}_{t_n} \|u(s)\|^2_{H^1}\,ds\Bigg)\\
				&\qquad\, \times \|e^{n+1}\|^2_{L^2} + Ckh^2 \|u_h^{n+1}\|_{L^2}^2 \|u(t_{n+1})\|_{H^2}^2 \\
				&\leq \frac{\nu k}{8} \|\partial_x\varepsilon^{n+1}\|_{L^2}^2 + \frac{32\alpha^2C_e^2}{\nu}\left(\int_{t_n}^{t_{n+1}}\|u(t_{n+1})-u(s)\|_{H^1}^2\,ds + \int_{t_n}^{t_{n+1}}\|u(s)\|_{H^1}^2\,ds\right)\\
				&\qquad\,\times \|\varepsilon^{n+1}\|_{L^2}^2  + Ch^4\left(\int_{t_n}^{t_{n+1}}\|u(t_{n+1})-u(s)\|_{H^1}^2\,ds + \int_{t_n}^{t_{n+1}}\|u(s)\|_{H^1}^2\,ds\right) \|u(t_{n+1})\|_{H^2}^2 \\
				&\qquad\, + Ckh^2\big(1 +\|u^{n+1}_h\|^2_{L^2}\big)\|u(t_{n+1})\|^2_{H^2}.
			\end{align*}
			
			Using \eqref{eq:reaction-growth}, we estimate $Z_5$ as follows:
			\[
			\begin{aligned}
				Z_5 &= - \beta \int^{t_{n+1}}_{t_n} \Big(u(t_{n+1})(1-u(t_{n+1}))(u(t_{n+1}) - \gamma) - u(s)(1-u(s))(u(s) - \gamma),\varepsilon^{n+1}\Big)\,ds \\  
				&\leq \beta\left(\frac{5}{2} + \gamma\right)\int_{t_n}^{t_{n+1}} \left(1+\|u(t_{n+1})\|_{L^\infty}^2+\|u(s)\|_{L^\infty}^2 \right)\|u(t_{n+1})-u(s)\|_{L^2} \|\varepsilon^{n+1}\|_{L^2}\,ds \\
				&\leq C\int_{t_n}^{t_{n+1}} \left( 1+\|u(t_{n+1})\|_{L^\infty}^2+\|u(s)\|_{L^\infty}^2 \right)^2 \|u(t_{n+1})-u(s)\|_{L^2}^2\,ds + \frac{k}{2}\|\varepsilon^{n+1}\|_{L^2}^2 \\
				&\leq C\int_{t_n}^{t_{n+1}} \left( 1+\|u(t_{n+1})\|_{L^\infty}^4+\|u(s)\|_{L^\infty}^4 \right) \|u(t_{n+1})-u(s)\|_{L^2}^2\,ds + \frac{k}{2}\|\varepsilon^{n+1}\|_{L^2}^2 \\
				&\le C \int_{t_n}^{t_{n+1}}\left(1+\|u(t_{n+1})\|_{H^1}^{4} +\|u(s)\|_{H^1}^{4}\right) \|u(t_{n+1})-u(s)\|_{L^2}^{2}\,ds + \frac{k}{2}\|\varepsilon^{n+1}\|_{L^2}^{2}.
			\end{aligned}
			\]
			
			Using \eqref{eq:one-sided-lipschitz}, \eqref{eq:reaction-growth}, and the decomposition \(\varepsilon^{n+1}=e^{n+1}-\theta^{n+1}\), we expand the integrand in $Z_6$ as follows:
			\begin{equation*}
				\begin{aligned}
					Z_6 &= \beta k \left(u(t_{n+1})(1-u(t_{n+1}))(u(t_{n+1}) - \gamma) - u^{n+1}_h(1-u^{n+1}_h)(u^{n+1}_h - \gamma),\varepsilon^{n+1} \right) \\
					&\leq \beta \mu k\|e^{n+1}\|_{L^2}^2 + Ck \left(1+\|u(t_{n+1})\|_{L^\infty}^2+\|u_h^{n+1}\|_{L^\infty}^2\right)\|e^{n+1}\|_{L^2}\|\theta^{n+1}\|_{L^2} \\
					&\leq \Bigg(\beta \mu + \frac{1}{2}\Bigg)k\|e^{n+1}\|_{L^2}^2 + Ck \left(1+\|u(t_{n+1})\|_{L^\infty}^2 + \|u_h^{n+1}\|_{L^\infty}^2 \right)^2 \|\theta^{n+1}\|_{L^2}^2 \\
					&\leq \Bigg(\beta \mu + \frac{1}{2}\Bigg)k\|e^{n+1}\|_{L^2}^2 + Ckh^4 \left(1+\|u(t_{n+1})\|_{L^\infty}^2 +  \|u_h^{n+1}\|_{L^\infty}^2 \right)^2 \|u(t_{n+1})\|_{H^2}^2 \\ 
					&\leq \Bigg(\beta \mu + \frac{1}{2}\Bigg)k\|e^{n+1}\|_{L^2}^2 + Ckh^4 \left(1+\|u(t_{n+1})\|_{L^\infty}^4 + \|u_h^{n+1}\|_{L^\infty}^4 \right) \|u(t_{n+1})\|_{H^2}^2. 
				\end{aligned}
			\end{equation*}
			
			Using the one-dimensional inverse inequality $\|u_h^{n+1}\|_{L^\infty} \leq Ch^{-1/2}\|u_h^{n+1}\|_{L^2}$, we obtain
			\begin{equation*}
				\begin{aligned}
					Z_6 &\leq (2\beta \mu + 1)k\|\varepsilon^{n+1}\|_{L^2}^2 + Ckh^4\|u(t_{n+1})\|_{H^2}^2 + Ckh^4 \left(1+\|u(t_{n+1})\|_{H^1}^4 \right) \|u(t_{n+1})\|_{H^2}^2 \\
					&\qquad + Ckh^2 \|u_h^{n+1}\|_{L^2}^4 \|u(t_{n+1})\|_{H^2}^2 \\
					&\leq (2\beta \mu + 1)k\|\varepsilon^{n+1}\|_{L^2}^2 + Ckh^4 \left( 1+\|u(t_{n+1})\|_{H^2}^6 \right) + Ckh^2 \|u_h^{n+1}\|_{L^2}^4 \|u(t_{n+1})\|_{H^2}^2 .
				\end{aligned}
			\end{equation*}
			
			Now, we estimate the noise term $Z_7$ as follows:
			\begin{align*}
				Z_7 &=\left(\int_{t_n}^{t_{n+1}} (\sigma(u(s))-\sigma(u_h^n))\,dW(s), \,\varepsilon^{n+1}-\varepsilon^n\right) +\left(\int_{t_n}^{t_{n+1}} (\sigma(u(s))-\sigma(u_h^n))\,dW(s), \,\varepsilon^n\right) \\
				&=\left(\int_{t_n}^{t_{n+1}} (\sigma(u(s))-\sigma(u(t_n)))\,dW(s), \,\varepsilon^{n+1}-\varepsilon^n\right) +\left((\sigma(u(t_n))-\sigma(u_h^n))\Delta W_n, \,\varepsilon^{n+1}-\varepsilon^n\right) \\
				&\qquad +\left(\int_{t_n}^{t_{n+1}} (\sigma(u(s))-\sigma(u_h^n))\,dW(s), \,\varepsilon^n\right) \\
				&\le 2\left\|\int_{t_n}^{t_{n+1}}(\sigma(u(s))-\sigma(u(t_n)))\,dW(s)\right\|_{L^2}^2 +2\|(\sigma(u(t_n))-\sigma(u_h^n))\Delta W_n\|_{L^2}^2 \\
				&\qquad +\frac14\|\varepsilon^{n+1}-\varepsilon^n\|_{L^2}^2 +\left(\int_{t_n}^{t_{n+1}} (\sigma(u(s))-\sigma(u_h^n))\,dW(s), \,\varepsilon^n\right).
			\end{align*}
			
			In preparation for applying the stochastic Gronwall inequality in Lemma~\ref{lem:discrete-stochastic-gronwall}, we add and subtract \(2k\|\sigma(u(t_n))-\sigma(u_h^n)\|_{L^2}^2\) on the right-hand side of the estimate for \(Z_7\). Using the Lipschitz condition \eqref{ass:noise-lipschitz}, we then derive
			\begin{align*}
				Z_7 &\le 2\left\|\int_{t_n}^{t_{n+1}}(\sigma(u(s))-\sigma(u(t_n)))\,dW(s)\right\|_{L^2}^2 +\frac14\|\varepsilon^{n+1}-\varepsilon^n\|_{L^2}^2 +2\|(\sigma(u(t_n))-\sigma(u_h^n))\Delta W_n\|_{L^2}^2\\
				&\quad -2k\|\sigma(u(t_n))-\sigma(u_h^n)\|_{L^2}^2 +2k\|\sigma(u(t_n))-\sigma(u_h^n)\|_{L^2}^2 +\left(\int_{t_n}^{t_{n+1}} (\sigma(u(s))-\sigma(u_h^n))\,dW(s), \,\varepsilon^n\right) \\
				&\le 2\left\|\int_{t_n}^{t_{n+1}}(\sigma(u(s))-\sigma(u(t_n)))\,dW(s)\right\|_{L^2}^2
				+\frac14\|\varepsilon^{n+1}-\varepsilon^n\|_{L^2}^2 +2\|(\sigma(u(t_n))-\sigma(u_h^n))\Delta W_n\|_{L^2}^2 \\
				&\quad -2k\|\sigma(u(t_n))-\sigma(u_h^n)\|_{L^2}^2  +2L_{\sigma}k\|e^n\|_{L^2}^2 +\left(\int_{t_n}^{t_{n+1}} (\sigma(u(s))-\sigma(u_h^n))\,dW(s), \,\varepsilon^n\right) \\
				&\le 2\left\|\int_{t_n}^{t_{n+1}}(\sigma(u(s))-\sigma(u(t_n)))\,dW(s)\right\|_{L^2}^2 +\frac14\|\varepsilon^{n+1}-\varepsilon^n\|_{L^2}^2 +2\|(\sigma(u(t_n))-\sigma(u_h^n))\Delta W_n\|_{L^2}^2\\
				&\quad-2k\|\sigma(u(t_n))-\sigma(u_h^n)\|_{L^2}^2 +2L_{\sigma}k\|\varepsilon^n\|_{L^2}^2+2L_{\sigma}C_{P}^2h^{2}k\|u(t_n)\|_{H^2}^2 \\
				&\quad +\left(\int_{t_n}^{t_{n+1}} (\sigma(u(s))-\sigma(u_h^n))\,dW(s), \,\varepsilon^n\right).
			\end{align*}
			
			After substituting the bounds for $Z_1, Z_2,\dots,Z_7$ into \eqref{eq:Z-equations} and absorbing the similar terms into the left-hand side, we obtain
			\begin{align}\label{eqn:3.27}
				&\frac12\big(\|\varepsilon^{n+1}\|_{L^2}^2-\|\varepsilon^n\|_{L^2}^2\big) + \frac{1}{4}\|\varepsilon^{n+1}-\varepsilon^n\|_{L^2}^2 + \frac{\nu k}{2}\|\partial_x\varepsilon^{n+1}\|_{L^2}^2 \\\nonumber
				&\leq \frac{40 \alpha^2 C^2_e}{\nu}\Bigg(\int_{t_n}^{t_{n+1}}\|u(t_{n+1})-u(s)\|_{H^1}^2\,ds + \int_{t_n}^{t_{n+1}} \|u(s)\|_{H^1}^2\,ds \Bigg)\|\varepsilon^{n+1}\|_{L^2}^2 \\\nonumber
				&\quad + k\Bigg(\frac{5}{2} + \frac{\nu}{8} + 2\beta \mu \Bigg)\|\varepsilon^{n+1}\|_{L^2}^2  + 2L_{\sigma}k\|\varepsilon^n\|_{L^2}^2 \\\nonumber
				&\quad+ Ckh^2 \left(1+\|u_h^{n+1}\|_{L^2}^4\right) \left(\|u(t_{n+1})\|_{H^2}^2 + \|u(t_n)\|_{H^2}^2 \right) + Ckh^4 \left( 1+\|u(t_{n+1})\|_{H^2}^6\right) \\\nonumber
				&\quad + C\int_{t_n}^{t_{n+1}} \left( 1+\|\partial_xu(s)\|_{L^2}^2 +\|u(t_{n+1})\|_{H^1}^2 \right) \|u(t_{n+1})-u(s)\|_{H^1}^2\,ds \\\nonumber
				&\quad+ C\int_{t_n}^{t_{n+1}} \left( 1+\|u(t_{n+1})\|_{H^1}^4 +\|u(s)\|_{H^1}^4\right) \|u(t_{n+1})-u(s)\|_{L^2}^2\,ds \\\nonumber
				&\quad +Ch^4\left(\int_{t_n}^{t_{n+1}}\|u(t_{n+1})-u(s)\|_{H^1}^2\,ds + \int_{t_n}^{t_{n+1}}\|u(s)\|_{H^1}^2\,ds\right) \|u(t_{n+1})\|_{H^2}^2\\\nonumber
				&\quad+ 2\left\| \int_{t_n}^{t_{n+1}}\bigl(\sigma(u(s))-\sigma(u(t_n))\bigr)\,dW(s) \right\|_{L^2}^2 \\\nonumber
				&\quad + 2\|(\sigma(u(t_n))-\sigma(u_h^n))\Delta W_n\|_{L^2}^2 - 2k\|\sigma(u(t_n))-\sigma(u_h^n)\|_{L^2}^2 \\\nonumber
				&\quad+\left(\int_{t_n}^{t_{n+1}}\bigl(\sigma(u(s))-\sigma(u_h^n)\bigr)\,dW(s), \,\varepsilon^n\right).
			\end{align}
			
			Taking the summation over \(n=0,\ldots,\ell\), for any \(0\leq \ell<M\), we obtain
			\begin{align}\label{eq:theo-before-gronwall}
				\|\varepsilon^{\ell+1}\|^2_{L^2} + \frac{1}{2} \sum^{\ell}_{n=0} \|\varepsilon^{n+1} - \varepsilon^n\|^2_{L^2} + \nu k \sum^{\ell}_{n=0} \|\partial_x \varepsilon^{n+1} \|^2_{L^2} \leq F_{\ell} + M_{\ell} + \sum^{\ell-1}_{n=0} G_n \|\varepsilon^{n+1}\|^2_{L^2},
			\end{align}
			where
			\begin{align*}
				F_{\ell} &:= \sum^{\ell}_{n=0} \Bigg[ Ckh^2 \left(1+\|u_h^{n+1}\|_{L^2}^4\right) \left(\|u(t_{n+1})\|_{H^2}^2 + \|u(t_n)\|_{H^2}^2 \right) + Ckh^4 \left( 1+\|u(t_{n+1})\|_{H^2}^6\right) \\
				&\quad + C\int_{t_n}^{t_{n+1}} \big(1 + \|\partial_xu(s)\|_{L^2}^2+\|u(t_{n+1})\|_{H^1}^2 \big) \|u(t_{n+1})-u(s)\|_{H^1}^2\,ds \\
				&\quad + C\int_{t_n}^{t_{n+1}} \left(1+\|u(t_{n+1})\|_{H^1}^4+\|u(s)\|_{H^1}^4 \right) \|u(t_{n+1})-u(s)\|_{L^2}^2\,ds \\
				&\quad +Ch^4\left(\int_{t_n}^{t_{n+1}}\|u(t_{n+1})-u(s)\|_{H^1}^2\,ds + \int_{t_n}^{t_{n+1}}\|u(s)\|_{H^1}^2\,ds\right) \|u(t_{n+1})\|_{H^2}^2\\
				&\quad + 4\left\|\int_{t_n}^{t_{n+1}} \bigl(\sigma(u(s))-\sigma(u(t_n))\bigr)\,dW(s) \right\|_{L^2}^2 \\
				&\quad + \frac{80\alpha^2C^2_e}{\nu}\Bigg(\int_{t_n}^{t_{n+1}}\|u(t_{n+1})-u(s)\|_{H^1}^2\,ds \Bigg)\|\varepsilon^{n+1}\|^2_{L^2} \Bigg]\\
				&\quad +  C\Bigg(k + \int_{t_\ell}^{t_{\ell+1}} \|u(s)\|_{H^1}^2\,ds \Bigg)\|\varepsilon^{\ell+1}\|_{L^2}^2,\\
				M_{\ell} &:= \sum^{\ell}_{n=0} Z_n,\\
				Z_{n} &:= 4\|(\sigma(u(t_n))-\sigma(u_h^n))\Delta W_n\|_{L^2}^2 -4k\|\sigma(u(t_n))-\sigma(u_h^n)\|_{L^2}^2 \\        
				&\quad\,+2\left(\int_{t_n}^{t_{n+1}} (\sigma(u(s))-\sigma(u_h^n))\,dW(s),\varepsilon^n\right),\\
				G_{n} &:= \frac{80\alpha^2C^2_e}{\nu} \int_{t_n}^{t_{n+1}} \|u(s)\|_{H^1}^2\,ds + k\Bigg(5 + \frac{\nu}{4} + 4\beta \mu + 4L_{\sigma}\Bigg).
			\end{align*}
			
			Assume for the moment that \(\{M_\ell\}_{\ell\geq 0}\) is a martingale; this will be verified at the end of the proof. Applying the stochastic Gronwall inequality in Lemma \ref{lem:discrete-stochastic-gronwall} to \eqref{eq:theo-before-gronwall}, for any \(q\in(0,1)\) with \(a=\frac{1+q}{2q}\), and \(b=\frac{1+q}{1-q}\), we obtain
			
			\begin{equation}\label{eq:discrete-gronwall-bound}
				\begin{aligned}
					&\left(\mathbb E\left[\sup_{0\leq \ell \leq M}\|\varepsilon^\ell\|_{L^2}^{2q}\right]\right)^{\frac{1}{2q}}+\left(\mathbb E\left[\left(\nu k\sum_{n=1}^{M}\|\partial_x\varepsilon^n\|_{L^2}^2\right)^q\right]\right)^{\frac{1}{2q}} \\
					&\leq2\left(1+\frac{1}{1-a q}\right)^{\frac{1}{2a q}}\left(\mathbb E\left[\exp\left(bq\sum_{n=0}^{M-1}G_n\right)\right]\right)^{\frac{1}{2b q}}\left(\mathbb E\left[\sup_{0\leq \ell < M}F_\ell\right]\right)^{\frac{1}{2}} .
				\end{aligned}
			\end{equation}
			
			It remains to estimate the terms on the right-hand side of \eqref{eq:discrete-gronwall-bound}. We first estimate the exponential factor involving \(G_n\). Set
			\begin{equation*}
				A:=\frac{80\alpha^2C_e^2bq}{\nu}, \qquad B:=bq\left(5+\frac{\nu}{4}+4\beta\mu+4L_\sigma\right).
			\end{equation*}
			This gives
			\begin{align*}
				\mathbb E\left[\exp\left(bq\sum_{n=0}^{M-1}G_n\right)\right] &= \exp(BT)\, \mathbb E\left[ \exp\left( A\int_0^T\|u(s)\|_{H^1}^2\,ds \right) \right] \\
				&= \exp(BT)\, \mathbb E\left[ \exp\left( A\int_0^T\|u(s)\|_{L^2}^2\,ds \right)\exp\left( A\int_0^T\|\partial_x u(s)\|_{L^2}^2\,ds \right) \right] \\
				&\leq \exp(BT)\, \Bigg(\mathbb E\left[ \exp\left( 2A\int_0^T\|u(s)\|_{L^2}^2\,ds \right)\right]\Bigg)^{\frac 12}\\
				&\quad\times\Bigg(\mathbb E\left[\exp\left( 2A\int_0^T\|\partial_x u(s)\|_{L^2}^2\,ds \right) \right]\Bigg)^{\frac 12}. 
			\end{align*}
			By the assumption on $K_\sigma$,
			\[
			K_\sigma \le \frac{\beta}{8L\kappa}= \frac{\beta\nu}{640L\alpha^2C_e^2bq}\min\left\{ \beta\gamma, \nu \right\} = \frac{\beta}{8L} \min\left\{ \frac{\beta \gamma}{A}, \frac{\nu}{A} \right\}.
			\]
			Therefore,
			\[
			\max\left\{\frac{A}{\beta \gamma},\frac{A}{\nu}\right\} \le \frac{\beta}{8K_\sigma L}.
			\]
			Then, by Lemma~\ref{lem:exponential-L2-estimate},
			\[
			\Bigg(\mathbb E\left[ \exp\left( 2A\int_0^T\|u(s)\|_{L^2}^2\,ds \right)\right]\Bigg)^{\frac 12}\Bigg(\mathbb E\left[\exp\left( 2A\int_0^T\|\partial_x u(s)\|_{L^2}^2\,ds \right) \right]\Bigg)^{\frac 12}\le C_{\ref{lem:exponential-L2-estimate}}.
			\]
			Consequently,
			\[
			\mathbb E\left[\exp\left(bq\sum_{n=0}^{M-1}G_n\right)\right]\leq C_{\ref{lem:exponential-L2-estimate}}\exp(BT)=:\widehat C_0.
			\]			
			Next, we estimate $\mathbb{E}\Big[\sup_{0\leq \ell < M} F_{\ell}\Big]$ as follows
			\begin{align}
				\mathbb E\left[\sup_{0\leq \ell < M} F_\ell\right] &\leq \sum_{n=0}^{M-1}\mathbb E\Big[Ckh^2 \left(1+\|u_h^{n+1}\|_{L^2}^4\right) \left(\|u(t_{n+1})\|_{H^2}^2+\|u(t_n)\|_{H^2}^2\right) \nonumber\\
				&\quad + Ckh^4 \left(1+\|u(t_{n+1})\|_{H^2}^6\right)\Big] \nonumber\\
				&\quad + \sum_{n=0}^{M-1}\mathbb E\Bigg[C\int_{t_n}^{t_{n+1}}\left(1+\|\partial_xu(s)\|_{L^2}^2 +\|u(t_{n+1})\|_{H^1}^2\right)\|u(t_{n+1})-u(s)\|_{H^1}^2\,ds \nonumber\\
				&\quad + C\int_{t_n}^{t_{n+1}}\left(1+\|u(t_{n+1})\|_{H^2}^4 +\|u(s)\|_{H^2}^4 \right)\|u(t_{n+1})-u(s)\|_{L^2}^2\,ds \Bigg] \nonumber\\
				&\quad +\sum_{n=0}^{M-1}\mathbb E\Bigg[Ch^4\left(\int_{t_n}^{t_{n+1}}\|u(t_{n+1})-u(s)\|_{H^1}^2\,ds + \int_{t_n}^{t_{n+1}}\|u(s)\|_{H^1}^2\,ds\right)\nonumber\\
				&\quad\quad\times \|u(t_{n+1})\|_{H^2}^2\Bigg] + \sum_{n=0}^{M-1}\mathbb E\left[4\left\|\int_{t_n}^{t_{n+1}}\bigl(\sigma(u(t_n))-\sigma(u(s))\bigr)\,dW(s)\right\|_{L^2}^2\right] \nonumber\\
				&\quad + \sum_{n=0}^{M-1}\mathbb E\Bigg[\frac{80\alpha^2C^2_e}{\nu}\Bigg(\int_{t_n}^{t_{n+1}}\|u(t_{n+1})-u(s)\|_{H^1}^2\,ds \Bigg)\|\varepsilon^{n+1}\|^2_{L^2}\Bigg] \nonumber\\
				&\quad + C\,\mathbb E\left[\sup_{0\leq \ell < M}\Bigg(k + \int_{t_\ell}^{t_{\ell+1}}\|u(t_{\ell+1})-u(s)\|_{H^1}^2\,ds + \int_{t_\ell}^{t_{\ell+1}} \|u(s)\|_{H^1}^2\,ds \Bigg)\|\varepsilon^{\ell+1}\|_{L^2}^2\right] \nonumber\\
				&=: T_1+T_2+T_3+T_4 + T_5 + T_6 \label{eq:discrete-gronwall-bound2}.
			\end{align}
			
			We estimate the terms \(T_1,\ldots,T_6\) separately. First, by Lemma~\ref{lem:Hm-q-moment} and Lemma~\ref{lem:stability-2}, we derive
			\begin{align*}
				T_1 &:= \sum_{n=0}^{M-1}\mathbb E\Big[Ckh^2\left(1+\|u_h^{n+1}\|_{L^2}^4\right) \left(\|u(t_{n+1})\|_{H^2}^2+\|u(t_n)\|_{H^2}^2\right) + Ckh^4 \left(1+\|u(t_{n+1})\|_{H^2}^6\right) \Big] \\
				&\leq Ckh^2\sum_{n=0}^{M-1}\Bigg[\left(\mathbb E\Big[\left(1+\|u_h^{n+1}\|_{L^2}^4\right)^2\Big]\right)^{1/2}\left(\mathbb E\Big[\left(\|u(t_{n+1})\|_{H^2}^2+\|u(t_n)\|_{H^2}^2\right)^2\Big]\right)^{1/2}\Bigg]\\
				&\quad\,+ Ckh^4\sum_{n=0}^{M-1}\mathbb E \Big[1+\|u(t_{n+1})\|_{H^2}^6\Big] \\
				&\leq Ckh^2\sum_{n=0}^{M-1}\Big[\left(\mathbb E\big[1+\|u_h^{n+1}\|_{L^2}^8\big]\right)^{1/2}\left(\mathbb E\big[\|u(t_{n+1})\|_{H^2}^4+\|u(t_n)\|_{H^2}^4\big]\right)^{1/2}\Big] \\
				&\quad\,+ Ckh^4\sum_{n=0}^{M-1}\mathbb E \Big[1+\|u(t_{n+1})\|_{H^2}^6\Big] \\
				&\leq CT\Big(1+C_{\ref{lem:stability-2}}^{0,8}\Big)^{\frac{1}{2}}\Big(2C^{2,4}_{\ref{lem:Hm-q-moment}}\Big)^{\frac{1}{2}}h^2 + CT\Big(1+C^{2,6}_{\ref{lem:Hm-q-moment}}\Big)h^4 =: C_{T_1}h^2.
			\end{align*}
			
			Using Lemma~\ref{lem:Hm-q-moment} and Lemma~\ref{lem:time-regularity}, we estimate \(T_2\) as follows:
			\begin{align*}
				T_2 &:= \sum_{n=0}^{M-1}\mathbb E\Bigg[C\int_{t_n}^{t_{n+1}}\left(1+\|\partial_xu(s)\|_{L^2}^2 +\|u(t_{n+1})\|_{H^1}^2\right)\|u(t_{n+1})-u(s)\|_{H^1}^2\,ds \\
				&\quad + C\int_{t_n}^{t_{n+1}}\left(1+\|u(t_{n+1})\|_{H^2}^4 +\|u(s)\|_{H^2}^4 \right)\|u(t_{n+1})-u(s)\|_{L^2}^2\,ds \Bigg] \\
				&\leq C \sum^{M-1}_{n=0} \int^{t_{n+1}}_{t_n} \mathbb{E} \big[\big(\|\partial_x u(s)\|^2_{L^2} + \|u(t_{n+1})\|^2_{H^1} + \|u(t_{n+1})\|^4_{H^2} + \|u(s)\|^4_{H^2}\big)\\
				&\quad\times\|u(t_{n+1}) - u(s)\|^2_{H^1}\big]\,ds + C \sum^{M-1}_{n=0} \int^{t_{n+1}}_{t_n} \mathbb{E} \big[\|u(t_{n+1}) - u(s)\|^2_{H^1}\big]\,ds\\
				&\leq C\sum_{n=0}^{M-1}\int_{t_n}^{t_{n+1}}\Big(\mathbb E\Big[\big(\|\partial_xu(s)\|_{L^2}^2+\|u(t_{n+1})\|_{H^1}^2+\|u(t_{n+1})\|_{H^2}^4+\|u(s)\|_{H^2}^4\big)^2\Big] \Big)^{1/2} \\
				&\qquad\times \big(\mathbb E\big[\|u(t_{n+1})-u(s)\|_{H^1}^4\big]\big)^{1/2}\,ds + C\sum_{n=0}^{M-1}\int_{t_n}^{t_{n+1}}\mathbb E\big[\|u(t_{n+1})-u(s)\|_{H^1}^2\big]\,ds \\
				&\leq C\sum_{n=0}^{M-1}\int_{t_n}^{t_{n+1}}\big(\mathbb E\big[\|\partial_xu(s)\|_{L^2}^4+\|u(t_{n+1})\|_{H^1}^4+\|u(t_{n+1})\|_{H^2}^8+\|u(s)\|_{H^2}^8\big]\big)^{1/2} \\
				&\qquad\times\big(\mathbb E\big[\|u(t_{n+1})-u(s)\|_{H^1}^4\big]\big)^{1/2}\,ds + C\sum_{n=0}^{M-1}\int_{t_n}^{t_{n+1}} \mathbb E\big[\|u(t_{n+1})-u(s)\|_{H^1}^2 \big]\,ds \\
				&\leq C \Big(2C^{1,4}_{\ref{lem:Hm-q-moment}} + 2C^{2,8}_{\ref{lem:Hm-q-moment}} \Big)^{\frac{1}{2}}\Big(C_{\ref{lem:time-regularity}}^{1,2}\Big)^{\frac{1}{2}} \frac{T}{2} k + C C_{\ref{lem:time-regularity}}^{1,1}\frac{T}{2}k =: C_{T_2}k.
			\end{align*}
			
			Using Lemma~\ref{lem:Hm-q-moment} and Lemma~\ref{lem:time-regularity}, we estimate \(T_3\) in the same manner as follows:
			\begin{align*}
				T_3 &:= \sum_{n=0}^{M-1}\mathbb E\left[Ch^4\left(\int_{t_n}^{t_{n+1}}\|u(t_{n+1}) -u(s)\|_{H^1}^2\,ds+\int_{t_n}^{t_{n+1}}\|u(s)\|_{H^1}^2\,ds\right)\|u(t_{n+1})\|_{H^2}^2 \right] \\
				&\leq Ch^4 \sum_{n=0}^{M-1}\int_{t_n}^{t_{n+1}}\mathbb E\left[\|u(t_{n+1})-u(s)\|_{H^1}^2\|u(t_{n+1})\|_{H^2}^2\right]ds \\
				&\quad + Ch^4 \sum_{n=0}^{M-1} \int_{t_n}^{t_{n+1}}\mathbb E\left[ \|u(s)\|_{H^1}^2 \|u(t_{n+1})\|_{H^2}^2 \right]ds \\
				&\leq Ch^4 \sum_{n=0}^{M-1} \int_{t_n}^{t_{n+1}} \left( \mathbb E\|u(t_{n+1})-u(s)\|_{H^1}^4 \right)^{1/2} \left( \mathbb E\|u(t_{n+1})\|_{H^2}^4\right)^{1/2}ds \\
				&\quad + Ch^4 \sum_{n=0}^{M-1} \int_{t_n}^{t_{n+1}}\left(\mathbb E\|u(s)\|_{H^1}^4\right)^{1/2}\left(\mathbb E\|u(t_{n+1})\|_{H^2}^4\right)^{1/2}ds \\
				&\leq CT\Big(C^{1,2}_{\ref{lem:time-regularity}}C^{2,4}_{\ref{lem:Hm-q-moment}}\Big)^{\frac12}kh^4 + CT\Big(C^{1,4}_{\ref{lem:Hm-q-moment}}C^{2,4}_{\ref{lem:Hm-q-moment}}\Big)^{\frac12}h^4 =: C_{T_3}h^4 .
			\end{align*}
			
			By It\^o's isometry, assumption~\eqref{ass:noise-lipschitz}, and Lemma~\ref{lem:time-regularity}, we have
			\begin{align*}
				T_4 &:= \sum_{n=0}^{M-1}\mathbb E\left[4\left\|\int_{t_n}^{t_{n+1}}\bigl(\sigma(u(t_n))-\sigma(u(s))\bigr)\,dW(s)\right\|_{L^2}^2\right]\\
				&\leq 4\sum_{n=0}^{M-1} \mathbb E\left[ \int_{t_n}^{t_{n+1}} \|\sigma(u(s))-\sigma(u(t_n))\|_{L^2}^2\,ds \right] \\
				&\leq 4L_{\sigma} \sum_{n=0}^{M-1} \mathbb E\left[ \int_{t_n}^{t_{n+1}} \|u(s)-u(t_n)\|_{L^2}^2\,ds \right] \\
				&\leq 2L_{\sigma}TC_{\ref{lem:time-regularity}}^{1,1}k =: C_{T_4}k.
			\end{align*}
			
			The estimate for $T_5$ follows from Lemmas~\ref{lem:Hm-q-moment}, \ref{lem:time-regularity}, and~\ref{lem:stability-2}
			\begin{align*}
				T_5 &:= \sum_{n=0}^{M-1}\mathbb E\Bigg[\frac{80\alpha^2C^2_e}{\nu}\Bigg(\int_{t_n}^{t_{n+1}}\|u(t_{n+1})-u(s)\|_{H^1}^2\,ds \Bigg)\|\varepsilon^{n+1}\|^2_{L^2}\Bigg] \\
				&\leq \frac{80\alpha^2C^2_e}{\nu}\mathbb E\Bigg[\Bigg(\sum_{n=0}^{M-1}\int_{t_n}^{t_{n+1}}\|u(t_{n+1})-u(s)\|_{H^1}^2\,ds \Bigg)\sup_{0 < \ell \leq M}\|\varepsilon^\ell\|^2_{L^2}\Bigg] \\
				&\leq \frac{80\alpha^2C^2_e}{\nu}\Bigg(\mathbb E\Bigg[\sum_{n=0}^{M-1}\int_{t_n}^{t_{n+1}}\|u(t_{n+1})-u(s)\|_{H^1}^2\,ds \Bigg]^2\Bigg)^{\frac12}\Bigg(\mathbb{E}\Bigg[\sup_{0 < \ell \leq M}\|\varepsilon^\ell\|^4_{L^2}\Bigg]\Bigg)^{\frac12} \\
				&\leq C\Bigg(T\sum_{n=0}^{M-1}\int_{t_n}^{t_{n+1}}\mathbb E\Big[\|u(t_{n+1})-u(s)\|_{H^1}^4\Big]\,ds \Bigg)^{\frac12}\Bigg(\mathbb{E}\Bigg[\sup_{0 < \ell \leq M}\|\varepsilon^\ell\|^4_{L^2}\Bigg]\Bigg)^{\frac12} \\
				&\leq C\Big(C^{1,2}_{\ref{lem:time-regularity}}TMk^3\Big)^{\frac12}\Big(C_{\ref{lem:Hm-q-moment}}^{0,4} + C_{\ref{lem:stability-2}}^{0,4}\Big)^{\frac{1}{2}} =: C_{T_5}k.
			\end{align*}
			
			By Lemmas~\ref{lem:Hm-q-moment} and~\ref{lem:stability-2}, we estimate \(T_6\) as follows:
			\begin{align*}
				T_6 &:= C\mathbb E\left[\sup_{0\leq \ell<M}\left(k+\int_{t_\ell}^{t_{\ell+1}}\|u(t_{\ell+1})-u(s)\|_{H^1}^2\,ds+\int_{t_\ell}^{t_{\ell+1}}\|u(s)\|_{H^1}^2\,ds\right)\|\varepsilon^{\ell+1}\|_{L^2}^2\right] \\
				&\leq Ck\,\mathbb E\left[\left(1+\sup_{0\leq t\leq T}\|u(t)\|_{H^1}^2\right)\sup_{0< \ell\leq M}\|\varepsilon^\ell\|_{L^2}^2\right] \\
				&\leq Ck\left(\mathbb E\left[\left(1+\sup_{0\leq t\leq T}\|u(t)\|_{H^1}^2\right)^2\right]\right)^{1/2}\left(\mathbb E\left[\sup_{0< \ell\leq M}\|\varepsilon^\ell\|_{L^2}^4\right]\right)^{1/2} \\
				&\leq Ck + C \Big(C_{\ref{lem:Hm-q-moment}}^{0,4}\Big)^{\frac{1}{2}} \Big(C_{\ref{lem:Hm-q-moment}}^{0,4} + C_{\ref{lem:stability-2}}^{0,4}\Big)^{\frac{1}{2}}k =: C_{T_6}k.
			\end{align*}
			
			Substituting the estimates for $T_1$, $T_2$,$\dots$,$T_6$ into \eqref{eq:discrete-gronwall-bound2} yields 
			\begin{align*}
				\Bigg(\mathbb E\Big[\sup_{0\leq \ell \leq M}\|\varepsilon^\ell\|_{L^2}^{2q}\Big]\Bigg)^{\frac{1}{2q}}+\left(\mathbb E\left[\left(\nu k\sum_{n=0}^{M}\|\partial_x \varepsilon^n\|_{L^2}^2\right)^q\right]\right)^{\frac{1}{2q}}\leq\widehat C\left(k^{\frac12}+h\right),
			\end{align*}
			where
			\begin{equation*}
				\widehat C := \Big(\widehat C_0\Big)^{\frac{1}{2bq}} \sqrt{\max\{C_{T_1},C_{T_2},C_{T_3},C_{T_4},C_{T_5},C_{T_6}\}}.
			\end{equation*}
			
			It remains to verify that \(\{M_\ell\}_{\ell=0}^{M-1}\) is a discrete martingale. Here, \(\mathcal F_{t_n}\) denotes the sigma-algebra of information available up to time \(t_n\) in the filtration associated with \(W\). In particular, \(u(t_n)\), \(u_h^n\), and \(\varepsilon^n\) are \(\mathcal F_{t_n}\)-measurable, while \(\Delta W_n\) is independent of \(\mathcal F_{t_n}\). Using the triangle inequality, It\^o's isometry, and assumption \eqref{ass:noise-lipschitz}, we obtain
			\begin{align*}
				\mathbb E[|M_\ell|] &\leq 8k\sum_{n=0}^{\ell}\mathbb E\left[\|\sigma(u(t_n))-\sigma(u_h^n)\|_{L^2}^2\right]  +2\mathbb E\left[\left|\sum_{n=0}^{\ell}\int_{t_n}^{t_{n+1}}\big(\sigma(u(s))-\sigma(u_h^n),\varepsilon^n\big)\,dW(s)\right|\right] \\
				&\leq 16TL_\sigma\left(C_{\ref{lem:Hm-q-moment}}^{0,2}+C_{\ref{lem:stability-2}}^{0,2}\right)  +2\left(\mathbb E\left[\sum_{n=0}^{\ell}\int_{t_n}^{t_{n+1}}\|\sigma(u(s))-\sigma(u_h^n)\|_{L^2}^2\|\varepsilon^n\|_{L^2}^2\,ds\right]\right)^{1/2} \\
				&\leq C\left(1+C_{\ref{lem:Hm-q-moment}}^{0,4}+C_{\ref{lem:stability-2}}^{0,4}\right)<\infty,
			\end{align*}
			where we used H\"older's inequality and the fourth-moment bounds in the last step. Moreover, \(Z_n\) is \(\mathcal F_{t_{n+1}}\)-measurable. By the independence of \(\Delta W_n\) from \(\mathcal F_{t_n}\) and the martingale property of the It\^o integral, we have
			\begin{align*}
				\mathbb E[Z_n\mid\mathcal F_{t_n}] &=4\|\sigma(u(t_n))-\sigma(u_h^n)\|_{L^2}^2\left(\mathbb E[|\Delta W_n|^2\mid\mathcal F_{t_n}]-k\right) \\
				&\quad +2\mathbb E\left[\int_{t_n}^{t_{n+1}}\big(\sigma(u(s))-\sigma(u_h^n),\varepsilon^n\big)\,dW(s)\,\middle|\,\mathcal F_{t_n}\right]=0.
			\end{align*}
			
			Consequently, for \(0\leq\ell\leq M-2\),
			\[
			\mathbb E[M_{\ell+1}\mid\mathcal F_{t_{\ell+1}}]
			=M_\ell+\mathbb E[Z_{\ell+1}\mid\mathcal F_{t_{\ell+1}}]=M_\ell.
			\]
			Thus, \(\{M_\ell\}_{\ell=0}^{M-1}\) is a discrete martingale with respect to \(\{\mathcal F_{t_{\ell+1}}\}_{\ell=0}^{M-1}\). The proof is complete.
		\end{proof}

		\subsection{Second-moment error estimates}\label{sec4.2}
		Theorem~\ref{thm:sub-second-moment-error} gives strong convergence in the discrete $L_\omega^r\ell_t^\infty L_x^2$- and $L_\omega^r\ell_t^2H_x^1$-norms for $0<r<2$. This restriction on $r$ comes from the stochastic Gronwall inequality in Lemma~\ref{lem:discrete-stochastic-gronwall}. In this subsection, we overcome this restriction by a higher-moment argument. Specifically, we derive a fourth-power error inequality and apply the stochastic Gronwall inequality with $q=\frac12$ to obtain the following second-moment estimate.
		
		\begin{theorem}\label{thm:second-moment-error}
			Let $u$ be the variational solution of \eqref{eq:week-solution}, and let $\{u_h^n\}_{n=0}^{M}$ be generated by the Main Algorithm with $u_h^0=P_hu_0$. Suppose that $u_0\in\bigcap_{i=0}^{2}L^{24\cdot3^{2-i}}(\Omega;H^i(D))$ and that $\sigma$ satisfies \eqref{ass:noise-lipschitz}. Assume further that $\mathbb E[\exp(2\kappa\|u_0\|_{L^2}^2)]<\infty$ and $K_\sigma\leq\frac{\beta}{8L\kappa}$, where
			\[
			\kappa=\frac{480\alpha^2C_e^2}{\nu}\max\left\{\frac{1}{\beta\gamma},\frac{1}{\nu}\right\}.
			\]
			Then
			\[
			\left(\mathbb E\left[\max_{1\leq n\leq M}\|u(t_n)-u_h^n\|_{L^2}^2\right]\right)^{1/2}
			\leq\widehat C_1\left(k^{1/2}+h\right),
			\]
			where $\widehat C_1=C(u_0,T,L_\sigma)>0$ is independent of $h$ and $k$.
		\end{theorem}
		\begin{proof}
			Multiplying both sides of $\eqref{eqn:3.27}$ by $\|\varepsilon^{n+1}\|^2_{L^2}$ and using the identity $2a(a-b) = a^2 - b^2 + (a-b)^2$, we obtain
			\begin{align}
				&\frac14\Big(\|\varepsilon^{n+1}\|_{L^2}^4-\|\varepsilon^n\|_{L^2}^4\Big)
				+ \frac14\Big(\|\varepsilon^{n+1}\|_{L^2}^2-\|\varepsilon^n\|_{L^2}^2\Big)^2
				+ \frac{\nu k}{2}\|\partial_x\varepsilon^{n+1}\|_{L^2}^2
				\|\varepsilon^{n+1}\|_{L^2}^2 \notag\\
				&\leq \frac{40\alpha^2 C_e^2}{\nu}
				\Bigg(
				\int_{t_n}^{t_{n+1}}\|u(t_{n+1})-u(s)\|_{H^1}^2\,ds
				+ \int_{t_n}^{t_{n+1}}\|u(s)\|_{H^1}^2\,ds
				\Bigg)\|\varepsilon^{n+1}\|_{L^2}^4 \notag\\
				&\quad + k\Bigg(\frac52+\frac{\nu}{8}+2\beta\mu\Bigg)
				\|\varepsilon^{n+1}\|_{L^2}^4
				+ 2L_{\sigma}k\|\varepsilon^n\|_{L^2}^2
				\|\varepsilon^{n+1}\|_{L^2}^2 \notag\\
				&\quad + Ckh^2\left(1+\|u_h^{n+1}\|_{L^2}^4\right)
				\left(\|u(t_{n+1})\|_{H^2}^2+\|u(t_n)\|_{H^2}^2\right)
				\|\varepsilon^{n+1}\|_{L^2}^2 \notag\\
				&\quad + Ckh^4\left(1+\|u(t_{n+1})\|_{H^2}^6\right)
				\|\varepsilon^{n+1}\|_{L^2}^2 \label{eqn:new-3.27}\\
				&\quad + C\int_{t_n}^{t_{n+1}}
				\left(1+\|\partial_xu(s)\|_{L^2}^2+\|u(t_{n+1})\|_{H^1}^2\right)
				\|u(t_{n+1})-u(s)\|_{H^1}^2
				\|\varepsilon^{n+1}\|_{L^2}^2\,ds \notag\\
				&\quad + C\int_{t_n}^{t_{n+1}}
				\left(1+\|u(t_{n+1})\|_{H^1}^4+\|u(s)\|_{H^1}^4\right)
				\|u(t_{n+1})-u(s)\|_{L^2}^2
				\|\varepsilon^{n+1}\|_{L^2}^2\,ds \notag\\
				&\quad + Ch^4\int_{t_n}^{t_{n+1}}
				\Big(\|u(t_{n+1})-u(s)\|_{H^1}^2+\|u(s)\|_{H^1}^2\Big)
				\|u(t_{n+1})\|_{H^2}^2
				\|\varepsilon^{n+1}\|_{L^2}^2\,ds \notag\\
				&\quad + 2\left\|
				\int_{t_n}^{t_{n+1}}
				\bigl(\sigma(u(s))-\sigma(u(t_n))\bigr)\,dW(s)
				\right\|_{L^2}^2\|\varepsilon^{n+1}\|_{L^2}^2 \notag\\
				&\quad + 2\|(\sigma(u(t_n))-\sigma(u_h^n))\Delta W_n\|_{L^2}^2
				\|\varepsilon^{n+1}\|_{L^2}^2 \notag\\
				&\quad + \left(
				\int_{t_n}^{t_{n+1}}
				\bigl(\sigma(u(s))-\sigma(u_h^n)\bigr)\,dW(s),
				\,\varepsilon^n
				\right)\|\varepsilon^{n+1}\|_{L^2}^2 \notag\\
				&=: X_1+X_2+\cdots+X_{11}. \notag
			\end{align}
			
			Using the Young inequality, we estimate $X_1$ as follows
			\begin{align*}
				X_1 &:= \frac{40 \alpha^2 C^2_e}{\nu}\Bigg(\int_{t_n}^{t_{n+1}}\|u(t_{n+1})-u(s)\|_{H^1}^2\,ds + \int_{t_n}^{t_{n+1}} \|u(s)\|_{H^1}^2\,ds \Bigg)\|\varepsilon^{n+1}\|_{L^2}^4 \\
				&\leq \frac{20 \alpha^2 C^2_e}{\nu}\Bigg(k + \int_{t_n}^{t_{n+1}}\|u(t_{n+1})-u(s)\|_{H^1}^4\,ds \Bigg)\|\varepsilon^{n+1}\|_{L^2}^4+ \frac{40 \alpha^2 C^2_e}{\nu}\left(\int_{t_n}^{t_{n+1}} \|u(s)\|_{H^1}^2\,ds\right)\|\varepsilon^{n+1}\|_{L^2}^4.
			\end{align*}
			
			Next, we estimate $X_2$ and $X_3$ by Young's inequality as follows:
			\begin{align*}
				X_2+X_3 :={}& k\left(\frac{5}{2}+\frac{\nu}{8}+2\beta\mu\right)\|\varepsilon^{n+1}\|_{L^2}^4+2L_\sigma k\|\varepsilon^n\|_{L^2}^2\|\varepsilon^{n+1}\|_{L^2}^2\\
				\leq{}&L_\sigma k\|\varepsilon^n\|_{L^2}^4+k\left(\frac{5}{2}+\frac{\nu}{8}+2\beta\mu+L_\sigma\right)\|\varepsilon^{n+1}\|_{L^2}^4.
			\end{align*}
			
			Similarly, for $X_4$ and $X_5$, Young's inequality yields
			\begin{align*}
				X_4 &:= Ckh^2 \left(1+\|u_h^{n+1}\|_{L^2}^4\right) \left(\|u(t_{n+1})\|_{H^2}^2 + \|u(t_n)\|_{H^2}^2 \right)\|\varepsilon^{n+1}\|_{L^2}^2\\
				&\leq Ckh^4 \left(1+\|u_h^{n+1}\|_{L^2}^8\right) \left(\|u(t_{n+1})\|_{H^2}^4 + \|u(t_n)\|_{H^2}^4 \right)+  k\|\varepsilon^{n+1}\|^4_{L^2},\\
				X_5 &= Ckh^4 \left( 1+\|u(t_{n+1})\|_{H^2}^6\right) \|\varepsilon^{n+1}\|_{L^2}^2 \\
				&\leq  Ckh^8 \left( 1+\|u(t_{n+1})\|_{H^2}^{12}\right) + k\|\varepsilon^{n+1}\|_{L^2}^4. 
			\end{align*}
			
			For $X_6$ and $X_7$, Young’s inequality gives
			\begin{align*}
				X_6 &:= C\int_{t_n}^{t_{n+1}} \left( 1+\|\partial_xu(s)\|_{L^2}^2 +\|u(t_{n+1})\|_{H^1}^2 \right) \|u(t_{n+1})-u(s)\|_{H^1}^2 \|\varepsilon^{n+1}\|_{L^2}^2\,ds \\
				&\leq C\int_{t_n}^{t_{n+1}} \left( 1+\|\partial_xu(s)\|_{L^2}^4 +\|u(t_{n+1})\|_{H^1}^4 \right) \|u(t_{n+1})-u(s)\|_{H^1}^4\,ds + k\|\varepsilon^{n+1}\|_{L^2}^4, \\   
				X_7 &:= C\int_{t_n}^{t_{n+1}} \left( 1+\|u(t_{n+1})\|_{H^1}^4 +\|u(s)\|_{H^1}^4\right) \|u(t_{n+1})-u(s)\|_{L^2}^2 \|\varepsilon^{n+1}\|_{L^2}^2\,ds  \\
				&\leq C\int_{t_n}^{t_{n+1}} \left( 1+\|u(t_{n+1})\|_{H^1}^8 +\|u(s)\|_{H^1}^8\right) \|u(t_{n+1})-u(s)\|_{L^2}^4\,ds  + k\|\varepsilon^{n+1}\|_{L^2}^4. 
			\end{align*}
			
			We rewrite the terms $X_9$, $X_{10}$, and $X_{11}$ by splitting $\|\varepsilon^{n+1}\|_{L^2}^2$ into $\|\varepsilon^{n}\|_{L^2}^2$ and $\|\varepsilon^{n+1}\|_{L^2}^2-\|\varepsilon^{n}\|_{L^2}^2$ as follows:
			\begin{align*}
				X_9 + X_{10} + X_{11} &:= 2\left\| \int_{t_n}^{t_{n+1}}\bigl(\sigma(u(s))-\sigma(u(t_n))\bigr)\,dW(s) \right\|_{L^2}^2 \|\varepsilon^{n+1}\|_{L^2}^2 \\
				&\quad + 2\|(\sigma(u(t_n))-\sigma(u_h^n))\Delta W_n\|_{L^2}^2\|\varepsilon^{n+1}\|_{L^2}^2
				\\
				&\quad+\left(\int_{t_n}^{t_{n+1}}\bigl(\sigma(u(s))-\sigma(u_h^n)\bigr)\,dW(s), \,\varepsilon^n\right)\|\varepsilon^{n+1}\|_{L^2}^2 \\
				&= 2\left\| \int_{t_n}^{t_{n+1}}\bigl(\sigma(u(s))-\sigma(u(t_n))\bigr)\,dW(s) \right\|_{L^2}^2 \|\varepsilon^{n}\|_{L^2}^2 \\
				&\quad + 2\|(\sigma(u(t_n))-\sigma(u_h^n))\Delta W_n\|_{L^2}^2\|\varepsilon^{n}\|_{L^2}^2
				\\
				&\quad+\left(\int_{t_n}^{t_{n+1}}\bigl(\sigma(u(s))-\sigma(u_h^n)\bigr)\,dW(s), \,\varepsilon^n\right)\|\varepsilon^{n}\|_{L^2}^2 \\
				&\quad+ 2\left\| \int_{t_n}^{t_{n+1}}\bigl(\sigma(u(s))-\sigma(u(t_n))\bigr)\,dW(s) \right\|_{L^2}^2 \big(\|\varepsilon^{n+1}\|_{L^2}^2 - \|\varepsilon^{n}\|_{L^2}^2\big)\\
				&\quad + 2\|(\sigma(u(t_n))-\sigma(u_h^n))\Delta W_n\|_{L^2}^2\big(\|\varepsilon^{n+1}\|_{L^2}^2 - \|\varepsilon^{n}\|_{L^2}^2\big)
				\\
				&\quad+\left(\int_{t_n}^{t_{n+1}}\bigl(\sigma(u(s))-\sigma(u_h^n)\bigr)\,dW(s), \,\varepsilon^n\right)\big(\|\varepsilon^{n+1}\|_{L^2}^2 - \|\varepsilon^{n}\|_{L^2}^2\big). 
			\end{align*}
			
			Applying the Cauchy--Schwarz and Young inequalities to the absolute values of the last three terms, we obtain
			\begin{align*}
				X_9 + X_{10} + X_{11} &\leq 2\left\| \int_{t_n}^{t_{n+1}}\bigl(\sigma(u(s))-\sigma(u(t_n))\bigr)\,dW(s) \right\|_{L^2}^2 \|\varepsilon^{n}\|_{L^2}^2\\
				&\quad +2\left\|\bigl(\sigma(u(t_n))-\sigma(u_h^n)\bigr)\Delta W_n \right\|_{L^2}^{2}\|\varepsilon^n\|_{L^2}^{2} \\
				&\quad +\left(\int_{t_n}^{t_{n+1}}\bigl(\sigma(u(s))-\sigma(u_h^n)\bigr)\,dW(s),\varepsilon^n\right)\|\varepsilon^n\|_{L^2}^{2} \\
				&\quad + 32\left\|\int_{t_n}^{t_{n+1}}\bigl(\sigma(u(s)) -\sigma(u(t_n))\bigr)\,dW(s)\right\|_{L^2}^{4} +32\left\|\bigl(\sigma(u(t_n))-\sigma(u_h^n)\bigr)\Delta W_n \right\|_{L^2}^{4} \\
				&\quad + 8\left\|\int_{t_n}^{t_{n+1}} \bigl(\sigma(u(s))-\sigma(u_h^n)\bigr)\,dW(s)\right\|_{L^2}^{2}\|\varepsilon^n\|_{L^2}^{2} +\frac{3}{32}\left(\|\varepsilon^{n+1}\|_{L^2}^{2}-\|\varepsilon^n\|_{L^2}^{2}\right)^2.
			\end{align*}
			
			Next, we add and subtract the corresponding quadratic variation terms to rewrite the stochastic terms in a martingale form
			\begin{align}
				X_9+X_{10}+X_{11} &\leq 2\Big[\left\|\bigl(\sigma(u(t_n))-\sigma(u_h^n)\bigr)\Delta W_n\right\|_{L^2}^2 -k\left\|\sigma(u(t_n))-\sigma(u_h^n)\right\|_{L^2}^2\Big]\|\varepsilon^n\|_{L^2}^2 \nonumber\\
				&\quad+ \left(\int_{t_n}^{t_{n+1}}\bigl(\sigma(u(s))-\sigma(u_h^n)\bigr)\,dW(s),\varepsilon^n\right)\|\varepsilon^n\|_{L^2}^2\nonumber\\
				&\quad+32\left\|\int_{t_n}^{t_{n+1}}\bigl(\sigma(u(s))-\sigma(u(t_n))\bigr)\,dW(s)\right\|_{L^2}^4 \nonumber\\
				&\quad+\Big[32\left\|\bigl(\sigma(u(t_n))-\sigma(u_h^n)\bigr)\Delta W_n\right\|_{L^2}^4 - 96k^2\left\|\sigma(u(t_n))-\sigma(u_h^n)\right\|_{L^2}^4 \Big]\nonumber\\
				&\quad+ 8\Bigg[\left\|\int_{t_n}^{t_{n+1}}\bigl(\sigma(u(s)) -\sigma(u_h^n)\bigr)\,dW(s)\right\|_{L^2}^2 - \int_{t_n}^{t_{n+1}}\left\|\sigma(u(s))-\sigma(u_h^n)\right\|_{L^2}^2\,ds\Bigg]\|\varepsilon^n\|_{L^2}^2 \nonumber\\
				&\quad+2k\left\|\sigma(u(t_n))-\sigma(u_h^n)\right\|_{L^2}^2\|\varepsilon^n\|_{L^2}^2 + 8\int_{t_n}^{t_{n+1}}\left\|\sigma(u(s))-\sigma(u_h^n)\right\|_{L^2}^2\|\varepsilon^n\|_{L^2}^2\,ds \label{term:01}\\
				&\quad+96k^2\left\|\sigma(u(t_n))-\sigma(u_h^n)\right\|_{L^2}^4 +\frac{3}{32}\left(\|\varepsilon^{n+1}\|_{L^2}^2-\|\varepsilon^n\|_{L^2}^2\right)^2 \label{term:02}\\\
				&\quad+2\left\|\int_{t_n}^{t_{n+1}}\bigl(\sigma(u(s))-\sigma(u(t_n))\bigr)\,dW(s)\right\|_{L^2}^2\|\varepsilon^n\|_{L^2}^2.\label{term:03}
			\end{align}
			
			We next estimate the terms in \eqref{term:01}--\eqref{term:03}. Using assumption \eqref{ass:noise-lipschitz}, Young's inequality, and the $L^2$-projection estimate \eqref{eq:projector-inequality}, we obtain
			\begin{align*}
				2k\|\sigma(u(t_n))-\sigma(u_h^n)\|_{L^2}^{2} \|\varepsilon^n\|_{L^2}^{2} &\le 2L_\sigma k \|u(t_n)-u_h^n\|_{L^2}^{2} \|\varepsilon^n\|_{L^2}^{2} \\
				&\le 4L_\sigma k\left(\|u(t_n)-P_hu(t_n)\|_{L^2}^{2}+\|\varepsilon^n\|_{L^2}^{2}\right)\|\varepsilon^n\|_{L^2}^{2} \\
				&\le2L_\sigma k\|u(t_n)-P_hu(t_n)\|_{L^2}^{4}+6L_\sigma k\|\varepsilon^n\|_{L^2}^{4} \\
				&\le2L_\sigma C_P^{4}kh^{8}\|u(t_n)\|_{H^2}^{4}+6L_\sigma k\|\varepsilon^n\|_{L^2}^{4}.
			\end{align*}
			
			Applying assumption \eqref{ass:noise-lipschitz} together with Young's inequality yields
			\begin{align*}
				8\int_{t_n}^{t_{n+1}}\left\|\sigma(u(s))- \sigma(u_h^n)\right\|_{L^2}^2\|\varepsilon^n\|_{L^2}^2\,ds 
				&\leq 4L_\sigma k\|\varepsilon^{n}\|_{L^{2}}^{4} +4L_\sigma\int_{t_n}^{t_{n+1}}\|u(s)-u_h^{n}\|_{L^{2}}^{4}\,ds \\
				&\leq 260L_{\sigma}k\|\varepsilon^n\|_{L^2}^{4} +32L_\sigma\int_{t_n}^{t_{n+1}}\|u(t_n) - u(s)\|_{L^{2}}^{4}\,ds \\
				&\quad + 256L_{\sigma}C_P^4kh^8\|u(t_n)\|^4_{H^2}.
			\end{align*}
			
			Similarly, by assumption \eqref{ass:noise-lipschitz} and the $L^2$-projection estimate \eqref{eq:projector-inequality}, we obtain
			\begin{align*}
				96k^2 \|\sigma(u(t_n))-\sigma(u_h^n)\|_{L^2}^{4} &\le 96L_\sigma^2 k^2 \|u(t_n)-u_h^n\|_{L^2}^{4}\\
				&\le 768L_\sigma^2 k^2 \left( \|u(t_n)-P_hu(t_n)\|_{L^2}^{4} +\|\varepsilon^n\|_{L^2}^{4}\right)\\
				&\le 768L_\sigma^2 C_P^4 k^2h^8 \|u(t_n)\|_{H^2}^{4} + 768L_\sigma^2 k^2 \|\varepsilon^n\|_{L^2}^{4}.
			\end{align*}
			
			By Young's inequality, we have
			\begin{align*}
				&2\left\|\int_{t_n}^{t_{n+1}}\bigl(\sigma(u(s))-\sigma(u(t_n))\bigr)\,dW(s)\right\|_{L^2}^2\|\varepsilon^n\|_{L^2}^2 \\
				&\qquad\leq\frac{1}{k}\left\|\int_{t_n}^{t_{n+1}}\bigl(\sigma(u(s))-\sigma(u(t_n))\bigr)\,dW(s)\right\|_{L^2}^{4}+k\|\varepsilon^n\|_{L^2}^{4}.
			\end{align*}
			
			Now, collecting the estimates for \(X_1,\ldots,X_{11}\), substituting them into \eqref{eqn:new-3.27}, absorbing the corresponding terms into the left-hand side of \eqref{eqn:new-3.27}, and summing over \(n=0,\ldots,\ell\) for any \(0\le \ell\le M-1\), we have
			\begin{align}
				&\|\varepsilon^{\ell+1}\|_{L^2}^4 + \frac{5}{8}\sum^{\ell}_{n=0}\Bigl(\|\varepsilon^{n+1}\|_{L^2}^{2}-\|\varepsilon^{n}\|_{L^2}^{2}\Bigr)^2+2\nu k\sum^{\ell}_{n=0}\|\partial_x\varepsilon^{n+1}\|_{L^2}^{2}\|\varepsilon^{n+1}\|_{L^2}^{2}\notag \\
				&\leq Ckh^4 \sum^{\ell}_{n=0}\left(1+\|u_h^{n+1}\|_{L^2}^8\right) \left(\|u(t_{n+1})\|_{H^2}^4 + \|u(t_n)\|_{H^2}^4 \right) \notag\\
				&\quad+ Ckh^8 \sum^{\ell}_{n=0}\left( 1+\|u(t_{n+1})\|_{H^2}^{12}\right) \notag \\
				&\quad + C\sum^{\ell}_{n=0}\int_{t_n}^{t_{n+1}} \left( 1+\|\partial_xu(s)\|_{L^2}^4 +\|u(t_{n+1})\|_{H^1}^4 \right) \|u(t_{n+1})-u(s)\|_{H^1}^4\,ds \notag \\
				&\quad + C\sum^{\ell}_{n=0}\int_{t_n}^{t_{n+1}} \left( 1+\|u(t_{n+1})\|_{H^1}^8 +\|u(s)\|_{H^1}^8\right) \|u(t_{n+1})-u(s)\|_{L^2}^4\,ds  \notag  \\
				&\quad +Ch^4\sum^{\ell}_{n=0}\int_{t_n}^{t_{n+1}}\big(\|u(t_{n+1})-u(s)\|_{H^1}^2+\|u(s)\|_{H^1}^2\big) \|u(t_{n+1})\|_{H^2}^2\|\varepsilon^{n+1}\|_{L^2}^2 \,ds \notag \\
				&\quad +4\left(\frac{1}{k}+32\right)\sum^{\ell}_{n=0}\left\|\int_{t_n}^{t_{n+1}}\bigl(\sigma(u(s))-\sigma(u(t_n))\bigr)\,dW(s)\right\|_{L^2}^{4}\notag\\
				&\quad +4 L_{\sigma}C^4_Pkh^8(258 +768L_{\sigma}k)\sum^{\ell}_{n=0}\|u(t_n)\|_{H^2}^{4}+128L_\sigma\sum^{\ell}_{n=0}\int_{t_n}^{t_{n+1}}\|u(t_n) - u(s)\|_{L^{2}}^{4}\,ds\notag\\
				&\quad +C\sum^{\ell}_{n=0}\Bigg(\int_{t_n}^{t_{n+1}}\|u(t_{n+1})-u(s)\|_{H^1}^4\,ds\Bigg)\|\varepsilon^{n+1}\|_{L^2}^4 \notag \\
				&\quad+ \sum^{\ell}_{n=0} 8\Big[\left\|\bigl(\sigma(u(t_n))-\sigma(u_h^n)\bigr)\Delta W_n\right\|_{L^2}^2-k\left\|\sigma(u(t_n))-\sigma(u_h^n)\right\|_{L^2}^2\Big]\|\varepsilon^n\|_{L^2}^2\notag \\
				&\quad+\sum^{\ell}_{n=0} 4\left(\int_{t_n}^{t_{n+1}}\bigl(\sigma(u(s))-\sigma(u_h^n)\bigr)\,dW(s),\varepsilon^n\right)\|\varepsilon^n\|_{L^2}^2\notag \\
				&\quad+\sum^{\ell}_{n=0} 4\Big[32\left\|\bigl(\sigma(u(t_n))-\sigma(u_h^n)\bigr)\Delta W_n\right\|_{L^2}^{4}-96k^2\left\|\sigma(u(t_n))-\sigma(u_h^n)\right\|_{L^2}^{4}\Big]\notag\\
				&\quad+\sum^{\ell}_{n=0} 32\Bigg[\left\|\int_{t_n}^{t_{n+1}}\bigl(\sigma(u(s))-\sigma(u_h^n)\bigr)\,dW(s)\right\|_{L^2}^{2}-\int_{t_n}^{t_{n+1}}\left\|\sigma(u(s))-\sigma(u_h^n)\right\|_{L^2}^{2}\,ds\Bigg]\|\varepsilon^n\|_{L^2}^2 \notag \\
				&\quad+4k\big(768L^2_{\sigma}k + 267L_{\sigma} + 1\big)\sum^{\ell}_{n=0}\|\varepsilon^{n}\|^4_{L^2} \notag \\
				&\quad+ 4k\sum^{\ell}_{n=0}\left(\frac{40\alpha^2C_e^2}{\nu k}\int_{t_n}^{t_{n+1}}\|u(s)\|_{H^1}^2\,ds + \frac{20\alpha^2C^2_e}{\nu} +\frac{15}{2} + \frac{\nu}{8} + 2\beta \mu + L_{\sigma}\right) \|\varepsilon^{n+1}\|^4_{L^2}. \notag 
			\end{align}
			
			To separate the resulting higher-order term in \(k\), we define. Set
			\begin{equation*}
				R_{n+1} := \frac{40\alpha^2 C_{e}^2}{\nu k}\int_{t_n}^{t_{n+1}}\|u(s)\|_{H^1}^{2}\,ds+\frac{20\alpha^2 C_{\varepsilon}^2}{\nu}+\frac{15}{2}+\frac{\nu}{8}+2\beta\mu+L_{\sigma}.
			\end{equation*}
			
			Separating the last term in the sum, using the identity $a^2-b^2 = (a+b)(a-b)$, and then applying Young’s inequality, we obtain
			\begin{align*} 
				4k\sum_{n=0}^{\ell}R_{n+1}\|\varepsilon^{n+1}\|_{L^2}^{4}={}&4k\sum_{n=0}^{\ell-1}R_{n+1}\|\varepsilon^{n+1}\|_{L^2}^{4}+4kR_{\ell+1}\|\varepsilon^{\ell+1}\|_{L^2}^{4}\\
				={}& 4k\sum_{n=0}^{\ell-1} R_{n+1}\|\varepsilon^{n+1}\|_{L^2}^{4} + 4kR_{\ell+1}\|\varepsilon^\ell\|_{L^2}^{4}\\	&\quad+4kR_{\ell+1}\left(\|\varepsilon^{\ell+1}\|_{L^2}^{2}+\|\varepsilon^\ell\|_{L^2}^{2}\right)\left( \|\varepsilon^{\ell+1}\|_{L^2}^{2}-\|\varepsilon^\ell\|_{L^2}^{2}\right)\\
				\le{}&4k\sum_{n=0}^{\ell}(R_n+R_{n+1})\|\varepsilon^n\|_{L^2}^{4}+32k^2R_{\ell+1}^{2}\left(\|\varepsilon^{\ell+1}\|_{L^2}^{2}+\|\varepsilon^\ell\|_{L^2}^{2}\right)^2 \\
				&\quad+ \frac{1}{8}\left(\|\varepsilon^{\ell+1}\|_{L^2}^{2}-\|\varepsilon^\ell\|_{L^2}^{2}\right)^2.
			\end{align*}
			
			Once again, we derive a new estimate for \eqref{eqn:new-3.27} as follows:
			\begin{equation}\label{eq:new-theo-before-gronwall}
				\begin{aligned}
					\|\varepsilon^{\ell+1}\|_{L^2}^{4} + \frac{5}{8}\sum_{n=0}^{\ell-1} \Bigl(\|\varepsilon^{n+1}&\|_{L^2}^{2}-\|\varepsilon^{n}\|_{L^2}^{2} \Bigr)^2 +2\nu k\sum_{n=0}^{\ell} \|\partial_x\varepsilon^{n+1}\|_{L^2}^{2}\|\varepsilon^{n+1}\|_{L^2}^{2}\\
					&\leq F_{\ell} + M_{\ell} + \sum^{\ell-1}_{n=0} G_n \|\varepsilon^{n+1}\|^4_{L^2},
				\end{aligned}
			\end{equation}
			where
			\begin{align*}
				F_{\ell} &:= Ckh^4 \sum^{\ell}_{n=0}\left(1+\|u_h^{n+1}\|_{L^2}^8\right) \left(\|u(t_{n+1})\|_{H^2}^4 + \|u(t_n)\|_{H^2}^4 \right) + Ckh^8 \sum^{\ell}_{n=0}\left( 1+\|u(t_{n+1})\|_{H^2}^{12}\right) \\
				&\quad + C\sum^{\ell}_{n=0}\int_{t_n}^{t_{n+1}} \left( 1+\|\partial_xu(s)\|_{L^2}^4 +\|u(t_{n+1})\|_{H^1}^4 \right) \|u(t_{n+1})-u(s)\|_{H^1}^4\,ds \\
				&\quad + C\sum^{\ell}_{n=0}\int_{t_n}^{t_{n+1}} \left( 1+\|u(t_{n+1})\|_{H^1}^8 +\|u(s)\|_{H^1}^8\right) \|u(t_{n+1})-u(s)\|_{L^2}^4\,ds   \\
				&\quad +Ch^4\sum^{\ell}_{n=0}\int_{t_n}^{t_{n+1}}\big(\|u(t_{n+1})-u(s)\|_{H^1}^2+\|u(s)\|_{H^1}^2\big) \|u(t_{n+1})\|_{H^2}^2\|\varepsilon^{n+1}\|_{L^2}^2 \,ds \\
				&\quad +4\left(\frac{1}{k}+32\right)\sum^{\ell}_{n=0}\left\|\int_{t_n}^{t_{n+1}}\bigl(\sigma(u(s))-\sigma(u(t_n))\bigr)\,dW(s)\right\|_{L^2}^{4}\\
				&\quad +4 L_{\sigma}C^4_Pkh^8(258 +768L_{\sigma}k)\sum^{\ell}_{n=0}\|u(t_n)\|_{H^2}^{4}+128L_\sigma\sum^{\ell}_{n=0}\int_{t_n}^{t_{n+1}}\|u(t_n) - u(s)\|_{L^{2}}^{4}\,ds\\
				&\quad +C\sum^{\ell}_{n=0}\Bigg(\int_{t_n}^{t_{n+1}}\|u(t_{n+1})-u(s)\|_{H^1}^4\,ds\Bigg)\|\varepsilon^{n+1}\|_{L^2}^4 \\
				&\quad +32k^2\left( \frac{40\alpha^2 C_e^2}{\nu k}\int_{t_\ell}^{t_{\ell+1}}\|u(s)\|_{H^1}^{2}\,ds+\frac{20\alpha^2 C_e^2}{\nu}+\frac{15}{2}+\frac{\nu}{8}+2\beta\mu+L_{\sigma}\right)^2\\
				&\quad \quad \times\left(\|\varepsilon^{\ell+1}\|_{L^2}^{2}+\|\varepsilon^\ell\|_{L^2}^{2}\right)^2,\\
				M_{\ell} &:= \sum^{\ell}_{n=0}\mathcal{M}_{n},\\
				\mathcal{M}_{n} &:= 8\Big[\left\|\bigl(\sigma(u(t_n))-\sigma(u_h^n)\bigr)\Delta W_n\right\|_{L^2}^2-k\left\|\sigma(u(t_n))-\sigma(u_h^n)\right\|_{L^2}^2\Big]\|\varepsilon^n\|_{L^2}^2\\
				&\quad+4\left(\int_{t_n}^{t_{n+1}}\bigl(\sigma(u(s))-\sigma(u_h^n)\bigr)\,dW(s),\varepsilon^n\right)\|\varepsilon^n\|_{L^2}^2\\
				&\quad+4\Big[32\left\|\bigl(\sigma(u(t_n))-\sigma(u_h^n)\bigr)\Delta W_n\right\|_{L^2}^{4}-96k^2\left\|\sigma(u(t_n))-\sigma(u_h^n)\right\|_{L^2}^{4}\Big]\\
				&\quad+32\Bigg[\left\|\int_{t_n}^{t_{n+1}}\bigl(\sigma(u(s))-\sigma(u_h^n)\bigr)\,dW(s)\right\|_{L^2}^{2}-\int_{t_n}^{t_{n+1}}\left\|\sigma(u(s))-\sigma(u_h^n)\right\|_{L^2}^{2}\,ds\Bigg]\|\varepsilon^n\|_{L^2}^2,\\
				G_{n} &:= 4\int^{t_{n+1}}_{t_n} \left(\frac{80\alpha^2C_e^2}{\nu}\|u(s)\|^2_{H^1} + \frac{40\alpha^2C_e^2}{\nu} + \frac{\nu}{4} +768L_{\sigma}^2k +269L_{\sigma} + 16\right)ds.
			\end{align*}
			
			Next, temporarily assuming the martingale property of $\{M_\ell : \ell \geq 0\}$, which will be verified below, we apply the stochastic Gronwall inequality in Lemma \ref{eq:discrete-stochastic-gronwall-assumption} to~\eqref{eq:new-theo-before-gronwall} with $q=\frac{1}{2}$, $\alpha=\frac{3}{2}$, and $\beta=3$, which yields
			\begin{equation}\label{eq:new-discrete-gronwall-bound}
				\left(\mathbb{E}\left[\sup_{0 \le \ell \le M}\|\varepsilon^\ell\|_{L^2}^{2}\right]\right)^{\frac12}\le\sqrt[3]{5}\left(\mathbb{E}\left[\exp\left(\frac32\sum_{n=0}^{M-1} G_n\right)\right]\right)^{\frac16}\left(\mathbb{E}\left[\sup_{0\le \ell<M}F_\ell\right]\right)^{\frac14}.
			\end{equation}
			
			We first estimate the exponential factor involving \(G_n\). Set
			\[
			A:=\frac{480\alpha^2C_e^2}{\nu} \quad\text{and} \quad B:=6\left( \frac{40\alpha^2C_e^2}{\nu} + \frac{\nu}{4} +768L_{\sigma}^2k +269L_{\sigma} + 16\right).
			\]
			
			This gives
			\begin{align*}
				\mathbb E\left[\exp\left(\frac{3}{2}\sum_{n=0}^{M-1}G_n\right)\right] &= \exp(BT)\, \mathbb E\left[ \exp\left( A\int_0^T\|u(s)\|_{H^1}^2\,ds \right) \right] \\
				&= \exp(BT)\, \mathbb E\left[ \exp\left( A\int_0^T\|u(s)\|_{L^2}^2\,ds \right)\exp\left( A\int_0^T\|\partial_x u(s)\|_{L^2}^2\,ds \right) \right] \\
				&\leq \exp(BT)\, \Bigg(\mathbb E\left[ \exp\left( 2A\int_0^T\|u(s)\|_{L^2}^2\,ds \right)\right]\Bigg)^{\frac 12}\\
				&\quad\times\Bigg(\mathbb E\left[\exp\left( 2A\int_0^T\|\partial_x u(s)\|_{L^2}^2\,ds \right) \right]\Bigg)^{\frac 12}. 
			\end{align*}
			
			By the assumption on $K_\sigma$,
			\[
			K_\sigma \le \frac{\beta}{8L\kappa}= \frac{\beta\nu}{3840\alpha^2C_e^2}\min\left\{ \beta\gamma, \nu \right\} = \frac{\beta}{8L} \min\left\{ \frac{\beta \gamma}{A}, \frac{\nu}{A} \right\}.
			\]
			
			Therefore,
			\[
			\max\left\{\frac{A}{\beta \gamma},\frac{A}{\nu}\right\} \le \frac{\beta}{8K_\sigma L}.
			\]
			
			Then, by Lemma~\ref{lem:exponential-L2-estimate},
			\[
			\Bigg(\mathbb E\left[ \exp\left( 2A\int_0^T\|u(s)\|_{L^2}^2\,ds \right)\right]\Bigg)^{\frac 12}\Bigg(\mathbb E\left[\exp\left( 2A\int_0^T\|\partial_x u(s)\|_{L^2}^2\,ds \right) \right]\Bigg)^{\frac 12}\le C_{\ref{lem:exponential-L2-estimate}}.
			\]
			Consequently,
			\[
			\sqrt[3]{5}\left(\mathbb{E}\left[\exp\left(\frac32\sum_{n=0}^{M-1} G_n\right)\right]\right)^{\frac16}\leq  \sqrt[3]{5} \left(\exp(BT)C_{\ref{lem:exponential-L2-estimate}}\right)^{\frac16}=:\widehat C_0.
			\]
			
			Next, we estimate $\mathbb{E}\left[\sup_{0\le \ell<M}F_{\ell}\right] =: T_1 + T_2 + \cdots + T_{9}$. By Lemmas \ref{lem:Hm-q-moment} and \ref{lem:stability-2}, we first estimate $T_1$ as follows:
			\begin{align*}
				T_1 &:= Ckh^4\sum_{n=0}^{M-1}\mathbb{E}\left[\left(1+\|u_h^{n+1}\|_{L^2}^{8}\right)\left(\|u(t_{n+1})\|_{H^2}^{4}+\|u(t_n)\|_{H^2}^{4}\right)\right]\\
				&\leq Ckh^4\sum_{n=0}^{M-1}\left(\mathbb{E}\left[\left(1+\|u_h^{n+1}\|_{L^2}^{8}\right)^2\right]\right)^{1/2}\left(\mathbb{E}\left[\left(\|u(t_{n+1})\|_{H^2}^{4}+\|u(t_n)\|_{H^2}^{4}\right)^2\right]\right)^{1/2}\\
				&\leq Ckh^4\sum_{n=0}^{M-1}\left(1+\mathbb{E}\Big[\|u_h^{n+1}\|_{L^2}^{16}\Big]\right)^{1/2}\left(\mathbb{E}\Big[\|u(t_{n+1})\|_{H^2}^{8}\Big]+\mathbb{E}\Big[\|u(t_n)\|_{H^2}^{8}\Big]\right)^{1/2} \\
				&\leq Ckh^4\sum_{n=0}^{M-1}\left(1+C_{\ref{lem:stability-2}}^{0,16}\right)^{1/2}\left(2C_{\ref{lem:Hm-q-moment}}^{2,8}\right)^{1/2} \\
				&\leq CTh^4 =: C_{T_1}h^4.
			\end{align*}
			
			For $T_2$, using Lemma \ref{lem:Hm-q-moment}, we have
			\begin{align*}
				T_2 &:= Ckh^8\sum_{n=0}^{M-1}\mathbb E\left[1+\|u(t_{n+1})\|_{H^2}^{12}\right] \leq CTh^8 + CTh^8C^{2,12}_{\ref{lem:Hm-q-moment}} =: C_{T_2}h^8.
			\end{align*}
			
			Using H\"older’s inequality, we estimate $T_3$ as follows:
			\begin{align*}
				T_3 &:= C\,\mathbb E\Bigg[\sum_{n=0}^{M-1}\int_{t_n}^{t_{n+1}}\left(1+\|\partial_xu(s)\|_{L^2}^4+\|u(t_{n+1})\|_{H^1}^4\right)\|u(t_{n+1})-u(s)\|_{H^1}^4\,ds\Bigg]\\
				&\leq C\sum_{n=0}^{M-1}\int_{t_n}^{t_{n+1}}\left\{\mathbb{E}\left[1+\|\partial_xu(s)\|_{L^2}^{16}+\|u(t_{n+1})\|_{H^1}^{16}\right]\right\}^{1/4}\left\{\mathbb{E}\left[\|u(t_{n+1})-u(s)\|_{H^1}^{16}\right]\right\}^{1/4}\,ds\\
				&\leq C\left(1+2C_{\ref{lem:Hm-q-moment}}^{1,16}\right)^{1/4}\sum_{n=0}^{M-1}\int_{t_n}^{t_{n+1}} \left(C^{1,8}_{\ref{lem:time-regularity}}\right)^{1/4}(t_{n+1}-s)^2\,ds \\
				&\leq CTk^2\left(1+2C_{\ref{lem:Hm-q-moment}}^{1,16}\right)^{1/4}\left(C^{1,8}_{\ref{lem:time-regularity}}\right)^{1/4} =: C_{T_3}k^2.
			\end{align*}
			
			Similarly for $T_4$,
			\begin{align*}
				T_4 &:= C\,\mathbb E\Bigg[\sum_{n=0}^{M-1}\int_{t_n}^{t_{n+1}}\left(1+\|u(t_{n+1})\|_{H^1}^8+\|u(s)\|_{H^1}^8\right)\|u(t_{n+1})-u(s)\|_{L^2}^4\,ds\Bigg]\\
				&\leq C\left(1+2C_{\ref{lem:Hm-q-moment}}^{1,16}\right)^{1/2}\sum_{n=0}^{M-1}\int_{t_n}^{t_{n+1}} \left(C^{0,8}_{\ref{lem:time-regularity}}\right)^{1/2}(t_{n+1}-s)^2\,ds  \\
				&\leq CTk^2\left(1+2C_{\ref{lem:Hm-q-moment}}^{1,16}\right)^{1/2}\left(C^{0,8}_{\ref{lem:time-regularity}}\right)^{1/2}=: C_{T_4}k^2.
			\end{align*}
			
			Next, we estimate $T_5$.
			\begin{align*}
				T_5 &:= Ch^4\mathbb E\Bigg[\sum_{n=0}^{M-1}\int_{t_n}^{t_{n+1}}\left(\|u(t_{n+1})-u(s)\|_{H^1}^2+\|u(s)\|_{H^1}^2\right)\|u(t_{n+1})\|_{H^2}^2\,ds \sup_{0< j\le M}\|\varepsilon^j\|_{L^2}^2\Bigg]\\
				&\le Ch^4\sum_{n=0}^{M-1}\int_{t_n}^{t_{n+1}}\left\{\mathbb{E}\left[\left(\|u(t_{n+1})-u(s)\|_{H^1}^{2}+\|u(s)\|_{H^1}^{2}\right)^2\|u(t_{n+1})\|_{H^2}^{4}\right]\right\}^{1/2}ds \\
				&\qquad\times\left(\mathbb{E}\left[\sup_{0<j\le M}\|\varepsilon^j\|_{L^2}^{4}\right]\right)^{1/2} \\
				&\le Ch^4\sum_{n=0}^{M-1}\int_{t_n}^{t_{n+1}}\left\{\mathbb{E}\left[\left(\|u(t_{n+1})-u(s)\|_{H^1}^{2}+\|u(s)\|_{H^1}^{2}\right)^4\right]\right\}^{1/4}\left\{\mathbb{E}\left[\|u(t_{n+1})\|_{H^2}^{8}\right]\right\}^{1/4}ds\\
				&\qquad\times\left(\mathbb{E}\left[\sup_{0<j\le M}\|\varepsilon^j\|_{L^2}^{4}\right]\right)^{1/2} \\
				&\le Ch^4\sum_{n=0}^{M-1}\int_{t_n}^{t_{n+1}}\left[C_{\ref{lem:time-regularity}}^{1,4}(t_{n+1}-s)^4+C_{\ref{lem:Hm-q-moment}}^{1,8}\right]^{1/4}\left(C_{\ref{lem:Hm-q-moment}}^{2,8}\right)^{1/4}ds \left(C^{0,4}_{\ref{lem:stability-2}} + C^{0,4}_{\ref{lem:Hm-q-moment}}\right)^{1/2} \\
				&=: C_{T_5}h^4.
			\end{align*}
			
			For $T_6$, the Burkh\"older–Davis–Gundy inequality yields
			\begin{align*}
				T_6 &:= 4\left(\frac1k+32\right)\sum_{n=0}^{M-1}\mathbb E\left[\left\|\int_{t_n}^{t_{n+1}}\bigl(\sigma(u(t_n))-\sigma(u(s))\bigr)\,dW(s)\right\|_{L^2}^4\right] \\
				&\leq C\left(\frac{1}{k}+32\right)\sum_{n=0}^{M-1}\mathbb{E}\left[\left(\int_{t_n}^{t_{n+1}}\|\sigma(u(t_n))-\sigma(u(s))\|_{L^2}^{2}\,ds\right)^2\right] \\
				&\leq C\left(\frac{1}{k}+32\right)\sum_{n=0}^{M-1}\mathbb{E}\left[\left(\int_{t_n}^{t_{n+1}}\|u(t_n)-u(s)\|_{L^2}^{2}\,ds\right)^2\right] \\
				&\leq C\left(\frac{1}{k}+32\right)k\sum_{n=0}^{M-1}\int_{t_n}^{t_{n+1}}C^{0,2}_{\ref{lem:time-regularity}}(t_n - s)^2\,ds \\
				&\leq C\left(\frac{1}{k}+32\right)Mk^4 =: C_{T_6}k^2.
			\end{align*}
			
			Applying Lemmas~\ref{lem:Hm-q-moment} and \ref{lem:time-regularity}, we estimate $T_7$ as follows:
			\begin{align*}
				T_7 &:= 4L_\sigma C_P^4 k h^8\left(2+768L_\sigma k\right)\sum^{M-1}_{n=0}\mathbb{E}\Big[\|u(t_n)\|^4_{H^2}\Big] + 128L_\sigma \sum_{n=0}^{M-1}\int_{t_n}^{t_{n+1}}\mathbb{E}\Big[\|u(t_n) - u(s)\|_{L^{2}}^{4}\Big]\,ds\\
				&\le 4L_\sigma C_P^4kh^8 \left(2+768L_\sigma k\right)\sum_{n=0}^{M-1}\mathbb{E}\Big[\|u(t_n)\|_{H^2}^{4}\Big] +128L_{\sigma}\sum_{n=0}^{M-1}\int_{t_n}^{t_{n+1}}\mathbb{E}\Big[\|u(t_n)-u(s)\|_{L^2}^{4}\Big]\,ds \\
				&\leq Ckh^8(2+Ck)C^{2,4}_{\ref{lem:Hm-q-moment}} + C\sum_{n=0}^{M-1}\int_{t_n}^{t_{n+1}}C^{0,2}_{\ref{lem:time-regularity}}(t_n - s)^2\,ds =: C_{T_7}(h^8 + k^2).
			\end{align*}
			
			The estimate for $T_8$ follows from Lemmas~\ref{lem:Hm-q-moment}, \ref{lem:time-regularity}, and~\ref{lem:stability-2}.
			\begin{align*}
				T_{8} &:= C\sum_{n=0}^{M-1}\mathbb E\Bigg[\Bigg(\int_{t_n}^{t_{n+1}}\|u(t_{n+1})-u(s)\|_{H^1}^4\,ds \Bigg)\|\varepsilon^{n+1}\|^4_{L^2}\Bigg] \\
				&\leq C\mathbb E\Bigg[\Bigg(\sum_{n=0}^{M-1}\int_{t_n}^{t_{n+1}}\|u(t_{n+1})-u(s)\|_{H^1}^4\,ds \Bigg)\sup_{1 \leq \ell \leq M}\|\varepsilon^\ell\|^4_{L^2}\Bigg] \\
				&\leq C\Bigg(\mathbb E\Bigg[\sum_{n=0}^{M-1}\int_{t_n}^{t_{n+1}}\|u(t_{n+1})-u(s)\|_{H^1}^4\,ds \Bigg]^2\Bigg)^{1/2}\Bigg(\mathbb{E}\Bigg[\sup_{1 \leq \ell \leq M}\|\varepsilon^\ell\|^8_{L^2}\Bigg]\Bigg)^{1/2} \\
				&\leq C\Bigg(T\sum_{n=0}^{M-1}\int_{t_n}^{t_{n+1}}\mathbb E\Big[\|u(t_{n+1})-u(s)\|_{H^1}^8\Big]\,ds \Bigg)^{1/2}\Bigg(\mathbb{E}\Bigg[\sup_{1 \leq \ell \leq M}\|\varepsilon^\ell\|^8_{L^2}\Bigg]\Bigg)^{1/2} \\
				&\leq C\Big(C^{1,4}_{\ref{lem:time-regularity}}TMk^5\Big)^{1/2}\Big(C_{\ref{lem:Hm-q-moment}}^{0,8} + C_{\ref{lem:stability-2}}^{0,8}\Big)^{1/2} =: C_{T_8}k^2.
			\end{align*}
			
			By Lemmas~\ref{lem:Hm-q-moment} and~\ref{lem:stability-2}, we estimate \(T_9\) as follows:
			\begin{align*}
				T_{9} &:= 32k^2\mathbb E\Bigg[\sup_{0\leq \ell<M}\left( \frac{40\alpha^2 C_e^2}{\nu k}\int_{t_\ell}^{t_{\ell+1}}\|u(s)\|_{H^1}^{2}\,ds+\frac{20\alpha^2 C_e^2}{\nu}+\frac{15}{2}+\frac{\nu}{8}+2\beta\mu+L_{\sigma}\right)^2\\
				&\quad \quad \times\left(\|\varepsilon^{\ell+1}\|_{L^2}^{2}+\|\varepsilon^\ell\|_{L^2}^{2}\right)^2\Bigg]\\
				&\leq Ck^2 \mathbb{E}\left[\sup_{0\leq \ell<M}\left(1+\frac{1}{k}\int^{t_{\ell+1}}_{t_{\ell}} \|u(s)\|^4_{H^1}\,ds\right)\sup_{0\leq\ell\leq M}\|\varepsilon^{\ell}\|^4_{L^2}\right] \\
				&\leq Ck^2 \left(\mathbb{E}\left[1+\sup_{0 \leq s \leq T}\|u(s)\|^8_{H^1}\right]\right)^{1/2}\left(\mathbb{E}\left[\sup_{0\leq\ell\leq M}\|\varepsilon^{\ell}\|^8_{L^2}\right]\right)^{1/2} \\
				&\leq Ck^2\Big(1 + C^{1,8}_{\ref{lem:Hm-q-moment}}\Big)^{1/2} \Big(C_{\ref{lem:Hm-q-moment}}^{0,8} + C_{\ref{lem:stability-2}}^{0,8}\Big)^{1/2} =: C_{T_9}k^2.
			\end{align*}
			
			Substituting the estimates for $T_1$, $T_2$, $\dots$, $T_{9}$ into \eqref{eq:new-discrete-gronwall-bound}, we obtain
			\begin{equation*}
				\left(\mathbb{E}\left[\sup_{0 \le \ell \le M}\|\varepsilon^\ell\|_{L^2}^{2}\right]\right)^{\frac12}\le \widehat  C_1\left(k^{1/2} + h\right),
			\end{equation*}
			where
			\begin{equation*}    
				\widehat C_1 := \widehat C_0 \sqrt{\max\{C_{T_1},C_{T_2},\dots,C_{T_{9}}\}}.
			\end{equation*}

			It remains to verify that $\{M_\ell\}_{\ell=0}^{M-1}$ is a discrete martingale. The Lipschitz condition on $\sigma$ and the eighth-moment bounds for $u$ and $u_h^n$ imply that $\mathcal M_n$ is integrable. Thus,
			\begin{equation*}    
				\mathbb E[|M_\ell|]\leq\sum_{n=0}^{\ell}\mathbb E[|\mathcal M_n|]<\infty.
			\end{equation*}
			
			Moreover, $\mathcal M_n$ is $\mathcal F_{t_{n+1}}$-measurable. Set $a_n:=\sigma(u(t_n))-\sigma(u_h^n)$ and
			\begin{equation}
				I_n:=\int_{t_n}^{t_{n+1}}\bigl(\sigma(u(s))-\sigma(u_h^n)\bigr)\,dW(s).
			\end{equation}
			
			Since $a_n$ and $\varepsilon^n$ are $\mathcal F_{t_n}$-measurable and $\Delta W_n$ is independent of $\mathcal F_{t_n}$, we have
			\begin{equation*}
				\mathbb E[|\Delta W_n|^2\mid\mathcal F_{t_n}]=k,
				\qquad
				\mathbb E[|\Delta W_n|^4\mid\mathcal F_{t_n}]=3k^2.
			\end{equation*}
			
			The martingale property of the It\^o integral and the conditional It\^o isometry also give
			\begin{equation*}
				\mathbb E[(I_n,\varepsilon^n)\mid\mathcal F_{t_n}]=0,
				\qquad
				\mathbb E[\|I_n\|_{L^2}^2\mid\mathcal F_{t_n}]
				=\mathbb E\left[\int_{t_n}^{t_{n+1}}\|\sigma(u(s))-\sigma(u_h^n)\|_{L^2}^2\,ds\,\middle|\,\mathcal F_{t_n}\right].
			\end{equation*}
			
			Consequently,
			\begin{align*}
				\mathbb E[\mathcal M_n\mid\mathcal F_{t_n}]
				&=8\|a_n\|_{L^2}^2\|\varepsilon^n\|_{L^2}^2\left(\mathbb E[|\Delta W_n|^2\mid\mathcal F_{t_n}]-k\right)\\
				&\quad+4\|\varepsilon^n\|_{L^2}^2\mathbb E[(I_n,\varepsilon^n)\mid\mathcal F_{t_n}]\\
				&\quad+128\|a_n\|_{L^2}^4\left(\mathbb E[|\Delta W_n|^4\mid\mathcal F_{t_n}]-3k^2\right)\\
				&\quad+32\|\varepsilon^n\|_{L^2}^2\mathbb E\left[\|I_n\|_{L^2}^2-\int_{t_n}^{t_{n+1}}\|\sigma(u(s))-\sigma(u_h^n)\|_{L^2}^2\,ds\,\middle|\,\mathcal F_{t_n}\right]=0.
			\end{align*}
			
			Therefore, for $0\leq\ell\leq M-2$,
			\[
			\mathbb E[M_{\ell+1}\mid\mathcal F_{t_{\ell+1}}]
			=M_\ell+\mathbb E[\mathcal M_{\ell+1}\mid\mathcal F_{t_{\ell+1}}]=M_\ell.
			\]
			Thus, $\{M_\ell\}_{\ell=0}^{M-1}$ is a discrete martingale with respect to $\{\mathcal F_{t_{\ell+1}}\}_{\ell=0}^{M-1}$. The proof is complete.

		\end{proof}

		\section{Numerical experiments}\label{sec5}
		In this section, we present numerical tests to examine the convergence rates of scheme \eqref{scheme:implicit-euler}. In all our experiments, we set $D=(0,1)$, $T=0.125=2^{-3}$, $\nu=0.1$, $\alpha=1$, $\beta=1$, and $\gamma=0.5$. We choose $\sigma(u)=0.5u/(1+u^2)$, which is globally Lipschitz and satisfies the linear-growth condition in $L^2(D)$. We use continuous piecewise linear finite elements for the spatial discretization and the implicit Euler--Maruyama scheme for the temporal discretization. To solve the nonlinear system at each time step, we employ the fixed-point iteration in \cite[Algorithm~3]{feng2024high} instead of Newton's method. For all the tests, we use the standard Monte Carlo method with $J=400$ independent samples to compute the expectations.
		
		Since the exact solution is unknown, the errors are computed between the numerical solution $u_h^n(\omega_j)$ and a reference solution $u_{\mathrm{ref}}^n(\omega_j)$ at the $\omega_j$-th sample. The reference solution is computed with $N_{\mathrm{ref}}=64$ spatial subintervals and $M_{\mathrm{ref}}=16384$ time steps, giving $h_{\mathrm{ref}}=2^{-6}$ and $k_{\mathrm{ref}}=2^{-17}$. The real-valued Wiener process is simulated using the minimal time step $k_{\mathrm{ref}}$, with increments $\Delta W_{\mathrm{ref}}^m=\sqrt{k_{\mathrm{ref}}}\,\xi_m$, where $\{\xi_m\}$ are independent standard normal random variables. The increments on coarser grids are obtained by summing the corresponding fine-grid increments, so that both solutions are driven by the same Wiener path. Here, $u_{\mathrm{ref}}^n$ denotes the reference solution evaluated at $t_n$. To compute the errors in strong norms, we use the following numerical integration formulas:
		\begin{align*}
			\|u-u_h\|_{L_\omega^2\ell_t^\infty L_x^2}
			&\approx\left[\frac1J\sum_{j=1}^{J}\max_{0\leq n\leq M}\|u_{\mathrm{ref}}^n(\omega_j)-u_h^n(\omega_j)\|_{L^2}^2\right]^{1/2},\\
			\|u-u_h\|_{L_\omega^2\ell_t^2H_x^1}
			&\approx\left[\frac1J\sum_{j=1}^{J}k\sum_{n=1}^{M}\|u_{\mathrm{ref}}^n(\omega_j)-u_h^n(\omega_j)\|_{H^1}^2\right]^{1/2}.
		\end{align*}
		
		\subsection*{Test 1: Periodic boundary conditions}
		In the first test, we examine the convergence rates under periodic boundary conditions and choose $u_0(x)=0.2+0.2\sin(2\pi x)$. To compute the $L_\omega^2\ell_t^\infty L_x^2$-errors, we take $h=2^{-i}$, $i=1,2,3,4$, and choose $k=h^4$, so that $k^{1/2}=h^2$. For the $L_\omega^2\ell_t^2H_x^1$-errors, we take $h=2^{-i}$, $i=2,3,4,5$, and choose $k=h^2$, so that $k^{1/2}=h$. The errors and convergence rates are reported in Tables~\ref{tab:periodic-coupled-linf-l2} and~\ref{tab:periodic-coupled-l2-h1}. We observe approximately second-order convergence in $h$ in the $L_\omega^2\ell_t^\infty L_x^2$-norm, which is higher than the first-order spatial rate established in our analysis. In addition, the results show first-order convergence in $h$ in the $L_\omega^2\ell_t^2H_x^1$-norm, consistent with the optimal spatial order for piecewise linear finite elements. These rates are obtained along the coupled refinements of $h$ and $k$.
		
		\begin{table}[H]
			\centering
			\small
			\setlength{\tabcolsep}{4pt}
			\begin{tabular}{@{}ccccc@{}}
				\toprule
				$h$ & $k=h^4$ & $\|u-u_h\|_{L_\omega^2\ell_t^\infty L_x^2}$ & Rate in $h$ & Rate in $k$ \\
				\midrule
				$2^{-1}$ & $2^{-4}$  & $1.4132604\times10^{-1}$ & --       & --       \\
				$2^{-2}$ & $2^{-8}$  & $3.130046\times10^{-2}$  & $2.1748$ & $0.5437$ \\
				$2^{-3}$ & $2^{-12}$ & $8.241740\times10^{-3}$  & $1.9252$ & $0.4813$ \\
				$2^{-4}$ & $2^{-16}$ & $2.007520\times10^{-3}$  & $2.0375$ & $0.5094$ \\
				\bottomrule
			\end{tabular}
			\caption{Periodic boundary conditions: $L_\omega^2\ell_t^\infty L_x^2$-errors and convergence rates along the coupled refinement $k=h^4$.}
			\label{tab:periodic-coupled-linf-l2}
		\end{table}
		
		\begin{table}[H]
			\centering
			\small
			\setlength{\tabcolsep}{4pt}
			\begin{tabular}{@{}ccccc@{}}
				\toprule
				$h$ & $k=h^2$ & $\|u-u_h\|_{L_\omega^2\ell_t^2H_x^1}$ & Rate in $h$ & Rate in $k$ \\
				\midrule
				$2^{-2}$ & $2^{-4}$  & $9.715564\times10^{-2}$ & --       & --       \\
				$2^{-3}$ & $2^{-6}$  & $5.442351\times10^{-2}$ & $0.8361$ & $0.4180$ \\
				$2^{-4}$ & $2^{-8}$  & $2.735688\times10^{-2}$ & $0.9923$ & $0.4962$ \\
				$2^{-5}$ & $2^{-10}$ & $1.232525\times10^{-2}$ & $1.1503$ & $0.5751$ \\
				\bottomrule
			\end{tabular}
			\caption{Periodic boundary conditions: $L_\omega^2\ell_t^2H_x^1$-errors and convergence rates along the coupled refinement $k=h^2$.}
			\label{tab:periodic-coupled-l2-h1}
		\end{table}

		\subsection*{Test 2: Homogeneous Dirichlet boundary conditions}
		In this test, we examine the convergence rates under homogeneous Dirichlet boundary conditions $u(t,0)=u(t,1)=0$ and choose $u_0(x)=0.2\sin(2\pi x)$. We follow the same strategies as in Test~1 to compute the errors, using the same numerical parameters, reference discretization, and Monte Carlo sample paths. The results are reported in Tables~\ref{tab:dirichlet-coupled-linf-l2} and~\ref{tab:dirichlet-coupled-l2-h1}. We observe similar results as in Test~1, with approximately second-order convergence in $h$ in the $L_\omega^2\ell_t^\infty L_x^2$-norm and first-order convergence in $h$ in the $L_\omega^2\ell_t^2H_x^1$-norm along the coupled refinements. Although our error analysis concerns periodic boundary conditions, this test suggests similar convergence behavior under homogeneous Dirichlet boundary conditions.
		
		\begin{table}[H]
			\centering
			\small
			\setlength{\tabcolsep}{4pt}
			\begin{tabular}{@{}ccccc@{}}
				\toprule
				$h$ & $k=h^4$ & $\|u-u_h\|_{L_\omega^2\ell_t^\infty L_x^2}$ & Rate in $h$ & Rate in $k$ \\
				\midrule
				$2^{-1}$ & $2^{-4}$  & $3.150835\times10^{-2}$ & --       & --       \\
				$2^{-2}$ & $2^{-8}$  & $8.660910\times10^{-3}$ & $1.8631$ & $0.4658$ \\
				$2^{-3}$ & $2^{-12}$ & $2.192550\times10^{-3}$ & $1.9819$ & $0.4955$ \\
				$2^{-4}$ & $2^{-16}$ & $5.287100\times10^{-4}$ & $2.0521$ & $0.5130$ \\
				\bottomrule
			\end{tabular}
			\caption{Homogeneous Dirichlet boundary conditions: $L_\omega^2\ell_t^\infty L_x^2$-errors and convergence rates along the coupled refinement $k=h^4$.}
			\label{tab:dirichlet-coupled-linf-l2}
		\end{table}
		
		\begin{table}[H]
			\centering
			\small
			\setlength{\tabcolsep}{4pt}
			\begin{tabular}{@{}ccccc@{}}
				\toprule
				$h$ & $k=h^2$ & $\|u-u_h\|_{L_\omega^2\ell_t^2H_x^1}$ & Rate in $h$ & Rate in $k$ \\
				\midrule
				$2^{-2}$ & $2^{-4}$  & $3.157355\times10^{-2}$ & --       & --       \\
				$2^{-3}$ & $2^{-6}$  & $1.631271\times10^{-2}$ & $0.9527$ & $0.4764$ \\
				$2^{-4}$ & $2^{-8}$  & $8.027440\times10^{-3}$ & $1.0230$ & $0.5115$ \\
				$2^{-5}$ & $2^{-10}$ & $3.598510\times10^{-3}$ & $1.1575$ & $0.5788$ \\
				\bottomrule
			\end{tabular}
			\caption{Homogeneous Dirichlet boundary conditions: $L_\omega^2\ell_t^2H_x^1$-errors and convergence rates along the coupled refinement $k=h^2$.}
			\label{tab:dirichlet-coupled-l2-h1}
		\end{table}

	\subsection*{Test 3: Influence of the stochastic forcing}
In the last test, we examine the influence of the noise using the same initial condition, periodic boundary conditions, and model parameters as in Test~1. We take $N=128$ and $M=1024$, giving $h=2^{-7}$ and $k=2^{-13}$. We compare the Monte Carlo approximation of $\mathbb E[u_h^{\mathrm{sto}}]$ for $\sigma(u)=0.5u/(1+u^2)$ with the deterministic solution obtained by setting $\sigma=0$ as well as the sample-path solutions.

\begin{figure}[H]
	\centering
	\includegraphics[width=\textwidth]{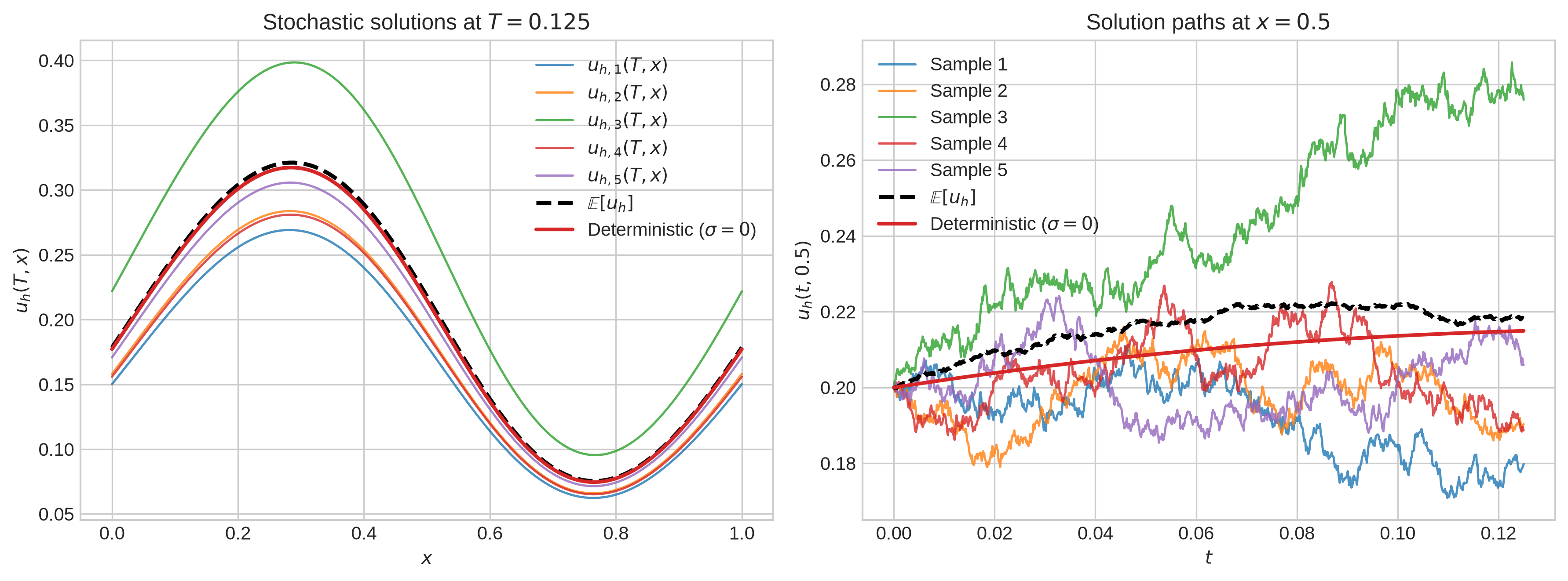}
	\caption{Influence of noise on the numerical solution under periodic boundary conditions. Left: the first five sample solutions, the Monte Carlo approximation of the expected solution, and the deterministic solution at $T=0.125$ as functions of $x$. Right: the corresponding solution paths at $x=0.5$ as functions of time.}
	\label{fig:noise-deterministic-comparison}
\end{figure}

Figure~\ref{fig:noise-deterministic-comparison} shows the first five sample solutions, the expected solution, and the deterministic solution at $T=0.125$ as functions of $x$ in the left plot and at $x=0.5$ as functions of time in the right plot. The sample solutions differ noticeably from the expected and deterministic solutions, while the latter two have similar spatial profiles. The sample paths also exhibit the irregular temporal behavior expected for solutions driven by a Wiener process.

		\section{Conclusion}\label{sec6}
		In this paper, we developed and analyzed a fully discrete finite element method for the stochastic Burgers--Huxley equation with multiplicative It\^o noise. Using moment bounds, exponential integrability estimates, and the discrete stochastic Gronwall inequality, we established sub-second-moment error estimates in the discrete maximum-in-time $L^2(D)$-norm and the discrete $L^2$-in-time $H^1(D)$-norm. A higher-moment argument then yielded second-moment error estimates in the discrete maximum-in-time $L^2(D)$-norm. The rates are of order $k^{1/2}$ in time and $h$ in space, with the spatial rate being optimal in the $H^1(D)$-norm. Numerical experiments showed first-order convergence in the $H^1$-based norm and second-order convergence in the $L^2$-based norm along coupled refinements, and illustrated the influence of multiplicative noise. Future work will focus on improving the spatial $L^2(D)$-error estimate, developing higher-order time discretizations, and extending the analysis to other nonlinear stochastic reaction--diffusion equations.

		\section*{Declarations of Funding}\,
		
		The author Liet Vo was supported by the National Science Foundation (NSF) under Grant No. DMS-2530211.
		
		\section*{Data Availability} \,
		
		Data sharing is not applicable to this article as no datasets were generated or analysed during the current study.
		
		\section*{Declaration of Conflict of Interest}\,
		
		The authors have no conflict of interest.

		\appendix
		\section{Auxiliary Results}\label{AppendixA}
		
		First, we state the following discrete stochastic Gronwall inequality from
		\cite[Theorem~1]{kruse2018discrete}, which plays a vital role in our error analysis in Section~3.
		\begin{lemma}\label{lem:discrete-stochastic-gronwall}
			Let \(\{M_n\}_{n\in\mathbb N}\) be an \(\{\mathcal F_n\}_{n\in\mathbb N}\)-martingale satisfying \(M_0=0\) on a filtered probability space \((\Omega,\mathcal F,\{\mathcal F_n\}_{n\in\mathbb N},\mathbb P)\). Let \(\{X_n\}_{n\in\mathbb N}\), \(\{F_n\}_{n\in\mathbb N}\), and \(\{G_n\}_{n\in\mathbb N}\) be sequences of nonnegative and adapted random variables with \(\mathbb E[X_0]<\infty\) such that
			\begin{equation*}
				X_n \leq F_n+M_n+\sum_{\ell=0}^{n-1}G_\ell X_\ell, \qquad \text{for all } n\in\mathbb N. \label{eq:discrete-stochastic-gronwall-assumption}
			\end{equation*}
			
			Then, for any \(q\in(0,1)\) and a pair of conjugate numbers \(a,b\in[1,\infty]\), i.e.,
			\[
			\frac1a+\frac1b=1,
			\]
			satisfying \(qa<1\), there holds
			\begin{equation*}
				\mathbb E\left[\sup_{0\leq \ell\leq n}X_\ell^q\right]\leq\left(1+\frac{1}{1-a q}\right)^{\frac1a}\left\|\prod_{\ell=0}^{n-1}(1+G_\ell)^q\right\|_{L^b(\Omega)}\left(\mathbb E\left[\sup_{0\leq \ell\leq n}F_\ell\right]\right)^q. \label{eq:discrete-stochastic-gronwall}
			\end{equation*}
		\end{lemma}

		\section{Theoretical PDE result proofs}\label{AppendixB}
		
		\subsection{Proof of Lemma \ref{lem:exponential-L2-estimate}}\label{AppendixB1}
		\begin{proof}
			
			Applying It\^o's formula to $\|u(t)\|^2_{L^2}$, we obtain    
			\begin{align}\label{eq:lem2.2-equation01}
				&	\|u(t)\|_{L^2}^2 + 2\nu \int_0^t \|\partial_x u(s)\|_{L^2}^2 \, ds + 2\beta\gamma \int_0^t \|u(s)\|_{L^2}^2 \, ds + 2\beta \int_0^t \|u(s)\|_{L^4}^4 \, ds \\\nonumber
				&= \|u_0\|_{L^2}^2 + 2\beta(1+\gamma)\int_0^t \|u(s)\|_{L^3}^3 \, ds  + 2 \int_0^t  \bigl(u(s), \sigma(u(s))\bigr)\, dW(s) + \int_0^t \|\sigma(u(s))\|_{L^2}^2 \, ds .
			\end{align}
			
			From Young's inequality,
			\begin{align*}
				2 \beta(1+ \gamma) \int^t_0 \|u(s)\|^3_{L^3}\,ds \leq \beta \int^t_0 \|u(s)\|^4_{L^4}\,ds + \beta(1+\gamma)^2 \int^t_0 \|u(s)\|^2_{L^2}\,ds.
			\end{align*}
			
			Combining \eqref{eq:lem2.2-equation01} with assumption \eqref{ass:noise-linear-growth}, we obtain
			\begin{align}\label{eq:lem2.2-young-apply}
				&\|u(t)\|_{L^2}^2 + 2\nu \int_0^t \|\partial_x u(s)\|_{L^2}^2 \, ds + 2\beta\gamma \int_0^t \|u(s)\|_{L^2}^2 \, ds + 2\beta \int_0^t \|u(s)\|_{L^4}^4 \, ds \\\nonumber
				&\le \|u_0\|_{L^2}^2 + \beta \int_0^t \|u(s)\|_{L^4}^4 \, ds + \bigl(\beta(1+\gamma)^2 + K_\sigma\bigr)\int_0^t \|u(s)\|_{L^2}^2 \, ds \\\nonumber
				&\qquad + tK_{\sigma} + 2 \int_0^t \bigl(u(s), \sigma(u(s))\bigr)\, dW(s).
			\end{align}
			
			Applying Young's inequality again, we have
			\begin{align*}
				\bigl(\beta(1+ \gamma)^2 + K_{\sigma}\bigr)
				\int^t_0 \|u(s)\|^2_{L^2}\,ds \leq \int^t_0 
				\left[ \int_D \frac{(\beta(1+\gamma)^2 + K_{\sigma})^2}{\beta}\,dx + \frac{\beta}{4} \|u(s)\|^4_{L^4} \right]ds.
			\end{align*}
			
			Multiplying both sides of \eqref{eq:lem2.2-young-apply} by \(\kappa>0\), we obtain
			\begin{align}\label{eq:lem2.2-estimate}
				&\kappa \|u(t)\|^2_{L^2} + 2\nu \kappa \int^t_0 \|\partial_x u(s)\|^2_{L^2}\,ds + 2\beta \gamma\kappa \int^t_0 \|u(s)\|^2_{L^2} \,ds + \frac{3}{4}\beta\kappa  \int^t_0 \|u(s)\|^4_{L^4}\, ds \\\nonumber
				&\leq \kappa \|u_0\|^2_{L^2} + \left(\frac{(\beta(1+\gamma)^2 + K_{\sigma})^2L}{\beta} + K_{\sigma}\right)T\kappa  +  2\kappa \int^t_0 \big(u(s), \sigma(u(s))\big) \,dW(s).
			\end{align}
			
			Define $Z_t := 2 \kappa \int^t_0 \bigl(u(s), \sigma(u(s))\bigr)\,dW(s)$. By assumption \eqref{ass:noise-linear-growth}, $Z_t$ is a martingale and its quadratic variation satisfies
			\begin{align*}
				\langle Z \rangle_t &= 4 \kappa^2 \int^t_0 |(u(s),\sigma(u(s))|^2 \,ds \\
				&\leq 4 \kappa^2 \int^t_0 \|u(s)\|^2_{L^2} \cdot \|\sigma(u(s))\|^2_{L^2} \,ds \\
				&\leq 4 \kappa^2 K_{\sigma} \int^t_0 \bigl(\|u(s)\|_{L^2}^2 + \|u(s)\|^4_{L^2} \bigr)\, ds \\
				&\leq 4 \kappa^2 K_{\sigma} \int^t_0 \left(\frac{3}{2}\|u(s)\|^4_{L^2} + \frac{1}{2} \right) ds \\
				&\leq 6 \kappa^2 K_{\sigma} L \int^t_0 \|u(s)\|^4_{L^4} \,ds + 2\kappa^2 K_{\sigma}T.
			\end{align*}
			
			Adding and subtracting $\langle Z \rangle_t$, we obtain
			\begin{align*}
				\kappa \|u(t)\|_{L^2}^2 &+ 2\nu\kappa \int_0^t \|\partial_x u(s)\|_{L^2}^2\,ds + 2\beta\gamma\kappa \int_0^t \|u(s)\|_{L^2}^2\,ds + \frac34 \beta\kappa \int_0^t \|u(s)\|_{L^4}^4\,ds \\
				&\le \kappa \|u_0\|_{L^2}^2 + \frac{\bigl(\beta(1+\gamma)^2 + K_{\sigma}\bigr)^2}{\beta}LT\kappa + \bigl(Z_t-\langle Z\rangle_t\bigr) + \langle Z\rangle_t \\
				&\le \kappa \|u_0\|_{L^2}^2 + \frac{\bigl(\beta(1+\gamma)^2 + K_{\sigma}\bigr)^2}{\beta}LT\kappa + 2 \kappa^2 K_{\sigma}T +\bigl(Z_t-\langle Z\rangle_t\bigr) \\
				&\qquad + 6\kappa^2 K_{\sigma}L \int_0^t \|u(s)\|_{L^4}^4\,ds .
			\end{align*}
			
			Therefore, 
			\begin{equation*}
				\begin{aligned}
					&\kappa \|u(t)\|^2_{L^2} + 2\nu \kappa \int^t_0 \|\partial_x u(s)\|^2_{L^2} \,ds + 2\beta \gamma\kappa \int^t_0 \|u(s)\|^2_{L^2} \,ds +\kappa \left(\frac{3}{4} \beta - 6 \kappa K_{\sigma} L\right)\int^t_0 \|u(s)\|^4_{L^4} \,ds \\
					&\leq \kappa \|u_0\|^2_{L^2} + \frac{(\beta(1+\gamma)^2 + K_{\sigma})^2}{\beta} LT\kappa + 2 \kappa^2 K_{\sigma}T +(Z_t - \langle Z \rangle_t).
				\end{aligned}
			\end{equation*}
			
			Provided that $\kappa$ is sufficiently small, for instance, $0 < \kappa \leq \frac{\beta}{8 K_{\sigma} L}$. We take expectations of the exponential and invoke H\"older's inequality to deduce
			\begin{align*}
				&\mathbb{E}\left(\exp\left[\kappa \|u(t)\|_{L^2}^{2}+ 2\nu\kappa \int_{0}^{t} \|\partial_x u(s)\|_{L^2}^{2}\,ds + 2\beta \gamma\kappa \int^t_0 \|u(s)\|^2_{L^2} \,ds\right]\right)\\
				&\le \exp\left(\frac{(\beta(1+\gamma)^2+K_\sigma)^2}{\beta}LT\kappa+2\kappa^2 K_\sigma T\right) \left(\mathbb{E}\left[\exp\left(2\kappa\|u_0\|_{L^2}^{2}\right)\right]\right)^{1/2} \\
				&\qquad\times\left(\mathbb{E}\left[\exp\left(2Z_t-\frac{\langle 2Z\rangle_t}{2}\right)\right]\right)^{1/2}, \qquad \text{for all } t\in[0,T].
			\end{align*}
			Moreover, for all \(t\in[0,T]\),
			\[
			\mathbb{E}\left[\exp\left(2Z_t-\frac{\langle 2Z\rangle_t}{2}\right)\right]\le 1.
			\]
			
			This completes the proof.
		\end{proof}
		
		\subsection{Proof of Lemma \ref{lem:Hm-q-moment}}\label{AppendixB2}
		\begin{proof}
			
			We prove \eqref{eq:Hm-moment-estimate} and \eqref{eq:Hm-sup-moment-estimate} by induction on the regularity index \(m\). We first consider the base case \(m=0\). Applying It\^o's formula to \(\|u\|_{L^2}^{q}\), we obtain
			\begin{align*}
				d\|u\|_{L^2}^q + \nu q \|u\|_{L^2}^{q-2}\|\partial_xu\|_{L^2}^2\,dt &= -\beta q \gamma \|u\|_{L^2}^{q}\,dt +(1+\gamma)\beta q \|u\|_{L^2}^{q-2}\|u\|_{L^3}^3\,dt \\
				&\qquad -\beta q \|u\|^{q-2}_{L^2}\|u\|^4_{L^4}\,dt + q\|u\|_{L^2}^{q-2}(u,\sigma(u))\,dW(t)\\
				&\qquad + \frac{q}{2}\|u\|_{L^2}^{q-2}\|\sigma(u)\|_{L^2}^2\,dt + \frac{q(q-2)}{2}\|u\|_{L^2}^{q-4}(u,\sigma(u))^2\,dt.
			\end{align*}
			
			Using Young's inequality, we estimate
			\begin{align*}
				(1+\gamma)\beta q\|u\|_{L^2}^{q-2} \|u\|^3_{L^3} &= \beta q \|u\|^{q-2}_{L^2} \|u\|^3_{L^3} + \gamma \beta q \|u\|^{q-2}_{L^2} \|u\|^3_{L^3} \\
				&\leq \frac{\beta q}{2} \|u\|^{q-2}_{L^2}(\|u\|^2_{L^2} + \|u\|^4_{L^4}) + \frac{\beta q}{2} \|u\|^{q-2}_{L^2} (\gamma^2\|u\|^2_{L^2} + \|u\|^4_{L^4}) \\
				&= \frac{\beta q (1+\gamma^2)}{2} \|u\|^q_{L^2} + \beta q \|u\|^{q-2}_{L^2} \|u\|^4_{L^4}.
			\end{align*}
			
			Using assumption \eqref{ass:noise-linear-growth} and the inequality $\|u\|^{q-2}_{L^2} \leq \frac{q-2}{q}\|u\|^q_{L^2} + \frac{2}{q}$, we obtain the estimates
			\begin{align*}
				\frac{q}{2}\|u\|^{q-2}_{L^2}\|\sigma(u)\|^2_{L^2} \leq \frac{qK_\sigma}{2}\bigl(\|u\|_{L^2}^{q-2}+\|u\|_{L^2}^{q}\bigr) \le (q-1)K_\sigma\bigl(1+\|u\|_{L^2}^{q}\bigr)
			\end{align*} 
			and
			\begin{align*}
				\frac{q(q-2)}{2}\|u\|^{q-4}_{L^2} \left|(u,\sigma(u))\right|^2 \leq  \frac{q(q-2)K_{\sigma} }{2}\bigl(\|u\|^{q-2}_{L^2} + \|u\|^q_{L^2}\bigr) \leq (q-1)(q-2)K_{\sigma}\bigl(1 + \|u\|^q_{L^2}\bigr).
			\end{align*}
			
			Combining the above estimates yields
			\begin{align}\label{eq:lem2.3-estimate01}
				d\|u\|_{L^2}^q + \nu q \|u\|_{L^2}^{q-2}\|\partial_x u\|_{L^2}^2\,dt \le q\|u\|_{L^2}^{q-2}(u,\sigma(u))\,dW(t) + C(1 + \|u\|_{L^2}^q)\,dt,
			\end{align}
			for a positive constant $C = \left(\frac{\beta q(1+\gamma^2)}{2} + (q-1)^2K_{\sigma}\right)$. Integrating over $[0,t]$ and taking expectations, we obtain, by Fubini's theorem,
			\begin{align*}
				\mathbb{E}\big[\|u(t)\|_{L^2}^q\big] + \nu q \int_0^t  \mathbb{E}\big[\|u(s)\|_{L^2}^{q-2}\|\partial_x u(s)\|_{L^2}^2\big] \, ds \le  \mathbb{E}\big[\|u_0\|_{L^2}^q\big] + C \int_0^t \bigl(1 + \mathbb{E}\big[\|u(s)\|_{L^2}^q\big] \bigr) \, ds,
			\end{align*}
			where we used that the stochastic integral
			\[
			M_0(t) := q\int_0^t \|u(s)\|_{L^2}^{q-2}(u(s),\sigma(u(s)))\,dW(s)
			\]
			is a martingale with zero expectation. By Gronwall's inequality,        
			\begin{align}\label{eq:lem2.3-L2-moment}
				\mathbb{E}\big[\|u(t)\|_{L^2}^q\big] \le \bigl(1 + \mathbb{E}\big[\|u_0\|_{L^2}^q\big]\bigr)e^{Ct}, \quad t \ge 0.
			\end{align}
			
			This establishes \eqref{eq:Hm-moment-estimate} for $m = 0$. Now, from \eqref{eq:lem2.3-estimate01}, we observe that
			\begin{align}\label{eq:u-sup-moment-pre-bdg}
				\mathbb{E}\!\left[\sup_{s \in [0,T]} \|u(s)\|_{L^2}^q \right] &\le \mathbb{E}\big[\|u_0\|_{L^2}^q\big] + \mathbb{E}\!\left[\sup_{s \in [0,T]} \left|M_0(s)\right| \right] + C \int_0^T \bigl(1 + \mathbb{E}\big[\|u(s)\|_{L^2}^q\big] \bigr)\, ds.
			\end{align}
			
			Using assumption \eqref{ass:noise-linear-growth}, we estimate the quadratic variation of $M_0$
			\begin{equation}\label{eq:quadratic-variation-M0-estimate}
				d\langle M_0 \rangle_s = q^2 \|u(s)\|_{L^2}^{2q-4} \left|(u(s), \sigma(u(s)))\right|^2\, ds \le C\|u(s)\|_{L^2}^{2q-2} \big(1 + \|u(s)\|^2_{L^2}\big)\, ds.
			\end{equation}
			
			We invoke Burkh\"older--Davis--Gundy's inequality (see \cite[Theorem 1]{davis1970intergrability} and \cite[Theorem 1.1]{burkholder2002best}) to infer
			\begin{align}\label{eq:M0-BDG-estimate}
				\mathbb{E}\!\left[\sup_{s \in [0,T]} \left|M_0(s)\right| \right] &\le C\mathbb{E}\left[ \langle M_0\rangle_T^{1/2}\right] \\\nonumber
				&\le C\, \mathbb{E}\!\left[\left( \int_0^T \|u(s)\|_{L^2}^{2q-2} \big(1 + \|u(s)\|^2_{L^2}\big)\, ds \right)^{\frac{1}{2}} \right] \\\nonumber
				&\le C\, \mathbb{E}\!\left[ \sup_{s \in [0,T]} \|u(s)\|_{L^2}^{q-1}\left(\int_0^T \big(1 + \|u(s)\|^2_{L^2}\big)\, ds \right)^{\frac{1}{2}} \right] \\\nonumber
				&\le \frac{1}{10} \,\mathbb{E}\! \left[ \sup_{s \in [0,T]} \|u(s)\|_{L^2}^{q}\right] +C\,\mathbb{E}\!\left[\left(\int_0^T \big(1 + \|u(s)\|^2_{L^2}\big)\, ds \right)^{\frac{q}{2}} \right] \\\nonumber
				&\leq \frac{1}{10} \,\mathbb{E}\! \left[ \sup_{s \in [0,T]} \|u(s)\|_{L^2}^{q}\right]             + C \int_0^T \bigl(1 + \mathbb{E}\big[\|u(s)\|_{L^2}^q\big]\bigr)\,ds.
			\end{align}
			
			Applying this estimate to \eqref{eq:u-sup-moment-pre-bdg} and using the estimate \eqref{eq:lem2.3-L2-moment}, we obtain
			\begin{equation*}
				\begin{aligned}                     
					\,\mathbb{E}\!\left[\sup_{s\in[0,T]}\|u(s)\|_{L^2}^q\right] &\le \mathbb{E} \big[\|u_0\|_{L^2}^q\big] + C \int_0^T \bigl(1 + \mathbb{E}\big[\|u(s)\|_{L^2}^q\big]\bigr)\,ds \\
					&\le \mathbb{E} \big[\|u_0\|_{L^2}^q\big] + C \bigl(1 + \mathbb{E} \big[\|u_0\|_{L^2}^q\big]\bigr) \\
					&\le C\bigl(1 + \mathbb{E}\big[\|u_0\|_{L^2}^{q}\big]\bigr).
				\end{aligned}
			\end{equation*} 
			
			This proves \eqref{eq:Hm-sup-moment-estimate} and completes the base case \(m=0\) for every \(q\geq 2\). We now proceed with the induction step. Assume that \eqref{eq:Hm-moment-estimate} and \eqref{eq:Hm-sup-moment-estimate} hold for all regularity levels \(0,\dots,m-1\), where \(m\geq 1\), and for every \(q\geq 2\). Applying It\^o's formula to \(\|\partial_x^m u\|_{L^2}^{q}\), we obtain
			\begin{align}\label{eqn:partial-m-u}
				d\|\partial_x^m u\|_{L^2}^q + \nu q \|\partial_x^m u\|_{L^2}^{q-2} \|\partial_x^{m+1} u\|_{L^2}^2 \, dt &=  -\alpha q \|\partial_x^m u\|_{L^2}^{q-2} \bigl(\partial_x^{m} (u \partial_x u),\partial_x^{m} u\bigr) \, dt \\\nonumber
				&\quad - \beta q \|\partial_x^m u\|_{L^2}^{q-2} \bigl(\partial_x^m u, \partial_x^m (u^3)\bigr) \, dt  \\\nonumber
				&\quad+ (1+\gamma)\beta q \|\partial_x^m u\|_{L^2}^{q-2}\bigl(\partial_x^m u, \partial_x^m (u^2)\bigr) \, dt \\\nonumber
				&\quad - \gamma \beta q \|\partial_x^m u\|_{L^2}^{q} \, dt  + q \|\partial_x^m u\|_{L^2}^{q-2} \bigl(\partial_x^m u, \partial_x^m \sigma(u)\bigr) \, dW(t) \\\nonumber
				&\quad + \frac{q}{2} \|\partial_x^m u\|_{L^2}^{q-2} \|\partial_x^m \sigma(u)\|_{L^2}^2 \, dt  \\\nonumber
				&\quad+ \frac{q(q-2)}{2} \|\partial_x^m u\|_{L^2}^{q-4}\bigl(\partial_x^m u, \partial_x^m \sigma(u)\bigr)^2 \, dt.
			\end{align}
			
			By assumption \eqref{ass:noise-linear-growth} and the same argument as in case $m = 0$, we obtain
			\begin{align*}
				\frac{q}{2} \|\partial_x^m u\|_{L^2}^{q-2} \|\partial_x^m \sigma(u)\|_{L^2}^2+ \frac{q(q-2)}{2}\|\partial_x^m u\|_{L^2}^{q-4} \bigl(\partial_x^m u,\partial_x^m \sigma(u)\bigr)^2\le (q-1)^2K_{\sigma}\big(1 +\|\partial_x^m u\|_{L^2}^q\big).
			\end{align*}
			
			For the first term in the right-hand side of \eqref{eqn:partial-m-u}, by Leibniz's rule and integration by parts, we obtain
			\begin{align*}
				\left|\bigl(\partial_x^{m} (u \partial_x u),\partial_x^{m} u\bigr)\right| &\leq \sum_{j=0}^{m} \binom{m}{j} \left|\int_D (\partial_x^j u)(\partial_x^{m-j+1}u)(\partial_x^m u)\,dx \right|\\
				&\leq C\left|\int_D (\partial_x u)(\partial_x^{m}u)^2\,dx\right| +  \sum_{j=2}^{m-1} \binom{m}{j} \left| \int_D (\partial_x^j u)(\partial_x^{m-j+1}u)(\partial_x^m u)\,dx \right|\\
				&=: T_1 + T_2.
			\end{align*}
			
			We notice that if \(m\leq 2\), then \(T_2=0\). We first estimate the term involving \(T_1\) as follows:        
			\begin{align*}
				\alpha q \|\partial_x^m u\|^{q-2}_{L^2} T_1 &\leq C\|\partial_x^m u\|^q_{L^2} \|\partial_x u\|_{L^\infty} \\
				&\leq C\|u\|^2_{H^{m-1}}\|\partial_x^m u\|^{q-2}_{L^2}\|\partial_x^{m+1}u\|_{L^2} \\
				&\leq \frac{\nu q}{4} \|\partial_x^m u\|^{q-2}_{L^2}\|\partial_x^{m+1} u\|^2_{L^2} + C\|\partial_x^m u\|^q_{L^2} + C \|u\|^{2q}_{H^{m-1}}.
			\end{align*}
			
			For each \(2\leq j\leq m-1\), set
			\[
			\ell_j:=\min\{j,m-j+1\},\qquad r_j:=\max\{j,m-j+1\}.
			\]
			
			Then \(\ell_j\leq m-2\) and \(r_j\leq m-1\). Therefore, by H\"older's inequality and the one-dimensional Sobolev embedding \eqref{ass:sobolev-embedding},
			\begin{align}\label{eq:estimate01}
				T_2 \leq \sum_{j=2}^{m-1} \binom{m}{j} \|\partial_x^{r_j} u\|_{L^{\infty}}\|\partial_x^{\ell_j}u\|_{L^2}\|\partial_x^m u\|_{L^2} \leq C\|u\|^2_{H^{m-1}}\|\partial_x^m u\|_{L^2}.
			\end{align}
			
			Consequently,
			\begin{align*}
				\alpha q \|\partial_x^m u\|^{q-2}_{L^2} T_2 &\leq C\|u\|^2_{H^{m-1}}\|\partial_x^m u\|_{L^2}^{q-1} \leq C \|\partial_x^{m} u\|^q_{L^2} + C\|u\|^{2q}_{H^{m-1}}.
			\end{align*}
			
			Combining the estimates for \(T_1\) and \(T_2\), we conclude that
			\begin{align*}
				\alpha q \|\partial_x^m u\|_{L^2}^{q-2} \left|\bigl(\partial_x^{m} (u \partial_x u),\partial_x^{m} u\bigr)\right| \leq \frac{\nu q}{4} \|\partial_x^m u\|^{q-2}_{L^2}\|\partial_x^{m+1} u\|^2_{L^2} + C\|\partial_x^m u\|^q_{L^2} + C \|u\|^{2q}_{H^{m-1}}.
			\end{align*}
			
			For the cubic reaction term, we have
			\begin{align*}
				&-\beta q \|\partial^m_x u \|^{q-2}_{L^2} \bigl(\partial^m_x u,\partial^m_x (u^3)\bigr) \\
				=& -3\beta q\|\partial_x^m u\|_{L^2}^{q-2} \int_D u^2|\partial_x^m u|^2\,dx -\beta q\|\partial_x^m u\|_{L^2}^{q-2} \left(\partial_x^m u,\partial^m_x(u^3) - 3u^2\partial^m_x u\right) \\
				\leq&\, \beta q\|\partial_x^m u\|_{L^2}^{q-2}\left| \left(\partial_x^m u,\partial^m_x(u^3) - 3u^2\partial^m_x u\right) \right| \\
				\leq&\, \beta q\|\partial_x^m u\|_{L^2}^{q-1}\|\partial^m_x(u^3) - 3u^2\partial^m_x u\|_{L^2} \\
				\leq&\, C \|\partial_x^m u\|_{L^2}^{q-1} \|u\|_{H^{m-1}}^3\\
				\leq&\, C\|\partial^m_x u\|^q_{L^2} + C\|u\|^{3q}_{H^{m-1}}.
			\end{align*}
			
			Indeed, 
			\begin{align}\label{eq:estimate-02}
				\partial^m_x(u^3) - 3u^2\partial^m_x u &= \sum_{i + j + k = m} \frac{m!}{i!j!k!} (\partial_x^i u)(\partial_x^j u)(\partial_x^k u) - 3u^2\partial^m_x u \\
				&= \sum_{\substack{i+j+k=m\\ 0\le i,j,k\le m-1}} \frac{m!}{i!j!k!} (\partial_x^i u)(\partial_x^j u)(\partial_x^k u). 
			\end{align}
			
			Hence, by the one-dimensional Sobolev embedding and the same argument as in \eqref{eq:estimate01}, we imply
			\[
			\left\|\partial_x^m(u^3)-3u^2\partial_x^m u\right\|_{L^2} \le C\|u\|_{H^{m-1}}^3.
			\]
			
			Similarly, for the quadratic reaction term, we yield
			\begin{align*}
				(1+\gamma)\beta q \|\partial_x^m u\|_{L^2}^{q-2}\bigl(\partial_x^m u, \partial_x^m (u^2)\bigr) &= \beta(1+\gamma) q \|\partial_x^m u\|^{q-2}_{L^2}\int_D u \left|\partial^m_x u\right|^2dx\\
				&\quad+\beta(1+\gamma) q \|\partial_x^m u\|^{q-2}_{L^2}\big(\partial^m_x u, \partial^m_x (u^2) - 2u\partial^m_x u\big).
			\end{align*}
			
			We estimate the two terms on the right-hand side separately. This gives
			\begin{align*}
				\beta(1+\gamma) q \|\partial_x^m u\|^{q-2}_{L^2}\int_D u \left|\partial^m_x u\right|^2\,dx &\leq C\|u\|_{H^{m-1}}\|\partial_x^m u\|^q_{L^2} \\
				&\leq C \|u\|^2_{H^{m-1}} \|\partial_x^m u\|^{q-2}_{L^2}\|\partial_x^{m+1}u\|_{L^2} \\
				&\leq \frac{\nu q}{4}\|\partial^m_x u\|^{q-2}_{L^2} \|\partial_x^{m+1} u\|^2_{L^2} + C\|\partial_x^m u\|^q_{L^2} + C\|u\|^{2q}_{H^{m-1}}.
			\end{align*}
			
			As the same argument in \eqref{eq:estimate01} and \eqref{eq:estimate-02}, we have
			\begin{align*}
				\beta(1+\gamma) q \|\partial_x^m u\|^{q-2}_{L^2}\big(\partial^m_x u, \partial^m_x (u^2) - 2u\partial^m_x u\big) \leq&\, C\|u\|_{H^{m-1}}^2 \|\partial_x^m u\|_{L^2}^{q-1} \\
				\leq&\, C\|\partial^m_x u\|^q_{L^2} + C\|u\|^{2q}_{H^{m-1}}.
			\end{align*}
			
			Collecting the above estimates and substituting them into \eqref{eqn:partial-m-u}, we obtain
			\begin{equation}\label{eq:estimate-m}
				\begin{aligned}
					d&\|\partial_x^m u\|_{L^2}^q + \frac{\nu q}{2}\|\partial_x^m u\|_{L^2}^{q-2} \|\partial_x^{m+1} u\|_{L^2}^2 \, dt \\
					&\leq q \|\partial_x^m u\|_{L^2}^{q-2} \bigl(\partial_x^m u, \partial_x^m \sigma(u)\bigr) \, dW(t) + C\|\partial_x^m u\|_{L^2}^q \,dt + C\|u\|_{H^{m-1}}^{3q} \,dt  + C.
				\end{aligned}
			\end{equation}
			
			Integrating over $[0,t]$ and taking expectations, we obtain, by Fubini's theorem,
			\begin{align*}
				&\mathbb{E}\big[\|\partial_x^m u(t)\|_{L^2}^q\big] + \frac{\nu q}{2}\int^t_0\mathbb{E}\big[\|\partial_x^m u(s)\|_{L^2}^{q-2} \|\partial_x^{m+1} u(s)\|_{L^2}^2\big] \, ds \\
				&\leq \mathbb{E}\big[\|\partial^m_x u_0\|^q_{L^2}\big] + C\int^t_0\Big(1 + \mathbb{E}\big[\|\partial_x^m u(s)\|_{L^2}^q\big] + \mathbb{E}\big[\|u(s)\|_{H^{m-1}}^{3q} \big]\Big)\,ds,
			\end{align*}
			where we used that the stochastic integral
			\begin{align*}
				M_m(t) := q \int^t_0\|\partial_x^m u\|_{L^2}^{q-2} \bigl(\partial_x^m u, \partial_x^m \sigma(u)\bigr) \, dW(s)
			\end{align*}
			is a martingale with zero expectation. By the induction hypothesis, we yield
			\begin{align*}
				\int^t_0 \mathbb{E}\big[\|u(s)\|_{H^{m-1}}^{3q}\big]\,ds \leq C\Big(\sum^{m-1}_{i=0} \mathbb{E}\big[\|u_0\|^{3^{m-i}q}_{H^i}\big] + 1\Big),
			\end{align*}
			implying
			\[\begin{aligned}
				\mathbb{E}\left[\|\partial_x^m u(t)\|_{L^2}^{q}\right] &+ \frac{\nu q}{2} \int_0^t \mathbb{E}\left[ \|\partial_x^m u(s)\|_{L^2}^{q-2} \|\partial_x^{m+1}u(s)\|_{L^2}^{2} \right]\,ds \\
				&\leq \mathbb{E}\left[\|\partial_x^m u_0\|_{L^2}^{q}\right] + C\int_0^t \mathbb{E}\left[ \|\partial_x^m u(s)\|_{L^2}^{q} \right]\,ds + C\sum_{i=0}^{m-1} \mathbb{E}\left[\|\partial_x^i u_0\|_{L^2}^{3^{m-i}q}\right]+C \\
				&\leq C\int_0^t \mathbb{E}\left[ \|\partial_x^m u(s)\|_{L^2}^{q} \right]\,ds + C\sum_{i=0}^{m} \mathbb{E}\left[ \|\partial_x^i u_0\|_{L^2}^{3^{m-i}q}\right] + C .
			\end{aligned}\]
			
			As a consequence of Gronwall's inequality, we obtain
			\begin{equation}\label{eq:estimate-03}
				\mathbb{E}\left[\|\partial_x^m u(t)\|_{L^2}^{q}\right] \leq Ce^{Ct} \left(\sum_{i=0}^{m}\mathbb{E}\left[\|\partial_x^i u_0\|_{L^2}^{3^{m-i}q}\right] + 1\right).
			\end{equation}
			
			This yields \eqref{eq:Hm-moment-estimate} for the general case \(m\geq 1\) and for all \(q\geq 2\). From \eqref{eq:estimate-m}, we see that
			\begin{equation}\label{eq:u-sup-moment-pre-bdg-c}
				\begin{aligned}
					\mathbb{E}\!\left[\sup_{s \in [0,T]} \|\partial_x^mu(s)\|_{L^2}^q \right] &\le \mathbb{E}\big[\|\partial_x^m u_0\|_{L^2}^q\big] + \mathbb{E}\!\left[\sup_{s \in [0,T]} |M_m(s)| \right] \\
					&\quad+ C \int_0^T \bigl(1 + \mathbb{E}\big[\|\partial_x^mu(s)\|_{L^2}^q\big] + \mathbb{E}\big[\| u(s)\|_{H^{m-1}}^{3q}\big] \bigr)\, ds.
				\end{aligned}
			\end{equation}
			
			We estimate the quadratic variation of $M_m$ as in the case $m = 0$
			\[
			d\langle M_m \rangle_s = q^2 \|\partial_x^mu(s)\|_{L^2}^{2q-4} \left|(\partial_x^m u(s), \partial_x^m\sigma(u(s)))\right|^2 ds \leq C\|\partial_x^m u(s)\|_{L^2}^{2q-2} \big(1 + \|\partial_x^m u(s)\|^2_{L^2}\big)\, ds.
			\]
			
			Using the same argument as in \eqref{eq:M0-BDG-estimate}, we yield
			\begin{align*}
				\mathbb{E}\!\left[\sup_{s \in [0,T]} |M_m(s)| \right] &\leq \frac{1}{10} \,\mathbb{E}\! \left[ \sup_{s \in [0,T]} \|\partial_x^m u(s)\|_{L^2}^{q}\right] + C \int_0^T \bigl(1 + \mathbb{E}\big[\|\partial_x^mu(s)\|_{L^2}^q\big]\bigr)\,ds.
			\end{align*}
			
			Applying this estimate to \eqref{eq:u-sup-moment-pre-bdg-c}, we obtain
			\begin{equation*}
				\begin{aligned}
					\frac{9}{10}\mathbb{E}\!\left[\sup_{s \in [0,T]} \|\partial_x^mu(s)\|_{L^2}^q \right] &\le \mathbb{E}\big[\|\partial^m_x u_0\|^q_{L^2}\big] + C \int_0^T \bigl(1 + \mathbb{E}\big[\|\partial_x^mu(s)\|_{L^2}^q\big] + \mathbb{E}\big[\|u(s)\|_{H^{m-1}}^{3q}\big]\bigr)\,ds.
				\end{aligned}
			\end{equation*}
			
			Using \eqref{eq:estimate-03}, we obtain
			\[
			\int_0^T\mathbb{E}\big[\|\partial_x^m u(s)\|_{L^2}^q\big] \,ds \le C\left(\sum_{i=0}^m\mathbb{E}\left[\|\partial_x^i u_0\|_{L^2}^{3^{m-i}q}\right]+1\right).
			\]
			
			Moreover, by the induction hypothesis applied to the order \(m-1\) with exponent \(3q\), we obtain
			\[ 
			\int_0^T\mathbb{E}\left[\|u(s)\|_{H^{m-1}}^{3q}\right]\,ds \le C\left(\sum_{i=0}^{m-1}\mathbb{E}\left[ \|\partial_x^i u_0\|_{L^2}^{3^{m-i}q}\right]+1\right).
			\]
			
			Consequently,
			\[
			\mathbb{E}\left[\sup_{s\in[0,T]} \|\partial_x^m u(s)\|_{L^2}^q\right]\le C\left(\sum_{i=0}^{m}\mathbb{E}\left[\|\partial_x^i u_0\|_{L^2}^{3^{m-i}q}\right]+1\right).
			\]
			
			This establishes \eqref{eq:Hm-sup-moment-estimate} for the general case \(m\ge 1\) and for all \(q\ge 2\), and hence completes the proof.
		\end{proof}
		
		\subsection{Proof of Lemma \ref{lem:time-regularity}}\label{AppendixB3}
		\begin{proof}
			By Theorem~\ref{theo:unique-sol}, there exists a unique variational solution $u$. Moreover, $u$ admits the following mild representation; see \cite[Theorem~6.5]{gawarecki2010stochastic}:
			\begin{equation}\label{eq:mild-solution}
				u(t)=S(t)u_0+\int_0^tS(t-r)N(u(r))\,dr+\int_0^tS(t-r)\sigma(u(r))\,dW(r),
			\end{equation}
			where $N(u):=-\alpha u\partial_xu+\beta c(u)$. 
			
			Using \eqref{eq:mild-solution} and the semigroup property, we obtain, for $0\leq s\leq t\leq T$,
			\[
			u(t)=S(t-s)u(s)+\int_s^tS(t-r)N(u(r))\,dr+\int_s^tS(t-r)\sigma(u(r))\,dW(r).
			\]
			
			Consequently,
			\[
			\begin{aligned}
				u(t)-u(s)&=(S(t-s)-\mathcal{I})u(s)+\int_s^tS(t-r)N(u(r))\,dr+\int_s^tS(t-r)\sigma(u(r))\,dW(r)\\
				&=:I_1+I_2+I_3.
			\end{aligned}
			\]
			
			Hence, $\mathbb E\Big[\|u(t)-u(s)\|_{H^m}^{2q}\Big] \leq C\sum_{i=1}^3\mathbb E\Big[\|I_i\|_{H^m}^{2q}\Big]$. We estimate the three terms separately. First, by the semigroup estimate
			\begin{align*}
				\|I_1\|_{H^m} &= \|(S(t-s)-\mathcal{I})u(s)\|_{H^m} \leq C\|\Lambda^{m}(S(t-s) - \mathcal{I})u(s)\|_{L^2} = C\|\Lambda^{-1}(\mathcal{I} - S(t-s))\Lambda^{m+1}u(s)\|_{L^2}.
			\end{align*}
			
			By \eqref{eq:semigroup-increment} with \(\rho=1/2\), it follows that
			\begin{align*}
				\|\Lambda^{-1}(\mathcal{I} - S(t-s))\Lambda^{m+1}u(s)\|_{L^2} &\leq C(t-s)^{1/2}\|\Lambda^{m+1}u(s)\|_{L^2} \leq C(t-s)^{1/2}\|u(s)\|_{H^{m+1}}.
			\end{align*}
			
			By Lemma \ref{lem:Hm-q-moment}, we yield
			\[
			\begin{aligned}
				\mathbb E\Big[\|I_1\|_{H^m}^{2q}\Big] &= \mathbb E\Big[\|(S(t-s)-\mathcal{I})u(s)\|_{H^m}^{2q}\Big] \leq C (t-s)^q \mathbb E\Big[\|u(s)\|_{H^{m+1}}^{2q}\Big] \leq C (t-s)^q.
			\end{aligned}
			\]
			
			For $I_2$, we observe that \(S(t)\) is bounded on \(H^m\), H\"older's inequality gives
			\[
			\begin{aligned}
				\mathbb E\Big[\|I_2\|_{H^m}^{2q}\Big] &= \mathbb{E} \Bigg[\left\|\int_s^t S(t-r)N(u(r))\,dr \right\|^{2q}_{H^m}\Bigg] \leq C (t-s)^{2q-1} \int_s^t \mathbb E\Big[\|N(u(r))\|_{H^m}^{2q}\Big]\,dr.
			\end{aligned}
			\]
			
			Moreover, by \eqref{eq:cubic-estimate} and \eqref{eq:burger-estimate}, we have
			\begin{align*}
				\|N(u)\|_{H^m}^{2q} &= \|-\alpha u\partial_xu+\beta c(u)\|_{H^m}^{2q} \\
				&\leq C\left(\|u\partial_x u\|^{2q}_{H^m} + \|c(u)\|^{2q}_{H^m} \right)\\
				&\leq C\left(1 + \|u\|^{4q}_{H^{m+1}} + \|u\|^{6q}_{H^{m+1}} \right)\\
				&\leq C\left(1+\|u\|_{H^{m+1}}^{6q}\right).
			\end{align*}
			
			By Lemma \ref{lem:Hm-q-moment}, we derive
			\[
			\begin{aligned}
				\mathbb E\Big[\|I_2\|_{H^m}^{2q}\Big] 
				&\leq C (t-s)^{2q-1} \int_s^t \mathbb E\left[ 1+\|u(r)\|_{H^{m+1}}^{6q} \right]dr \\
				&\le C(t-s)^{2q-1} \int_s^t \Big(1+\sup_{\tau\in[0,T]}\mathbb E\Big[\|u(\tau)\|_{H^{m+1}}^{6q}\Big] \Big)\,dr\\
				&\leq C (t-s)^{2q}.
			\end{aligned}
			\]
			
			Finally, to estimate \(I_3\), we apply the Burkh\"older--Davis--Gundy inequality and use the boundedness of \(S(t)\) on \(H^m\),
			\begin{align*} 
				\mathbb E\big[\|I_3\|_{H^m}^{2q}\big] &= \mathbb{E}  \Bigg[\left\|\int_s^t S(t-r)\sigma(u(r))\,dW(r)\right\|^{2q}_{H^m}\Bigg]\\
				&\leq C\mathbb E \left[\left(\int_s^t\|S(t-r)\sigma(u(r))\|_{H^m}^2\,dr \right)^q\right] \\
				&\leq C\mathbb E\left[\left(\int_s^t\|\sigma(u(r))\|_{H^m}^2\,dr\right)^q\right] \\
				&\leq C (t-s)^{q-1} \int_s^t \mathbb E\big[\|\sigma(u(r))\|_{H^m}^{2q}\big]\,dr.
			\end{align*}
			
			By the linear growth assumption \eqref{ass:noise-linear-growth}, we have
			\[
			\begin{aligned}
				\mathbb E\Big[\|I_3\|_{H^m}^{2q}\Big] &\leq C (t-s)^{q-1} \int_s^t \mathbb E\left[ 1+\|u(r)\|_{H^{m}}^{2q} \right]dr \\
				&\le C(t-s)^{q-1} \int_s^t \Big(1+\sup_{\tau\in[0,T]}\mathbb E\Big[\|u(\tau)\|_{H^{m}}^{2q}\Big] \Big)\,dr\\
				&\leq C (t-s)^{q}.
			\end{aligned}
			\]
			
			Combining the estimates for \(I_1\), \(I_2\), and \(I_3\), we conclude that
			\[
			\mathbb E\big[\|u(t)-u(s)\|_{H^m}^{2q}\big] \leq C(t-s)^q.
			\]
			
			This completes the proof.
		\end{proof}

		\bibliographystyle{abbrv}
		\bibliography{references_SBH}	
		
	\end{document}